\documentclass[a4,10pt]{article}
\usepackage{amsmath,amscd,amssymb}

\newcommand{\nbige}{\mathcal{E}}

\newcommand{\nbigh}{\mathcal{H}}

\newcommand{\nbigl}{\mathcal{L}}

\newcommand{\nbigo}{\mathcal{O}}
\newcommand{\nbigp}{\mathcal{P}}

\newcommand{\nbigs}{\mathcal{S}}

\newcommand{\proj}{\mathbb{P}}
\newcommand{\seisuu}{\mathbb{Z}}
\newcommand{\rnum}{{\boldsymbol Q}}
\newcommand{\nnum}{{\boldsymbol N}}
\newcommand{\cnum}{{\boldsymbol C}}
\newcommand{\real}{{\boldsymbol R}}
\newcommand{\hyperh}{\mathbb{H}}

\newcommand{\DD}{\mathbb{D}}
\newcommand{\EE}{\mathbb{E}}

\newcommand{\gminie}{\mathfrak e}

\newcommand{\gminim}{\mathfrak m}

\newcommand{\gminip}{\mathfrak p}

\newcommand{\gminit}{\mathfrak t}

\newcommand{\vece}{{\boldsymbol e}}

\newcommand{\vecv}{{\boldsymbol v}}

\newcommand{\vecalpha}{{\boldsymbol \alpha}}

\newcommand{\vecn}{{\boldsymbol n}}

\newcommand{\lrarr}{\longrightarrow}

\newcommand{\pf}{{\bf Proof}\hspace{.1in}}
\newcommand{\qed}{\mbox{\rule{1.2mm}{3mm}}}

\def\End{\mathop{\rm End}\nolimits}

\def\Image{\mathop{\rm Im}\nolimits}

\def\Re{\mathop{\rm Re}\nolimits}

\def\Gr{\mathop{\rm Gr}\nolimits}

\def\rank{\mathop{\rm rank}\nolimits}

\def\Ker{\mathop{\rm Ker}\nolimits}

\def\Gr{\mathop{\rm Gr}\nolimits}

\def\Res{\mathop{\rm Res}\nolimits}

\def\tr{\mathop{\rm tr}\nolimits}
\def\vol{\mathop{\rm dvol}\nolimits}
\def\dvol{\mathop{\rm dvol}\nolimits}

\def\id{\mathop{\rm id}\nolimits}

\newcommand{\del}{\partial}
\newcommand{\delbar}{\overline{\del}}

\newcommand{\Deltabar}{\overline{\Delta}}
\newcommand{\Deltabarast}{\Deltabar^{\ast}}

\newcommand{\nbar}{\underline{n}}

\newcommand{\harmonicbundle}{(E,\delbar_E,\theta,h)}

\newcommand{\prolong}[1]{{}^{\diamond}{#1}}

\newcommand{\KMS}{{\mathcal{KMS}}}
\newcommand{\KMSoverline}{\overline{\mathcal{KMS}}}

\newcommand{\Sp}{{\mathcal Sp}}

\newcommand{\paramap}{\gminip}
\newcommand{\eigenmap}{\gminie}
\newcommand{\multiplicity}{\gminim}

\newcommand{\lefttop}[1]{{}^{#1}}

\newcommand{\closedopen}[2]{[#1,#2[}
\newcommand{\openclosed}[2]{]#1,#2]}

\newcommand{\PH}{{\mathcal P}{\mathcal H}}
\newcommand{\bdmath}{\begin{displaymath}}
\newcommand{\edmath}{\end{displaymath}}

\newcommand{\beqn}{\begin{equation}}
\newcommand{\eeqn}{\end{equation}}

\newcommand{\beqnarray}{\begin{eqnarray}}
\newcommand{\eeqnarray}{\end{eqnarray}}

\newcommand{\bitemize}{\begin{itemize}}
\newcommand{\eitemize}{\end{itemize}}

\newcommand{\benumerate}{\begin{enumerate}}
\newcommand{\eenumerate}{\end{enumerate}}

\newcommand{\bdescriprion}{\begin{description}}
\newcommand{\edescriprion}{\end{description}}
\newtheorem{thm}{Theorem}[section]
\newtheorem{cor}{Corollary}[section]

\newtheorem{rem}{Remark}[section]
\newtheorem{lem}{Lemma}[section]
\newtheorem{prop}{Proposition}[section]
\newtheorem{df}{Definition}[section]

\newtheorem{condition}{Condition}[section]

\def\longuparrow{\vcenter{%
   \lineskip0pt\lineskiplimit0pt\baselineskip0pt\ialign{%
   \hfil{##}\hfil\crcr\hbox to 0pt{\hss$\uparrow$\hss}\cr%
   \hbox to 0pt{\hss\vrule width.4pt depth 0pt height 1em\hss}\cr}}}
\def\longdownarrow{\vcenter{%
   \lineskip0pt\lineskiplimit0pt\baselineskip0pt\ialign{\hfil{##}\hfil%
   \crcr\hbox to 0pt{\hss\vrule width.4pt depth 0pt height 1em\hss}\cr%
   \hbox to 0pt{\hss$\downarrow$\hss}\cr}}}

\begin{document}

\title{A characterization of semisimple local system\\
 by tame pure imaginary pluri-harmonic metric }
\author{Takuro Mochizuki}
\date{}

\maketitle

\begin{abstract}
 
Let $L$ be a local system on a smooth quasi projective variety
over $\cnum$.
We see that $L$ is semisimple if and only if
there exists a tame pure imaginary pluri-harmonic metric on $L$.
Although it is a rather minor refinement of a result of Jost and Zuo,
it is significant for the study of
harmonic bundles and pure twistor $D$-modules.
As one of the application, we show that
the semisimplicity of local systems are preserved
by the pull back via a morphism of quasi projective varieties.

This paper will be included as an appendix
to our previous paper,
``{\it  Asymptotic behaviour of tame harmonic bundles and
an application to pure twistor $D$-modules}'',
math.DG/0312230. \\


\end{abstract}

\section{Introduction}

\subsection{Main result}

%c1.tex

\subsubsection{Main theorem}

Let $X$ be a smooth projective variety over $\cnum$,
and $D$ be a normal crossing divisor of $X$.
Let $(E,\nabla)$ be a flat bundle on $X-D$.
It is our main purpose to show the following theorem,
which gives a characterization of semisimplicity of flat bundles
by the existence of a pure imaginary pluri-harmonic metric.

\begin{thm}  \mbox{{}}\label{thm;04.1.4.1}
\begin{itemize}
\item
The flat bundle $(E,\nabla)$ is semisimple
if and only if there exists a tame pure imaginary pluri-harmonic metric
$h$ on $(E,\nabla)$. 
\item
If $(E,\nabla)$ is simple,
then the tame pure imaginary pluri-harmonic metric $h$ is 
uniquely determined up to positive constant multiplication.
\end{itemize}
\end{thm}

Let us explain a tame pure imaginary pluri-harmonic metric.
From a pluri-harmonic metric $h$ on $(E,\nabla)$,
we obtain the harmonic bundle $\harmonicbundle$.
Roughly speaking,
$\harmonicbundle$ is tame
if there exists a holomorphic bundle $\tilde{E}$
and a regular Higgs field
$\tilde{\theta}\in End(\tilde{E})\otimes\Omega^{1,0}_X(\log D)$
such that $(\tilde{E},\tilde{\theta})_{|X-D}=(E,\theta)$.
If any eigenvalues of the residues of $\tilde{\theta}$ are
pure imaginary,
then $\harmonicbundle$ is called pure imaginary.
We remark that the eigenvalues of the residues
are independent of a choice of 
a prolongment $(\tilde{E},\tilde{\theta})$.

\subsubsection{Some remarks}

It can be said that 
Theorem \ref{thm;04.1.4.1} is a partial refinement
of the result of Jost-Zuo in \cite{JZ2},
which is technically minor but significant for our application.
Let us explain for more detail.

In \cite{JZ2}, Jost and Zuo discussed the existence
of tame pluri-harmonic (twisted) maps
from a complement of a normal crossing divisor in 
a compact Kahler manifold
to the symmetric space of non-compact types.
On the other hand, we only consider $GL(r)/U(r)$
and a quasi projective variety 
as the target space and the domain respectively.
However we can impose the pure imaginary condition
to the behaviour of the tame pluri-harmonic twisted map
at infinity.
We can also derive the `only if' part.
It is the meaning of `partial refinement'.

The refinement itself is rather technically minor
in the following sense.
We can observe `the pure imaginary condition' easily from
the argument in \cite{JZ2}.
The `only if' part is a rather easy consequence of
the observation of Sabbah given in \cite{sabbah}.

However, it seems significant for the theory of pure twistor $D$-modules.
Briefly speaking, Theorem \ref{thm;04.1.4.1} gives
a characterization of semisimple local system on a quasi projective
variety.
From Theorem \ref{thm;04.1.4.1},
we can derive the correspondence of
semisimple perverse sheaves and
a `pure imaginary' pure twistor $D$-modules,
although the latter has not appeared in the literature.
By using it together with Sabbah's theory \cite{sabbah}
and our result in \cite{mochi2},
it seems possible to show the regular holonomic version
of Kashiwara's conjecture
(see Introduction in \cite{mochi2}).
We will discuss it elsewhere. 
In this paper, we give only the following theorem
as an easy application,
which is an affirmative answer to a question posed by Kashiwara.

\begin{thm}[Theorem \ref{thm;04.1.30.100}]
Let $X$ and $Y$ be irreducible quasi projective varieties over $\cnum$.
Let $F:X\lrarr Y$ be a morphism.
Let $L$ be a semisimple local system on $Y$.
Then the pull back $F^{-1}(L)$ is also semisimple.
\hfill\qed
\end{thm}

As is noted above,
our main result (Theorem \ref{thm;04.1.4.1}) can be regarded
as a technically minor and partial refinement of
the result of \cite{JZ2},
once we completely understand the proof of the existence theorem
of pluri-harmonic metric.
Hence the author should explain why he writes
this paper, which looks rather long.
He feels that the argument of Jost-Zuo 
seems the argument for the experts of harmonic maps,
and that 
it may not seem so easy to understand for non-specialists, at a sight.
However, Theorem \ref{thm;04.1.4.1} is one of the most important
key steps in our application of harmonic bundle
to pure twistor $D$-modules.
Hence the author thinks it appropriate
to give a detailed proof of Theorem \ref{thm;04.1.4.1},
which is available for a wide range of the readers.
We will start from elementary facts,
and we give the proofs of some rather well known results
when the author does not know an appropriate reference.
We give the details of the estimates,
because it is rather delicate to deal with the infinite energy.
We have to control the divergent term carefully.
It is one of the reason why the paper is rather long.

As is explained,
the most part of the paper is an effect of our effort
to understand \cite{JZ2},
and the most essential ideas
for the existence part of Theorem \ref{thm;04.1.4.1}
are due to Jost and Zuo,
although we do not follow their arguments straightforwardly.
Needless to say,
the author is responsible for any mistakes contained in this paper.

\subsection{The outline of the paper}

%c1.2.tex

\subsubsection{Section \ref{section;04.1.30.200}}

In the subsection \ref{subsection;04.1.30.201},
we recall some standard facts
just for our reference of the later discussion.
In the subsection \ref{subsection;04.1.30.202},
we recall the elementary geometry of the symmetric space $\PH(r)$
of the positive definite hermitian metrics.
Lemma \ref{lem;a12.29.10} is one of the key lemmas,
although it is elementary.
We also give the comparison of 
the distance of the hermitian metrics in $\PH(r)$
and the norm of the identity with respect to the two metrics
in the subsubsection \ref{subsubsection;04.1.30.300}.

In the subsection \ref{subsection;04.1.30.203},
we discuss a twisted map associated with
a commuting tuple of endomorphisms.
The result will be used in the subsection
\ref{subsection;04.1.1.2}
to construct the twisted map of $(\Deltabarast)^2$
whose energy is controlled.

In the subsection \ref{subsection;04.1.30.204},
we recall some standard facts on twisted harmonic maps.
The Bochner type formula
in the subsubsection \ref{subsubsection;04.1.8.30}
is due to Corlette.
A variation of Bochner type formula is given
in the subsubsection \ref{subsubsection;04.1.30.350}.
They will be used in the proof of the pluri-harmonicity
(see the subsections
\ref{subsection;04.1.30.400}--\ref{subsection;04.1.30.401}).
It is important for our argument to consider two kinds of
Bochner type formula.

\subsubsection{Section \ref{section;04.1.30.205}}

In the section \ref{section;04.1.30.205},
we give the definition of tame pure imaginary harmonic bundles,
and some of useful properties.
In the subsection \ref{subsection;04.1.30.206},
we give a definition of pure imaginary property
of tame harmonic bundles.
In the subsection \ref{subsection;04.1.30.207},
we see the estimate of the energy functions of 
a tame pure imaginary harmonic bundles
on a punctured disc.
We give a characterization of tame and pure imaginary properties
by an increasing order of the energy
in the subsubsection \ref{subsubsection;04.1.29.300}.

In the subsection \ref{subsection;04.1.30.208},
we show that the underlying flat bundle
of a tame pure imaginary harmonic bundle is semisimple.
It is a consequence of the observation,
which is essentially due to Sabbah.

In the subsection \ref{subsection;04.1.30.209},
we see the maximum principle for
the distance of two tame pure imaginary harmonic metrics
on a punctured disc.
The result will be useful to control the energy
(the section \ref{section;04.1.30.410}).
It will be also used to show the tameness
(the subsection \ref{subsection;04.1.30.411}).

In the subsection \ref{subsection;04.1.30.210},
we show the uniqueness of tame pure imaginary harmonic bundle
of a flat bundle on a quasi projective variety.
Note that it essentially follows from
the Kobayashi-Hitchin correspondence
of Simpson and Biquard (\cite{s2} and \cite{b}).
However the detailed proof for uniqueness for
in the case of parabolic flat bundle
seems to be omitted there.
Thus we give the detailed proof within our necessity.
We essentially follow the argument of Corlette.

\subsubsection{Section \ref{section;04.1.30.211}}

In the subsection \ref{subsection;04.1.6.10},
we discuss the Dirichlet problem of 
a tame pure imaginary harmonic bundle on a punctured disc.
The argument seems essentially due to Lohkamp \cite{lohkamp}
and Jost-Zuo \cite{JZ2}.

In the subsection \ref{subsection;04.1.30.212},
we discuss the family version of the Dirichlet problem.
We give the estimate of the differentials,
which is given by
the maximum principle (the subsection \ref{subsection;04.1.30.209})
and the estimate in the subsection \ref{subsection;04.1.30.207}.
The idea is essentially due to Jost and Zuo (\cite{JZ1}).
The result will be used for the construction
of the twisted map whose energy is controlled
(the section \ref{section;04.1.30.410}).

\subsubsection{Section \ref{section;04.1.30.410}}

We construct the twisted map 
on the complement of a normal crossing divisor
in a compact Kahler manifold.
We essentially follow the method of Jost-Zuo (\cite{JZ1} and \cite{JZ2}).

In the subsection \ref{subsection;04.1.8.20},
we construct the twisted map
around smooth points of a divisor,
by solving the family of the Dirichlet problem.
Then the energy of the map is controlled
by the result of the subsection \ref{subsection;04.1.30.212}.
We also give the lower bound of the energy
for arbitrary twisted map around smooth points of a divisor
in the subsubsection \ref{subsubsection;04.1.15.1}.
Essentially it is a consequence of
Lemma \ref{lem;a12.29.10}.
However we need some care to control the divergent term precisely.

In the subsection \ref{subsection;04.1.1.2},
we construct the twisted map around the intersection
of the divisors.
The results in the subsection \ref{subsection;04.1.30.203}
and the subsubsection \ref{subsubsection;04.1.7.20}
are used.
We also give the lower bound of the energy
of arbitrary twisted maps.
Again, we have to be careful to control the divergent term.

In the subsection \ref{subsection;04.1.11.5},
we give the decomposition of the complement of a normal crossing divisor
in a compact Kahler manifold,
and we obtain the twisted map whose energy is controlled.
We also give the lower bound of the energy
of arbitrary twisted map.
They are direct consequences of the results
in the subsections
\ref{subsection;04.1.8.20}--\ref{subsection;04.1.1.2}.

\subsubsection{Section \ref{section;04.1.30.220}}

In the subsection \ref{subsection;04.1.30.221},
we obtain the harmonic metric
of a semisimple flat bundle on a quasi projective variety.
The argument is essentially same
as that in the Dirichlet problem on a punctured disc
(the subsection \ref{subsection;04.1.6.10}),
except for the use of the argument in \cite{jy4}.
In the subsubsections
\ref{subsubsection;04.1.15.15}--\ref{subsubsection;04.1.23.1},
we give the detailed estimate of the energy of the resulted harmonic metric.
They are necessary for the later discussion.

In the subsection \ref{subsection;04.1.30.400},
we show that $\delbar\theta$ and $\theta^2$ are $L^2$,
where $\theta$ is the associated $(1,0)$-form
for the resulted harmonic metric.
We use the Bochner type formula
in the subsubsection \ref{subsubsection;04.1.8.30}.
If the integral of the Bochner type formula
would vanish, then we would obtain the pluri-harmonicity.
However it does not seem easy to show such vanishing directly.
(See the convergence (\ref{eq;04.1.30.500}), for example.)

Hence, in the subsection \ref{subsection;04.1.30.401},
we use another kind of Bochner type formula
given in the subsubsection \ref{subsubsection;04.1.30.350}.
It is rather easy to show the vanishing of the integral in this time,
by using the $L^2$-property of $\delbar\theta$
obtained in the subsection \ref{subsection;04.1.30.400}.
As a result, we obtain the pluri-harmonicity
of the harmonic metric obtained in the subsection
\ref{subsection;04.1.30.221}.

In the subsection \ref{subsection;04.1.30.411},
we show the tameness and the pure imaginary property
for the resulted pluri-harmonic metric.
It is rather easy consequence of
the estimate of the energy
given in the subsubsection \ref{subsubsection;04.1.23.1},
the characterization of tame pure imaginary harmonic bundle
on a punctured disc
(the subsubsection \ref{subsubsection;04.1.29.300}),
the maximum principle (the subsection \ref{subsection;04.1.30.210})
and Hartogs type theorem.

Thus we obtain the existence theorem
of a tame and pure imaginary pluri-harmonic harmonic metric
for any semisimple flat bundle on a quasi projective surface,
or more generally, the complement of a normal crossing divisor
in a compact Kahler surface.

In the subsection \ref{subsection;04.1.30.412},
we show the existence of tame and pure imaginary pluri-harmonic
metric in the higher dimensional case,
by reducing the problem
to the case of quasi projective surface.
Here we need the quasi projectivity.
To generalize the argument
in the sections \ref{section;04.1.30.410}--\ref{section;04.1.30.220}
in higher dimensional case,
it seems that we need some additional arguments.
For example, naively speaking,
we need the family version of the argument
in the subsection \ref{subsection;04.1.1.2}.
However, perhaps, the author feels that
it seems not so straightforward,
for we cannot use the maximum principle for 
the family version of the map $F$ on $Y$,
for example.

\subsubsection{Section \ref{section;04.1.30.230}}

As a simple application,
we show Theorem \ref{thm;04.1.30.100}.
We have only to show that
the pull back of a tame pure imaginary harmonic bundle
is also a tame pure imaginary harmonic bundle.

\subsection{Acknowledgement}

%c1.1.tex

The author is grateful to the colleagues in Osaka City University.
He specially thanks Mikiya Masuda for his encouragement and supports.
The author thanks Yoshifumi Tsuchimoto and Akira Ishii
for their constant encouragements.
The author thanks Masaki Kashiwara for his attractive question.
The author thanks the financial support by Japan Society for
the Promotion of Science.

\section{Preliminary}

\label{section;04.1.30.200}

\subsection{Notation}

%c10.tex

\subsubsection{Sets}

We will use the following notation:

\begin{tabular}{llll}
$\seisuu$: & the set of the integers, &
$\seisuu_{>0}$: &the set of the positive integers,\\
$\rnum$: &the set of the rational numbers,&
$\rnum_{>0}$: &the set of the positive rational numbers,\\
$\real$: & the set of the real numbers, &
$\real_{>0}$: &the set of the positive real numbers,\\
$\cnum$: &the set of the complex numbers, &
$\nbar$: & the set $\{1,2,\ldots,n\}$,\\
$M(r)$: & the set of $r\times r$-matrices,&
$\nbigh_r$: & the set of $r\times r$-hermitian matrices,\\

\end{tabular}\\

We denote the set of positive hermitian metric
of $V$ by $\nbigp\nbigh(V)$.
We often identify it with the set of the positive hermitian matrices
by taking an appropriate base of $V$.

We put $[a,b]:=\{x\in\real\,|\,a\leq x\leq b\}$,
$\closedopen{a}{b}:=\{x\in\real\,|\,a\leq x<b\}$,
$\openclosed{a}{b}:=\{x\in\real\,|\,a<x\leq b\}$
for any $a,b\in\real$.

\subsubsection{A disc, a punctured disc and some products}

For any positive number $C>0$ and $z_0\in\cnum$,
the open disc $\bigl\{z\in\cnum\,\big|\,|z-z_0|<C\bigr\}$
is denoted by $\Delta(z_0,C)$,
and the punctured disc $\Delta(z_0,C)-\{z_0\}$
is denoted by $\Delta^{\ast}(z_0,C)$.
When $z_0=0$, $\Delta(0,C)$ and $\Delta^{\ast}(0,C)$
are often denoted by $\Delta(C)$ and $\Delta^{\ast}(C)$.
Moreover, if $C=1$, $\Delta(1)$ and $\Delta^{\ast}(1)$
are often denoted by $\Delta$ and $\Delta^{\ast}$.
If we emphasize the variable,
we describe as $\Delta_z$, $\Delta_i$.
For example,
$\Delta_z\times \Delta_w=\{(z,w)\in\Delta\times \Delta\}$,
and $\Delta_1\times \Delta_2=\{(z_1,z_2)\in\Delta\times \Delta\}$.
We often use the notation $\cnum_{\lambda}$ and $\cnum_{\mu}$
to denote the complex planes
$\bigl\{\lambda\in\cnum\bigr\}$
and $\bigl\{\mu\in\cnum\bigr\}$.

Unfortunately, the notation $\Delta$ is also used to denote the
Laplacian. The author hopes that there will be no confusion.

\subsection{Miscellaneous}
\label{subsection;04.1.30.201}

%c11.tex

\subsubsection{Differentiability of Lipschitz continuous functions}

We use the coordinate $(x,y_1,\ldots,y_l)$
for $\real\times\real^l$.
\begin{lem}
Let $f$ be a Lipschitz continuous function on $\real\times\real^l$.
Then $\frac{\del f}{\del x}$ is defined 
almost everywhere.
Namely we have the measurable function $F$
such that the following holds almost everywhere:
\[
 F(x,y)=\lim_{h\to 0}\frac{f(x+h,y)-f(x,y)}{h}.
\]
Moreover $\frac{\del f}{\del x}$ is bounded.
\end{lem}
\pf
In the case $l=0$,
the differentiability of an absolute continuous function,
and thus a Lipschitz continuous function,
is well known.
Let us consider the general case.
Let $h_i$ be any sequence of real numbers such that
$h_i\to 0$.
We put $F_i:=h_i^{-1}\cdot\bigl(f(x+h_i,y)-f(x,y)\bigr)$,
and then  we obtain the sequence of the measurable functions
$\{F_i\}$.
It is well known that $\overline{\lim}F_i$
and $\underline{\lim}F_i$ are measurable.
Thus the set
$S:=\bigl\{(x,y)\,\big|\,
 \overline{\lim}F_i(x,y)\neq \underline{\lim}F_i(x,y)
 \bigr\}$
is measurable.
By using the result in the case $l=0$,
we can easily derive that the measure of $S$ is $0$.
Hence we obtain the measurable function $F:=\lim F_i$.
It is easy to check that $F$ has the desired properties.
The boundedness follows from 
$|f(x+h,y)-f(x,y)|\leq C\cdot h$ for some constant $C$.
See any appropriate text book of measure theory
for the facts we used in the argument.
\hfill\qed

\vspace{.1in}
Let $f$ be a Lipschitz continuous function
on $\real\times\real^l$.
We obtain the measurable function $\del f/\del x$,
which is bounded.
It naturally gives the distribution.

On the other hand,
$f$ naturally gives the distribution.
Hence we obtain the differential of $f$
with respect to the variable $x$ as the distribution,
which we denote by $D_xf$.
\begin{lem}
We have $D_xf=\del f/\del x$ as the distribution.
\end{lem}
\pf
Let $\phi$ be a test function.
We have the following equality:
\begin{equation}\label{eq;a12.27.1}
 \int_{\real^{l+1}}
 \frac{f(x+h,y)-f(x,y)}{h}\cdot\phi(x,y)
=\int_{\real^{l+1}}
 f(x,y)\cdot\frac{\phi(x-h,y)-\phi(x,y)}{h}.
\end{equation}
Since $f$ is bounded, the right hand side of (\ref{eq;a12.27.1})
converges to $-\int_{\real^{l+1}}f\cdot (\del\phi/\del x)$,
due to the dominated convergence theorem.
Since $f$ is Lipschitz,
there exists a positive constant $C$
such that
$ \bigl|
 h^{-1}\cdot\bigl(f(x+h,y)-f(x,y)\bigr)
 \bigr|
\leq C$ holds for any $h$.
Hence the left hand side of (\ref{eq;a12.27.1})
converges to $\int_{\real^{l+1}}(\del f/\del x)\cdot \phi$.
Thus we are done.
\hfill\qed

\vspace{.1in}

\begin{cor}
Let $f$ be a Lipschitz function on $\real^l$.
Then $f$ is locally $L_1^p$ for any real number $p\geq 1$.
\hfill\qed
\end{cor}

%c11.1.tex

\subsubsection{Elementary linear algebra}

Let $V$ be an $r$-dimensional vector space,
and $h$ be a hermitian metric of $V$.
Let $\vecv$ be a base of $V$.
We put $H:=H(h,\vecv)$.
Let $f$ be an endomorphism of $V$,
and then we have the matrix $A$
such that $f\cdot\vecv=\vecv\cdot A$.
Let $f^{\dagger}$ denote the adjoint of $f$
with respect to the metric $h$.
Then we have
$f^{\dagger}\vecv=
 \vecv\cdot \overline{H}^{-1}\cdot\lefttop{t}\overline{A}\cdot\overline{H}$.

We have the norm $|f|_h$ of $f$ with respect to the metric $h$.
\begin{lem}\label{lem;04.1.6.1}
The following holds:
\[
 |f|_h^2
=\tr\bigl(A\cdot\overline{H}^{-1}\cdot
  \lefttop{t}\overline{A}\cdot \overline{H}\bigr).
\]
\end{lem}
\pf
We have $|f|_h^2=\tr(f\cdot f^{\dagger})$.
Then the claim immediately follows.
\hfill\qed

%c11.2.tex

\subsubsection{Hartogs type theorem}

\begin{lem} \label{lem;04.1.21.1}
Let $F$ be a holomorphic function on $\Deltabar^2-\{z_2=0\}$.
Assume that $F_{|\pi_2^{-1}(Q)}$ is holomorphic 
for almost every $Q\in \{z_2=0\}$.
Then $F$ is holomorphic on $\Delta$.
\end{lem}
\pf
We put as follows:
\[
 G(z_1,z_2):=\int_{|\zeta|=1}\frac{G(z_1,\zeta)}{(\zeta-z_2)}\cdot 
 \frac{d\zeta}{2\pi\sqrt{-1}}.
\]
Then $G(z_1,z_2)$ gives a holomorphic function on $\Deltabar^2$.
Due to the assumption,
there exists a dense subset $Y\subset\Deltabar^2-\{z_2=0\}$
such that
$G_{|Y}=F_{|Y}$.
Hence we obtain $G=F$ on $\Deltabar^2-\{z_2=0\}$.
\hfill\qed

\vspace{.1in}
The following lemma will be used later,
which can be shown similarly.
\begin{lem} \label{lem;04.1.21.2}
Let $F$ be a holomorphic function on $\Deltabar^2-\{z_2=0\}$.
Assume that there exists an open subset $U\subset \{z_2=0\}$
such that $F_{|\pi_2^{-1}(Q)}$ is holomorphic on $Q\times\Deltabar$.
Then $F$ is holomorphic on $\Deltabar^2$.
\hfill\qed
\end{lem}

%c11.3.tex

\subsubsection{Lefschetz hyperplane theorem for the fundamental group}

\begin{lem}
Let $X$ be a smooth projective surface,
and $D$ be a normal crossing divisor of $X$.
Let $H_1$ be a sufficiently ample smooth divisor
such that $H_1\cup D$ is normal crossing.
Then $\pi_1(H_1\setminus D)\lrarr\pi_1(X-D)$ is surjective.
\end{lem}
\pf
We give only a sketch of a proof.
Let $D_1$ be an ample smooth divisor of $X$
such that $D_0=D_1\cup D$ is very ample and normal crossing.
Let $D=D_2\cup\cdots \cup D_l$ be the irreducible decomposition.
Let $N_i$ denote the tubular neighbourhood of $D_i$
$(i=1,\ldots,l)$.
We put $N_0:=\bigcup_{i=1}^l N_i$.
Recall that we can decompose $X$ into $N_0$
and $i$-handles $(i=2,3,4)$
(See \cite{milnor}).
Hence the inclusion $N_0-D\lrarr X-D$ induces
the surjection $\pi_1(N_0-D)\lrarr \pi_1(X-D)$.

Let $s$ be a section of $\nbigo(D_0)$
satisfying that
the zero set $H_1:=s^{-1}(0)$ is smooth,
the $H_1\cap D_0$ is contained in the smooth part of $D_0$,
and $H_1\cup D$ is normal crossing.

We take the Lefschetz pencil $X'$
of $H_1$ and $D_0$.
We take the desingularization $\tilde{X}$ of $X'$.
We have the birational morphism
$p:\tilde{X}\lrarr X$
and the morphism $\pi:\tilde{X}\lrarr \proj^1$.
We may assume that
$\pi^{-1}(0)=p^{-1}(D_0)=p^{-1}(D)\cup p^{-1}(D_1)$.
For any sufficiently small $\epsilon>0$
and for any point $t\in \Delta^{\ast}(\epsilon)$,
$\pi^{-1}(t)$ is smooth.

Let us consider the inclusion
$\iota_1:\pi^{-1}(\Delta(\epsilon))\setminus p^{-1}(D)
\lrarr
 p^{-1}(N_0)\setminus p^{-1}(D)$.
Let $P_1,\ldots,P_m$ denote the points of $H_1\cap D$.
We can take continuous maps
$\varphi_i:S^1\times S^1\lrarr
\del\pi^{-1}(\Delta(\epsilon))$ such that
$p^{-1}(N_0)\setminus p^{-1}(D)$
is homotopy equivalent to the topological space obtained from
$\pi^{-1}(\Delta(\epsilon)-p^{-1}(D))$
and $\coprod_{i=1}^m S^1\times D^2$
via the attaching maps $\varphi_i$ $(i=1,\ldots,m)$.
Hence $\iota_1$ induce the surjection
$\pi_1\bigl(\pi^{-1}(\Delta(\epsilon))\setminus p^{-1}(D)\bigr)
\lrarr
 \pi_1\bigl(p^{-1}(N_0)\setminus p^{-1}(D)\bigr)$.

The fiber bundle
$\pi^{-1}\bigl(\Delta^{\ast}(\epsilon)\bigr)
 \setminus p^{-1}(D_0)\lrarr \Delta^{\ast}(\epsilon)$.
induces the exact sequence:
\begin{equation}\label{eq;04.1.22.1}
\begin{CD}
 \pi_1\bigl( \pi^{-1}(t)\setminus p^{-1}(D_0) \bigr)
@>>>
 \pi_1\bigl(
 \pi^{-1}\bigl(\Delta^{\ast}(\epsilon)\bigr)
\setminus p^{-1}(D_0)
 \bigr)
@>{\pi_{\ast}}>>
 \pi_1\bigl(\Delta^{\ast}(\epsilon)\bigr)
@>>> 1.
\end{CD}
\end{equation}
The inclusion
$\iota_2:\pi^{-1}\bigl(\Delta^{\ast}(\epsilon)\bigr)
\setminus p^{-1}(D_0)
\lrarr \pi^{-1}\bigl(\Delta(\epsilon)\bigr)\setminus p^{-1}(D)$
induces the surjection:
\begin{equation}\label{eq;04.1.22.2}
\iota_{2\,\ast}: \pi_1\bigl(
 \pi^{-1}(\Delta^{\ast}(\epsilon)) \setminus p^{-1}(D_0)
 \bigr)
\lrarr
 \pi_1\bigl(
 \pi^{-1}(\Delta(\epsilon))\setminus p^{-1}(D)
 \bigr).
\end{equation}
We can take a loop $\gamma$ around $D_1$ which is 
mapped to the generator of $\pi_1\bigl(\Delta^{\ast}(\epsilon)\bigr)$
via $\pi_{\ast}$ in (\ref{eq;04.1.22.1}).
On the other hand, $\gamma$ is mapped to $0$
via the map $\iota_{2\,\ast}$.
Hence we can conclude that
the inclusion
$\pi^{-1}(t)\setminus p^{-1}(D_0)\lrarr
 \pi^{-1}(\Delta(\epsilon))\setminus p^{-1}(D)$
induces the surjection of the fundamental groups.

In all,
the natural morphism
$\pi^{-1}(t)\setminus p^{-1}(D_0)\lrarr X-D$ 
induces the surjection of the fundamental groups.
Since the morphism factors through
$H_t\setminus D$,
we obtain the surjectivity of
$\pi_1\bigl(H_t\setminus D\bigr)\lrarr \pi_1(X-D)$
for any sufficiently small $t\neq 0$.
Since we can isotopically deform $H_t\setminus D$
to $H_1\setminus D$ in $X-D$,
we can conclude that
$\pi_1(H_1\setminus D)\lrarr \pi_1(X-D)$ is surjective.
\hfill\qed

\begin{lem}\label{lem;04.1.22.10}
Let $X$ be a smooth projective variety,
and $D$ be a normal crossing divisor of $X$.
If $H$ is sufficiently ample smooth divisor such that
$H\cup D$ is normal crossing,
then the map
$\pi_1(H\setminus D)\lrarr \pi_1(X\setminus D)$ is onto.
\end{lem}
\pf
We give only a sketch of a proof.
We use an induction on a dimension of $X$.

Let $D=\bigcup_{i=1}^l D_i$ be the irreducible decomposition.
Let $N_i$ be a tubular neighbourhood of $D_i$.
Let $N_H$ denote the tubular neighbourhood of $H$.
We put $N=N_H\cup\bigl(\bigcup_{i=1}^l N_i\bigr)$.
The morphism $\pi_1(N\setminus D)\lrarr \pi_1(X-D)$ is surjective
\cite{milnor}.

Let $\gamma$ be an element of $\pi_1(N\setminus D)$.
Since each connected component $D_i\cap D_j$
intersects with $H$, 
we may decompose $\gamma$
into the product of the paths $\gamma_{\alpha}$
such that they are represented by closed paths
$\tilde{\gamma}_{\alpha}$ satisfying
$\tilde{\gamma}_{\alpha}\bigl([0,a]\bigr)\subset N_H\setminus D$,
$\tilde{\gamma}_{\alpha}\bigl([a,b]\bigr)\subset N_{i_{\alpha}}\setminus D$
and
$\tilde{\gamma}_{\alpha}\bigl([b,1]\bigr)\subset N_{H}\setminus D$.

Due to the hypothesis of our induction,
the inclusion
$D_{i_{\alpha}}\cap N_H\setminus \bigl(\bigcup_{j\neq i_{\alpha}} N_j\bigr)
\subset
 D_{i_{\alpha}}\setminus \bigl(\bigcup_{j\neq i_{\alpha}}N_j\bigr)$
induces the surjection of the fundamental group.
Since 
$N_{i_{\alpha}}\cap N_H\setminus\bigl(\bigcup_{j\neq i_{\alpha}}N_j\bigr)$
and $N_{i_{\alpha}}\cap \bigl(\bigcup_{j\neq i_{\alpha}}N_j\bigr)$
are disc bundles over
$D_{i_{\alpha}}\cap N_H\setminus \bigl(\bigcup_{j\neq i_{\alpha}} N_j\bigr)$
and
$ D_{i_{\alpha}}\setminus \bigl(\bigcup_{j\neq i_{\alpha}}N_j\bigr)$
respectively.
Thus it is easy to see that
$\tilde{\gamma}_{i_{\alpha}}$ is homotopic to a closed path
in $N_H\setminus D$.
Hence we obtain the surjectivity desired.
\hfill\qed

\begin{rem}
In the proof, we have to take the tubular neighbourhoods $N_i$ 
and $N_H$ cleanly.
We omit to give the detail.
\hfill\qed
\end{rem}

\begin{cor}
Let $X$ be a quasi projective variety.
Let $L$ be a local system on $X$.
Let $H$ be a sufficiently ample hypersurface of $X$.
Then $L$ is semisimple if and only if
the restriction $L_{|H}$ is semisimple.
\hfill\qed
\end{cor}

\subsection{Elementary geometry of $GL(r)/U(r)$}
\label{subsection;04.1.30.202}

%c12.tex
\subsubsection{The $GL(r)$-invariant metric}

We have the standard left action $\kappa$ of $GL(r)$ on $\PH(r)$:
\[
 GL(r)\times \PH(r)\lrarr \PH(r),
\quad
 (g,H)\longmapsto \kappa(g,H)=g\cdot H\cdot \lefttop{t}\bar{g}.
\]
For any point $H\in\PH(r)$,
the tangent space $T_H\PH(r)$ is naturally identified with
the vector space $\nbigh(r)$.
Let $I_r$ denote the identity matrix.
We have the positive definite metric of $T_{I_r}\PH(r)$
given by $(A,B)_{I_r}=\tr(A\cdot B)=\tr(A\cdot \lefttop{t}\bar{B})$.
It is easy to see the metric is invariant
with respect to the $U(r)$-action on $T_{I_r}\PH(r)$.

Let $H$ be any point of $\PH(r)$
and let $g$ be an element of $GL(r)$
such that $H=g\cdot \lefttop{t}\bar{g}$.
Then the metric of $T_{H}\PH(r)$ is given as follows:
\begin{equation}\label{eq;04.1.4.2}
(A,B)_{H}=\bigl(\kappa(g^{-1})_{\ast}A,\kappa(g^{-1})_{\ast}B\bigr)_{I_r}
=\tr\Bigl(
 g^{-1}A\lefttop{t}\bar{g}^{-1}\cdot g^{-1}B\lefttop{t}\bar{g}^{-1}
 \Bigr)
=\tr\bigl(H^{-1}AH^{-1}B\bigr).
\end{equation}
Since $(\cdot,\cdot)_{I_r}$ is $U(r)$-invariant,
the metric $(\cdot,\cdot)_{H}$ on $T_H\PH(r)$ is well defined.
Thus we have the $GL(r)$-invariant Riemannian metric
of $\PH(r)$.

It is well known that $\PH(r)$ with the metric above is
a symmetric space with non-positive curvature.
We denote the induced distance by $d_{\PH(r)}$.
We often use the simple notation $d$ to denote $d_{\PH(r)}$,
if there are no confusion.

Let $X$ be a manifold,
and $\Psi:X\lrarr\PH(r)$ be a differentiable map.
Let $P$ be a point of $X$,
and $v$ be an element of the tangent space $T_PX$.
\begin{lem}\label{lem;04.1.13.1}
We have the following formula:
\[
 \bigl|d\Psi(v)\bigr|^2_{T_{\Psi(P)}\PH(r)}
=\tr\Bigl(
 \Psi(P)^{-1}\cdot d\Psi(v)\cdot \Psi(P)^{-1}\cdot d\Psi(v)
 \Bigr)
=\tr\Bigl(
 \overline{\Psi(P)}^{-1}\cdot
 \overline{d\Psi(v)}
 \cdot 
 \overline{\Psi(P)}^{-1}\cdot
 \overline{d\Psi(v)}
 \Bigr).
\]
\end{lem}
\pf
It follows from (\ref{eq;04.1.4.2}).
\hfill\qed

\subsubsection{The geodesics and some elementary estimates of the distances} 
\label{subsubsection;a12.28.10}

Let $\vecalpha=(\alpha_1,\ldots,\alpha_r)$ be a tuple of real numbers.
Let $\gamma_{\vecalpha}(t)$ denotes the diagonal matrices
whose $(i,i)$-th component is $e^{\alpha_i\cdot t}$.
When we regard $\gamma_{\vecalpha}$ as a path in $\PH(r)$,
it is well known that $\gamma_{\vecalpha}$ is a geodesic.
We put $s:=\sqrt{\sum \alpha_i^2}\cdot t$,
and then $|s|$ gives the arc length from $I_r$.
We put $\widetilde{\gamma}_{\vecalpha}(s):=\gamma_{\vecalpha}(t)$.

Let us consider the case $\alpha_1>\alpha_2>\cdots>\alpha_r$.
Let $k$ be an upper triangular matrices
whose diagonal entries are $1$.
i.e.,
$k_{i\,j}=0$ unless $i\leq  j$
and $k_{i\,i}=1$.
The following lemma can be checked directly.
\begin{lem}\label{lem;a12.28.1}
We put $C_2:=\min\bigl\{\alpha_i-\alpha_{i+1}\bigr\}>0$.
Then the following holds:
\[
 \big|\big|
 I-\gamma_{\vecalpha}(t)^{-1/2}\cdot k\cdot \gamma_{\vecalpha}(t)^{1/2}
 \big|\big|
\leq C_1\cdot e^{-C_2\cdot t}.
\]
Here $||\cdot||$ denote the norm of $M(r)$.
\hfill\qed
\end{lem}

\begin{lem}\label{lem;a12.28.2}
There exist positive numbers $C_1$ and $C_2$,
independent of $k$,
such that the following holds:
\[
 d_{\PH(r)}\Bigl(
 \gamma_{\vecalpha}(t),\,
 \kappa\bigl(k,\gamma_{\vecalpha}(t)\bigr)
 \Bigr)\leq C_1\cdot e^{-C_2\cdot t}.
\]
\end{lem}
\pf
We put
 $k(t):=
 \gamma_{\vecalpha}(t)^{-1/2}\cdot k\cdot \gamma_{\vecalpha}(t)^{1/2}$.
We have the following:
\[
 d_{\PH(r)}\Bigl(
 \gamma_{\vecalpha}(t),\,\kappa\bigl(k,\gamma_{\vecalpha}(t)\bigr)
 \Bigr)
=d_{\PH(r)}\Bigl(
 I_r,\,k(t)\!\cdot\! \lefttop{t}\overline{k(t)}
 \Bigr).
\]
Then the claim immediately follows from
Lemma \ref{lem;a12.28.1}.
\hfill\qed

\begin{lem} \label{lem;a12.28.3}
Let $s, s_0$ be non-negative numbers such that $s\geq s_0$.
Let $k$ be as above.
Then we have the following inequality
for some positive numbers $C_1$ and $C_2$
which are independent of $k$ and $s_0$:
\[
 \Bigl|
 d\Bigl(
 \widetilde{\gamma}_{\vecalpha}(s),\,
 \kappa\bigl(k,
 \widetilde{\gamma}_{\vecalpha}(s_0)
 \bigr)
 \Bigr)-s+s_0
 \Bigr|
\leq C_1\cdot e^{-C_2\cdot s}.
\]
\end{lem}
\pf
We have the following triangle inequality:
\[
 \Bigl|
 d\Bigl(\widetilde{\gamma}_{\vecalpha}(s),\,
 \kappa(k,\widetilde{\gamma}_{\vecalpha}(s_0))\Bigr)
-d\Bigl(
 \kappa\bigl(k,\widetilde{\gamma}_{\vecalpha}(s)\bigr),\,
 \kappa\bigl(k,\widetilde{\gamma}_{\vecalpha}(s_0)\bigr)\Bigr)
 \Bigr|
 \leq
d\Bigl(
 \widetilde{\gamma}_{\vecalpha}(s),\,
  \kappa\bigl(k,\widetilde{\gamma}_{\vecalpha}(s)\bigr)\Bigr).
\]
Due to Lemma \ref{lem;a12.28.2},
the right hand side is dominated by $C_1\cdot e^{-C_2s}$
for some positive constants $C_1$ and $C_2$.
On the other hand, 
we have $d\Bigl(
 \kappa\bigl(k,\widetilde{\gamma}_{\vecalpha}(s)\bigr),\,
 \kappa\bigl(k,\widetilde{\gamma}_{\vecalpha}(s_0)\bigr)\Bigr)=s-s_0$.
Thus we are done.
\hfill\qed

\vspace{.1in}
Let $B_{\widetilde{\gamma}_{\vecalpha}}$ be the Busemann function
given as follows (see \cite{Eberlein} or \cite{sakai}, for example):
\[
 B_{\widetilde{\gamma}_{\vecalpha}}(x):=
 \lim_{s\to\infty}
 \Bigl(
 d\bigl(\widetilde{\gamma}_{\vecalpha}(s),x\bigr)-s
 \Bigr).
\]
Then Lemma \ref{lem;a12.28.3} is reformulated as follows:
\begin{equation}\label{eq;a12.28.4}
 B_{\widetilde{\gamma}_{\vecalpha}}\bigl(
 \kappa\bigl(k, \widetilde{\gamma}_{\vecalpha}(s_0)\bigr)\bigr)
=-s_0
=B_{\widetilde{\gamma}_{\vecalpha}}\bigl(
 \widetilde{\gamma}_{\vecalpha}(s_0)
 \bigr).
\end{equation}

\begin{lem}\label{lem;a12.28.6}
We have the following inequality:
\[
 d\bigl(I_r,\kappa\bigl(k,\widetilde{\gamma}_{\vecalpha}(s_0)\bigr)\bigr)
\geq d\bigl(I_r,\widetilde{\gamma}_{\vecalpha}(s_0) \bigr)=s_0.
\]
\end{lem}
\pf
It follows from (\ref{eq;a12.28.4}).
Note that the horospheres and $\widetilde{\gamma}$ are orthogonal.
See \cite{sakai} for more detail.
\hfill\qed

\subsubsection{An estimate of infimum}
\label{subsubsection;04.1.6.2}

Let $A$ be an element of $GL(r)$.
Let $a_1,\ldots, a_r$ be eigenvalues of $A$.
We put as follows:
\[
 \rho(A):=\Bigl(
 \sum \bigl(\log|a_i|^2\bigr)^2
 \Bigr)^{1/2}.
\]

\begin{lem}\label{lem;a12.28.7}
Assume $|a_i|>|a_{i+1}|$ for any $i$.
Then we have the following inequality,
for any $H\in\PH(r)$:
\begin{equation}\label{eq;a12.28.5}
 \rho(A)\leq
 d_{\PH(r)}\Bigl( \kappa(A,H),\,H \Bigr).
\end{equation}
\end{lem}
\pf
We have only to show the inequality (\ref{eq;a12.28.5})
in the case $H=I_r$ for any $A$.
For any element $U\in U(r)$,
we have the following:
\[
 d\Bigl(
 \kappa\bigl(UAU^{-1},I_r\bigr),\,I_r
 \Bigr)
=d\Bigl(
 \kappa\bigl(A,I_r\bigr),\,I_r
 \Bigr),
\quad
 \rho\bigl(UAU^{-1}\bigr)=\rho(A).
\]
Hence we may assume that $A$ is an upper triangular matrices
such that whose $(i,i)$-th entries are $a_i$.
Then we can decompose $A$ into the product $A_1\cdot A_2$
such that the following holds:
\begin{itemize}
\item
$A_1$ is the upper triangular matrix
whose $(i,i)$-th entries are $|a_i|$.
\item
$A_2$ is the diagonal matrix
such that the absolute values of the $(i,i)$-th entries are $1$.
In particular, $A_2$ is unitary.
\end{itemize}
Then it is easy to see that
$\rho(A)=\rho(A_1)$
and
$d\bigl(\kappa\bigl(A,I_r\bigr),I_r\bigr)
=d\bigl(\kappa\bigl(A_1,I_r\bigr),I_r\bigr)$.
Hence we may assume
$a_i=e^{\alpha_i}$ for some real numbers $\alpha_i$
from the beginning.
Note we have $\alpha_1>\cdots>\alpha_r$,
due to our assumption
$|a_i|>|a_{i+1}|$ for any $i$.

We decompose $A$ into the product $K\cdot A_0$
such that the following holds:
\begin{itemize}
\item
 $K$ is the upper triangular matrix whose diagonal entries are $1$.
\item
 $A_0$ is the diagonal matrix whose $(i,i)$-th component is $a_i$.
\end{itemize}

Due to Lemma \ref{lem;a12.28.6},
we have the following inequality:
\[
 d\bigl(I,\,A\!\cdot\! \lefttop{t}\!\bar{A}\bigr)
\geq
 d\bigl(I,\,A_0\!\cdot\!\lefttop{t}\!\bar{A}_0\bigr)
=\Bigl(
 \sum\bigl(\log|a_i|^2\bigr)^2
 \Bigr)^{1/2}
=\rho(A).
\]
Thus we are done.
\hfill\qed

\begin{cor}\label{cor;a12.28.9}
Let $A$ be any element of $GL(r)$.
Then we have the inequality
$\rho(A)\leq d\bigl(\kappa(A,H),H\bigr)$
for any $H\in\PH(r)$.
\end{cor}
\pf
Let $a_1,\ldots,a_r$ be eigenvalues of $A$.
In the case $|a_i|>|a_{i+1}|$,
the claim is already shown in Lemma \ref{lem;a12.28.7}.
Let us take a sequence $\{A^{(n)}\}$
in $GL(r)$ such that the following holds:
\begin{itemize}
\item
 The sequence $\{A^{(n)}\}$ converges to $A$.
\item
 Let $a^{(n)}_1,\ldots,a^{(n)}_r$ be eigenvalues of $A^{(n)}$.
 Then the inequality $|a^{(n)}_i|>|a^{(n)}_{i+1}|$ hold for any $i$.
\end{itemize}
Then we have the inequalities
$\rho(A^{(i)})\leq d\bigl(\kappa(A^{(i)},H),H\bigr)$.
We also have the convergences
$\rho(A^{(i)})\lrarr \rho(A)$
and
$d\bigl(\kappa(A^{(i)},H),H\bigr)\lrarr d\bigl(\kappa(A,H),H\bigr)$.
Thus we are done.
\hfill\qed

\begin{lem}\label{lem;a12.28.11}
Let $A$ be any element of $GL(r)$.
Then we have the following equality:
\begin{equation}\label{eq;a12.28.8}
 \rho(A)=
 \inf\bigl\{
 d(\kappa(A,H),H)\,\big|\,
 H\in\PH(r)
 \bigr\}.
\end{equation}
\end{lem}
\pf
We have already shown the inequality $\leq$ in (\ref{eq;a12.28.8})
(Corollary \ref{cor;a12.28.9}).
Let us show the inequality $\geq$.
We may assume that $A$ is an upper triangular matrix
whose $(i,i)$-th entries are $a_i$.
Let $\vecalpha$ be any element of $\real^r$
such that
$\alpha_i>\alpha_{i+1}$ for any $i$,
and let us consider the geodesic
$\gamma_{\vecalpha}(t)$
(the subsubsection \ref{subsubsection;a12.28.10}).
We put
$A(t):=
 \gamma_{\vecalpha}(t)^{-1/2}\cdot A\cdot
 \gamma_{\vecalpha}(t)^{1/2}$.
Then $A(t)$ converges to the diagonal matrix $C$
whose $(i,i)$-th entries are $a_i$.
Since $A(t)\cdot \lefttop{t}\overline{A(t)}$
converges to $C\cdot \overline{C}$,
we have the following convergence:
\[
 d\Bigl(
 \kappa\bigl(A,\gamma_{\vecalpha}(t)\bigr),\,
 \gamma_{\vecalpha}(t)\Bigr)
=d\bigl(
 A(t)\!\cdot\! \lefttop{t}\overline{A(t)},\,I_r
 \bigr)
\lrarr d\bigl(C\cdot \overline{C},\,I_r\bigr)=\rho(A).
\]
Thus we are done.
\hfill\qed

\subsubsection{A lower bound of the energy}

Let $A$ be an element of $GL(r)$,
and $H$ be an element of $\PH(r)$.
Let $\varphi:[0,2\pi]\lrarr\PH(r)$ be an $L_1^2$-map
such that $\varphi(0)=H$ and $\varphi(2\pi)=\kappa(A,H)$.
Let $g$ be a positive continuous function on $[0,2\pi]$.
The following inequality will often be used.
\begin{lem}\label{lem;a12.29.10}
\[
\int_0^{2\pi}
 \left|
 \frac{\del \varphi}{\del t}
 \right|^2\cdot g\cdot dt
\geq
 \rho(A)^{2}\cdot
 \left(
 \int g^{-1}dt
 \right)^{-1}.
\]
\end{lem}
\pf
Due to the Schwarz's inequality,
we have the following:
\[
 \left(
 \int_0^{2\pi}\left|
 \frac{\del \varphi}{\del t}
 \right|
 \right)^2
\leq
 \int_0^{2\pi}
 \left|\frac{\del \varphi}{\del t}\right|^2
 \cdot g\cdot dt
\times
 \int_0^{2\pi}
 g^{-1}\cdot dt.
\]
Then the claim follows from Lemma \ref{lem;a12.28.11}.
\hfill\qed

%c12.1.tex

\subsubsection{Comparison of the norm and the distance}

\label{subsubsection;04.1.30.300}

For any elements $H_1$ and $H_2$ of $\PH(r)$,
we have an element $g\in GL(r)$
such that $\kappa(g,H_1)=I_r$
and that $\kappa(g,H_2)$ is the diagonal matrix.
The set of the eigenvalues of $\kappa(g,H_2)$ is independent
of a choice of $g$.
Let $e^{\alpha_1},\ldots,e^{\alpha_r}$ be 
the eigenvalues of $\kappa(g,H_2)$.
We put as follows:
\[
 \delta(H_1,H_2):=
\left(
 \sum_{i=1}^r \Bigl(
 \frac{e^{\alpha_i}-e^{-\alpha_i}}{2}
 \Bigr)^2
\right)^{1/2}.
\]
On the other hand,
we have the distance
$d_{\PH(r)}(H_1,H_2)=\bigl(\sum \alpha_i^2\bigr)^{1/2}$.

For any real number $R$, we put as follows:
\[
 C(R):=\frac{e^{R}-e^{-R}}{2R}.
\]
If $0\leq x\leq R$,
we have $x\leq C(R)\cdot x$.
We also note that $C(R)\to 1$ when $R\to 0$.

\begin{lem}\mbox{{}}\label{lem;a12.28.12}
\begin{itemize}
\item
 We have the inequality:
\[
 d_{\PH(r)}(H_1,H_2)
\leq
 \delta(H_1,H_2).
\]
\item
 If $d(H_1,H_2)\leq R$,
 we have the inequality:
\[
 \delta(H_1,H_2)\leq
 C(R)\cdot d_{\PH(r)}(H_1,H_2).
\]
\end{itemize}
\end{lem}
\pf
It can be checked elementarily.
\hfill\qed

\vspace{.1in}
We reformulate Lemma \ref{lem;a12.28.12} as follows:
Let $V$ be an $r$-dimensional vector space,
and let $h_1$ and $h_2$ be hermitian metrics of $V$.
The identity map induces the map
$\Phi:(V,h_1)\lrarr (V,h_2)$.
We have the norms $|\Phi|$ and $\big|\Phi^{-1}\big|$.

\begin{lem}\mbox{{}}
\begin{itemize}
\item
 The following inequality holds:
 \[
  d_{\PH(r)}(h_1,h_2)^2
\leq
 \frac{
 \big|\Phi\big|^2
+\big|\Phi^{-1}\big|^{2}
-2r}{4}.
 \]
\item
 If $d_{\PH(r)}(h_1,h_2)\leq R$,
 the following inequality holds:
\[
 \frac{
  \big|\Phi\big|^2+\big|\Phi^{-1}\big|^2-2r}{4}
\leq
 C(R)^2\cdot
 d_{\PH(r)}(h_1,h_2)^2.
\]
\end{itemize}
\end{lem}
\pf
If we take an appropriate base of $V$,
$h_1$ and $h_2$ are represented by 
the identity matrix $I_r$
and the diagonal matrices whose diagonal entries are
$e^{\alpha_1},\ldots,e^{\alpha_r}$.
It is easy to check that
$|\Phi|^2=\sum_{i=1}^r e^{2\alpha_i}$
and $|\Phi^{-1}|^2=\sum_{i=1}^r e^{-2\alpha_i}$.
Then it immediately follows
$4^{-1}\bigl(
 |\Phi|^2+|\Phi^{-1}|^2-2r
 \bigr)=\delta(H_1,H_2)^2$.
Thus the claims follow from Lemma \ref{lem;a12.28.12}.
\hfill\qed

\subsection{Maps associated to commuting tuple of endomorphisms}
\label{subsection;04.1.30.203}

%c12.2.tex

\subsubsection{Preliminary}

Let $V$ be an $r$-dimensional vector space.
Let $\vecv=(v_i)$ be a frame of $V$.
For any endomorphism $f$ of $V$,
we have the matrix $A(f)\in M_r$ determined by 
$f\cdot \vecv=\vecv\cdot A(f)$,
i.e.,
$f(v_i)=\sum A(f)_{j\,i}\cdot v_j$.

Let $f_i$ $(i=1,2)$ be elements of $End(V)$
such that $f_1\circ f_2=f_2\circ f_1$.
We decompose $f_i$ into the product
of the unipotent part $f_i^u$ and the semisimple part $f_i^s$.
Then $f_1^u$, $f_1^s$, $f_2^u$ and $f_2^s$ are mutually commutative.
There exists an appropriate frame $\vecv$ of $V$
such that $A(f_i^s)$ $(i=1,2)$ are diagonal matrices,
and $A(f_i^u)$ are upper triangular matrices.
Then the matrices $A\bigl(\log f_i^u\bigr)$
are also upper triangular matrices.
In the following,
we identify the endomorphisms and the matrices
via the frame $\vecv$ above.

We put $\eta_i^n:=(2\pi)^{-1}\cdot \log f_i^u$.
We also have $\eta_i^s$ satisfying
$\exp\bigl(2\pi\eta_i^s\bigr)=f_i^s$
such that $0\leq \Image(\alpha)<1$ holds
for any eigenvalues $\alpha$ of $\eta_i^s$.
We put $\eta_i:=\eta_i^s+\eta_i^n$,
and then we have $\exp\bigl(2\pi\eta_i\bigr)=f_i$.

Let $\varphi:\real^2\lrarr GL(V)$ be the morphism
given by
$ \varphi(\theta_1,\theta_2)
:=\exp\bigl(
 \theta_1\!\cdot\!\eta_1+\theta_2\!\cdot\!\eta_2
 \bigr)$.
We also put as
$ \varphi^s(\theta_1,\theta_2)
:=\exp\bigl(
 \theta_1\!\cdot\!\eta_1^s+\theta_2\!\cdot\!\eta_2^s
 \bigr)$ and
$\varphi^u(\theta_1,\theta_2)
:=\exp\bigl(
 \theta_1\!\cdot\!\eta_1^n+\theta_2\!\cdot\!\eta_2^n
 \bigr)$.
We have $\varphi=\varphi^s\cdot \varphi^u$.
Under the identification $GL(V)=GL(r)$ via the frame $\vecv$ above,
$\varphi^u(\theta_1,\theta_2)$ are upper triangular
whose diagonal entries are $1$,
and $\varphi^s(\theta_1,\theta_2)$ are diagonal matrices.

\subsubsection{Construction}

\label{subsubsection;04.1.27.100}

Let us take an element $\vecalpha=(\alpha_1,\ldots,\alpha_r)\in\real^r$
such that $\alpha_i>\alpha_{i+1}$ for any $i$.
We put $\beta:=\min\bigl\{\alpha_i-\alpha_{i+1}\bigr\}>0$.
We have the $C^{\infty}$-map $F:\real\times\real^2\lrarr\PH(r)$
given as follows:
\[
 F(t,\theta_1,\theta_2):=
 \kappa\Bigl(
 \varphi(\theta_1,\theta_2),\gamma_{\vecalpha}(t)
 \Bigr).
\]
We put
$\widetilde{\varphi}(t,\theta_1,\theta_2)=
 \gamma_{\vecalpha}(t)^{-1/2}\cdot\varphi(\theta_1,\theta_2)\cdot \gamma_{\vecalpha}(t)^{1/2}$.
Similarly we obtain $\widetilde{\varphi}^u$ and $\widetilde{\varphi}^s$.
We have $\widetilde{\varphi}^s=\varphi^s$.
We also have the following:
\[
 \gamma_{\vecalpha}(t)^{-1/2}\cdot F(t,\theta_1,\theta_2)\cdot \gamma_{\vecalpha}(t)^{1/2}
=\widetilde{\varphi}(t,\theta_1,\theta_2)
\cdot\lefttop{t}\overline{\widetilde{\varphi}(t,\theta_1,\theta_2)}.
\]

\begin{lem}\label{lem;a12.28.15}
For any $C_0$, there exists a positive number $C_1$
such that
the following inequalities hold,
in the case $|\theta_1|+|\theta_2|<C_0$:
\[
 \big|\big|\widetilde{\varphi}^u(t,\theta_1,\theta_2)-I_r\big|\big|
 \leq C_1\cdot e^{-\beta\cdot t}.
\]
\[
\bigl|\bigl|
 \widetilde{\varphi}(t,\theta_1,\theta_2)-
 \varphi^s(\theta_1,\theta_2) \bigr|\bigr|
\leq C_1\cdot e^{-\beta\cdot t}.
\]
\[
 d_{\PH(r)}\Bigl(
 \gamma_{\vecalpha}(t)^{-1/2}\!\cdot\! F(t,\theta_1,\theta_2)\!\cdot\!
 \gamma_{\vecalpha}(t)^{1/2},\,\,
 \varphi^s(\theta_1,\theta_2)\!\cdot\!
 \lefttop{t}\overline{\varphi^s(\theta_1,\theta_2)}
 \Bigr)
\leq C_1\cdot e^{-\beta\cdot t}.
\]
\end{lem}
\pf
Since $\varphi^u(\theta_1,\theta_2)$ are upper triangular
whose diagonal entries are $1$,
it is easy to check.
\hfill\qed

\subsubsection{Estimate of derivatives}

We see the estimate of the derivative of $F$.
We have 
$ \frac{\del F}{\del \theta_i}=
 \eta_i\cdot F+F\cdot\lefttop{t}\overline{\eta_i}$.

\begin{lem}\label{lem;a12.30.5}
For any positive number $C_0$,
there exists a positive number $C_1$ such that
the following holds in the case $|\theta_1|+|\theta_2|<C_0$:
\[
 \left|\left|
 \gamma_{\vecalpha}(t)^{-1/2}\cdot\frac{\del F}{\del \theta_i}\cdot \gamma_{\vecalpha}(t)^{1/2}
-\varphi^s(\theta_1,\theta_2)\cdot
 \bigl(
 2\Re\eta_i^s
 \bigr)\cdot \lefttop{t}\overline{\varphi^s(\theta_1,\theta_2)}
 \right|\right|
\leq
 C_1\cdot e^{-\beta\cdot t}.
\]
As a result,
we have the following:
\[
 \left|
 \frac{\del F}{\del \theta_i}
 \right|
=2\left|  \Re\eta_i^s \right|+O\bigl(e^{-\beta t}\bigr)
=\frac{\rho(f_i)}{2\pi}+O\bigl(e^{-\beta t}\bigr).
\]
\end{lem}
\pf
We have the following convergence when $t\to\infty$:
\begin{multline}
 \gamma_{\vecalpha}(t)^{-1/2}\cdot
 \frac{\del F}{\del\theta_i}\cdot
  \gamma_{\vecalpha}(t)^{1/2}\\
=\bigl(
 \gamma_{\vecalpha}(t)^{-1/2}\cdot\eta_i\cdot \gamma_{\vecalpha}(t)^{1/2}
 \bigr)
\cdot
 \bigl(
 \gamma_{\vecalpha}(t)^{-1/2}\cdot F\cdot\gamma_{\vecalpha}(t)^{-1/2}
 \bigr)
+\bigl(
 \gamma_{\vecalpha}(t)^{-1/2}\cdot F\cdot \gamma_{\vecalpha}(t)^{-1/2}
 \bigr)
 \cdot
\lefttop{t}
 \bigl(
  \overline{\gamma_{\vecalpha}(t)^{-1/2}\cdot 
    \eta_i\cdot \gamma_{\vecalpha}(t)^{1/2}}\bigr)\\
\lrarr \eta_i^s\cdot
 \varphi^s(\theta_1,\theta_2)
\lefttop{t}\overline{\varphi^s(\theta_1,\theta_2)}
+\varphi^s(\theta_1,\theta_2)
\lefttop{t}\overline{\varphi^s(\theta_1,\theta_2)}
\cdot \lefttop{t}\overline{\eta_i^s}
=2\varphi^s(\theta_1,\theta_2)\cdot
 \Re(\eta_i^s)\cdot
 \lefttop{t}\overline{\varphi^s(\theta_1,\theta_2)}.
\end{multline}
Here the convergence is estimated by $e^{-\beta t}$.
Thus we are done.
\hfill\qed

\begin{lem}\label{lem;a12.30.6}
We have the following:
\[
 \left|
 \frac{\del F}{\del t}
 \right|^2
=\left|
 \frac{d \gamma_{\vecalpha}(t)}{dt}
 \right|^2
=\sum\alpha_i^2.
\]
\end{lem}
\pf
It follows from the $GL(r)$-invariance of the metric of
$\PH(r)$.
\hfill\qed

%c12.22.tex

\subsubsection{Extension}
\label{subsubsection;a12.30.10}

We use the real coordinate $\zeta_i=\xi_i+\sqrt{-1}\eta_i$ $(i=1,2)$
for $\hyperh^2$.
Let $A$, $B$, and $C_i$ $(i=1,2)$ be positive numbers.
Let us consider the morphism
$\Phi:\hyperh^2\lrarr\PH(r)$ given as follows:
\begin{equation} \label{eq;04.1.22.6}
 \Phi(\zeta_1,\zeta_2):=
 \gamma_{\vecalpha}\bigl(B\cdot\log C_1\eta_1\bigr)^{-1/2}\cdot
 F\bigl(B\cdot\log(C_1\eta_1+C_2\eta_2+A),\xi_1,\xi_2\bigr)
 \cdot
 \gamma_{\vecalpha}\bigl(B\cdot\log C_1\eta_1\bigr).
\end{equation}
We also put as follows:
\[
 \varphi^{\star}(\eta_1,\xi_1,\xi_2)
:=\tilde{\varphi}\bigl(B\cdot\log C_1\eta_1,\xi_1,\xi_2\bigr).
\]
Then we have the following estimate,
due to Lemma \ref{lem;a12.28.15}:
\[
 \bigl|\bigl|
 \varphi^{\star}(\eta_1,\xi_1,\xi_2)
-\varphi^s(\xi_1,\xi_2)
 \bigr|\bigr|
=O\bigl(\eta_1^{-\beta\cdot B}\bigr).
\]

\begin{lem}\label{lem;a12.30.1}
We have the following:
\[
\lim_{\eta_1\to\infty}\Phi(\zeta_1,\zeta_2)
=\varphi^s(\xi_1,\xi_2)\cdot\lefttop{t}\overline{\varphi^s(\xi_1,\xi_2)}.
\]
\end{lem}
\pf
We have the following:
\[
 \Phi(\zeta_1,\zeta_2)=
\varphi^{\star}(\eta_1,\xi_1,\xi_2)\cdot
\gamma_{\vecalpha}\bigl(
 B\cdot \log \bigl(C_1\eta_1+C_2\eta_2+A\bigr)
-B\cdot\log\bigl(C_1\eta_1\bigr)
 \bigr)
\cdot\lefttop{t}\overline{\varphi^{\star}(\eta_1,\xi_1,\xi_2)}.
\]
Then it is easy to check the claim.
\hfill\qed

\vspace{.1in}
We take the isomorphism $\real_{\geq 0}\lrarr \closedopen{0}{1}$
given by $\eta_1\longmapsto \kappa=(1+\eta_1)^{-1}\cdot\eta_1$.
It induces the isomorphism
$\hyperh^2\simeq\bigl(\real\times\closedopen{0}{1}\bigr)\times \hyperh$,
and hence we obtain the map
$\tilde{\Phi}:\bigl(\real\times\closedopen{0}{1}\bigr)
      \times \hyperh\lrarr\PH(r)$.
\begin{lem}
The map $\tilde{\Phi}$ can be naturally extended
to the continuous map $\real\times[0,1]\times\hyperh\lrarr\PH(r)$.
We denote the extended map by $\tilde{\Phi}$.
On $\real\times\{1\}\times\hyperh$,
we have
$\tilde{\Phi}(\xi_1,1,\zeta_2)
=\varphi^s(\xi_1,\xi_2)\cdot\lefttop{t}\overline{\varphi^s(\xi_1,\xi_2)}$.
\end{lem}
\pf
It follows from Lemma \ref{lem;a12.30.1}.
\hfill\qed

\vspace{.1in}
We have the map
$\tilde{f}_2(\kappa):\closedopen{0}{1}
 \lrarr GL(r)$ given as follows:
\[
 \tilde{f}_2(\kappa)=
\gamma_{\vecalpha}\Bigl(
 B\cdot\log\frac{c_1\kappa}{1-\kappa}
 \Bigr)^{-1/2}
\cdot
 f_2\cdot
 \gamma_{\vecalpha}\Bigl(
 B\cdot\log\frac{c_1\kappa}{1-\kappa}
 \Bigr)^{1/2}.
\]
It is naturally extended to the continuous map
$[0,1]\lrarr GL(r)$,
which we denote also by $\tilde{f}_2(\kappa)$.
We have $\tilde{f}_2(1)=f_2^s$.

Let $\tilde{\Phi}_{(\xi_1,\kappa)}$
denote the  restriction of $\tilde{\Phi}$ to
$\{(\xi_1,\kappa)\}\times\hyperh$.
We have the $\seisuu$-action
on $\hyperh$ by $n\cdot (\xi_2,\eta_2)\lrarr (\xi_2+n,\eta_2)$.
For $(\xi_1,\kappa)\in \real\times[0,1]$,
the action of $\seisuu$ on $\PH(r)$ is given by
$\tilde{f}_2(\kappa)$.
Then we can regard $\tilde{\Phi}_{(\xi_1,\kappa)}$
as the twisted map with respect to the actions,
i.e.,
$\tilde{\Phi}_{(\xi_1,\kappa)}:
 \Deltabar^{\ast}\lrarr\PH(r)/\langle \tilde{f}_2(\kappa)\rangle$.
(See the subsubsection \ref{subsubsection;04.1.22.5}
for twisted maps.)
Thus we obtain the continuous family of $C^{\infty}$-twisted maps
$\bigl\{\tilde{\Phi}_{(\xi_1,\kappa)}\,
 \big|\,(\xi_1,\kappa)\in \real\times[0,1]\bigr\}$.
\label{subsubsection;a12.30.11}

\begin{lem}
We have the estimate,
which is independent of $(\xi_1,\kappa)$:
\[
 \left|
 \frac{\del \tilde{\Phi}_{(\xi_1,\kappa)}}{\del \xi_2}
 \right|
=\frac{\rho(f_2)^2}{4\pi^2}
+O\Bigl(\bigl(C_1\eta_1+C_2\eta_2+A\bigr)^{-B\cdot \beta}\Bigr),
\]
\[
 \left|
 \frac{\del \tilde{\Phi}_{(\xi_1,\kappa)}}{\del \eta_2}
 \right|
=O\Bigl(
 \bigl(C_1\eta_1+C_2\eta_2+A\bigr)^{-1}
 \Bigr).
\]
Here we put $\eta_1=\kappa\cdot (1-\kappa)^{-1}$.
\end{lem}
\pf
It from Lemma \ref{lem;a12.30.5} and Lemma \ref{lem;a12.30.6}.
\hfill\qed

\vspace{.1in}
Let $g_0:=(2\eta_2+b)^{-2}\cdot \bigl(d\xi_2^2+d\eta_2^2\bigr)$
be the Poincar\'{e} metric of $\Deltabarast$,
where $b$ denotes a positive constant.

\begin{lem}\label{lem;a12.30.12}
Assume the constant $B$ in {\rm(\ref{eq;04.1.22.6})}
is sufficiently large.
There exists a positive constant $C$, which is independent of
$(\xi_1,\kappa)$,
such that the following holds:
\begin{equation}\label{eq;a12.30.7}
 \int_{\Deltabarast}
 \left|
 e\bigl(\tilde{\Phi}_{(\xi_1,\kappa)}\bigr)
 -\frac{\rho(f_2)^2}{4\pi^2}\cdot (2\eta_2+b)^2
 \right|
 \dvol_{g_0}<C.
\end{equation}
\end{lem}
\pf
It follows from the previous lemma.
\hfill\qed

\subsection{Preliminary for harmonic maps and harmonic bundles}

\label{subsection;04.1.30.204}

%c13.2.tex

\subsubsection{The energy function is dominated by the energy}

Recall the section 8 in \cite{ES}.
Let $(N,g_N)$ be a Riemannian manifold with non-positive curvature.
Let $g$ be a Riemannian metric of $\Delta(R)^n$,
and let $f:\Delta(R)^n\lrarr N$ be a differentiable map.
Let $e(f)$ be the energy function of $f$
with respect to the metrics $g$ and $g_N$.
Assume that $e(f)$ is integrable,
and let $E(f)$ denote the energy.
\begin{lem}\label{lem;a12.29.4}
Let us take a positive number $R'<R$,
and then there exists a positive number $C>0$
such that the following inequality holds 
for any $f$ and for any point $P\in\Delta(R')^n$:
\[
 e(f)(P)\leq C\cdot E(f).
\]
\end{lem}
\pf
We only give a remark.
In \cite{ES},
it is discussed for the map $f:(M,g_M)\lrarr (N,g_N)$
for a compact Riemannian manifold $(M,g_M)$.
Since the argument is local,
it can be applied in our case.
\hfill\qed

%c13.22.tex

\subsubsection{Twisted map and twisted harmonic map}
\label{subsubsection;04.1.22.5}

Recall the correspondence of harmonic metrics
and twisted harmonic maps.
Let $X$ be a $C^{\infty}$-manifold and
let $(E,\nabla)$ be a flat connection on $X$.
Let $h$ be a continuous hermitian metric of $E$.
Let $\pi:\tilde{X}\lrarr X$ be a universal covering.
Then we obtain the pull backs $\pi^{-1}(E,\nabla,h)$.
Once we take a flat frame $\vecv$ of $\pi^{-1}(E,\nabla)$,
then the metric $\pi^{-1}h$ induces the $C^{\infty}$-map
$\tilde{\Psi}_h:\tilde{X}\lrarr \PH(r)$.
If $h$ is $C^{\infty}$, then $\tilde{\Psi}_h$ is also $C^{\infty}$.

We have the $\pi_1(X)$-action on $\tilde{X}$.
The monodromy induces endomorphism $\rho:\pi_1(X)\lrarr GL(r)$,
which induces the $\pi_1(X)$-action on $\PH(r)$.
Since the map $\tilde{\Psi}_h$ is equivariant,
we obtain the continuous map
$\Psi_h:X\lrarr \PH(r)\big/\rho\bigl(\pi_1(X)\bigr)$.
Although $\PH(r)\big/\rho\bigl(\pi_1(X)\bigr)$ is not a manifold,
we say $\Psi_h$ is $C^{\infty}$ (resp. harmonic)
if $\tilde{\Psi_h}$ is $C^{\infty}$ (resp. harmonic).
We do not explicitly distinguish $\Psi_h$ and $\tilde{\Psi}_h$.

%c13.tex

\subsubsection{The induced one form}
\label{subsubsection;04.1.5.10}

We recall some formalism following \cite{s5}.
Let $X$ be a complex manifold
with Kahler metric $g$.
Let $V$ be a $C^{\infty}$-vector bundle on $X$,
and $\nabla$ be a flat connection of $V$.
We have the decomposition $\nabla=d'+d''$
into the sum of the $(1,0)$-part and the $(0,1)$-part.

Let $K$ be a hermitian metric of $V$.
Then the differential operators $\delta'$ and $\delta''$
is given by the conditions
that $\delta'+d''$ and $d'+\delta''$ are unitary connections.
We put as follows:
\[
 \del=\frac{d'+\delta'}{2},
 \quad
 \delbar=\frac{d''+\delta''}{2},
 \quad
 \theta=\frac{d'-\delta'}{2},
 \quad
 \theta^{\dagger}=\frac{d''-\delta''}{2}.
\]
\begin{lem}
The operator $\del+\delbar$ is a unitary connection.
The one forms $\theta$ and $\theta^{\dagger}$ are mutually adjoint,
and hence $\theta+\theta^{\dagger}$ is self adjoint.
\hfill\qed
\end{lem}

\begin{lem}
We have the following vanishings:
\begin{equation}\label{eq;a12.28.20}
 (\del+\delbar)(\theta+\theta^{\dagger})=0,
\quad\quad
 (\del+\delbar)^2+(\theta+\theta^{\dagger})^2=0.
\end{equation}
The first equality means
the vanishing of the commutator
$[\del+\delbar,\,\,\theta+\theta^{\dagger}]$
of the operators.
\end{lem}
\pf
We decompose the equality $\nabla^2=0$
into the self adjoint part and the anti-self adjoint part.
Then we obtain the equalities.
\hfill\qed

\begin{cor}\label{cor;04.1.11.2}
We have the following vanishings:
\[
 \del\theta=0,\quad \delbar\theta^{\dagger}=0,
\quad \del\theta^{\dagger}+\delbar\theta=0,
\]
\[
 \del^2+\theta^2=0,
\quad
 \delbar^2+\theta^{\dagger\,2}=0,
\quad
 \del\delbar+\delbar\del+[\theta,\theta^{\dagger}]=0.
\]
\end{cor}
\pf
We have only to decompose the equalities (\ref{eq;a12.28.20})
into $(2,0)$-parts, $(1,1)$-parts and $(0,2)$-parts.
\hfill\qed

\begin{rem}
If we have $\theta^2=0$ and $\delbar\theta=0$,
we also have the other vanishings
$(\theta^{\dagger})^2=\delbar^2=\del^2=\del\theta^{\dagger}=0$.
In other words, the metric $K$ is pluri-harmonic.
\hfill\qed
\end{rem}

\subsubsection{The description by matrices}

Let $\vecv$ be a flat frame of $V$.
Let $H=(h_{i\,j})$ be a $\PH(r)$-valued function
determined by $H=H(K,\vecv)$.
We also have the $M(r)$-valued $(1,0)$-form $\Theta$
and the $M(r)$-valued $(0,1)$-form $\Theta^{\dagger}$
determined by 
$\theta\vecv=\vecv\cdot\Theta$ and
$\theta^{\dagger}\vecv=\vecv\cdot\Theta^{\dagger}$.
We have the following:
\begin{multline}
 \delbar h_{i\,j}=
 K(d''v_i,v_j)+K(v_i,\delta'v_j)
=K\bigl(v_i,(\delta'-d')v_j\bigr) \\
=-K(v_i,2\theta v_j)=-2\sum_k K(v_i,\Theta_{k\,j}v_k)
=-2\sum_k h_{i\,k}\cdot\overline{\Theta}_{k\,j}.
\end{multline}
Thus we obtain
$\del \overline{h_{i\,j}}=-2\overline{h_{i\,k}}\cdot\Theta_{k\,j}$.
Namely we have the relation
$\Theta=-\frac{1}{2}\cdot\bar{H}^{-1}\cdot\del \overline{H}$.
We also have the following:
\begin{multline}
 \del h_{i\,j}=
 K(d'v_i,v_j)+K(v_i,\delta''v_j)
=K\bigl(v_i,(\delta''-d'')v_j\bigr)\\
=-2K(v_i,\theta^{\dagger}v_j)=-2\sum_k K(v_i,\Theta^{\dagger}_{k\,j}v_k)
=-2\sum_k h_{i\,k}\overline{\Theta}^{\dagger}_{k\,j}.
\end{multline}
Thus we obtain the relation
$\delbar \overline{h}_{i\,j}=\sum
-2\overline{h}_{i\,k}\cdot\Theta^{\dagger}_{k\,j}$.
Namely we have the relation:
$ \Theta^{\dagger}
=-\frac{1}{2}\cdot\overline{H}^{-1}\cdot \del \overline{H}$.
Thus we obtain the following relation:
\[
 \Theta+\Theta^{\dagger}=-\frac{1}{2}\cdot\overline{H}^{-1}\cdot d\overline{H}.
\]

Let $\Psi_K$ denote the twisted map $X\lrarr \PH(r)/\pi_1(X)$
associated to $K$.
Let $\vece=(e_1,\ldots,e_{2n})$ be an orthonormal base of $T_PX$.
We have the metric of $End(E)\otimes \Omega^1_X$
induced by the metrics $K$ and $g$.
By using Lemma \ref{lem;04.1.6.1} and Lemma \ref{lem;04.1.13.1},
we obtain the following:
\begin{equation}\label{eq;a12.28.45}
 8\cdot||\theta||_{K,g}^2
=4\cdot||\theta+\theta^{\dagger}||^2_{K,g}
=4\sum_i
 ||\Theta(e_i)+\Theta^{\dagger}(e_i)||^2_K
=\sum_{i}\tr\Bigl(
 \overline{H}^{-1}\cdot d\overline{H}(e_i)
\cdot \overline{H}^{-1}\cdot d\overline{H}(e_i)
 \Bigr)
=e(\Psi_K).
\end{equation}

%c13.1.tex

\subsubsection{Bochner type formula due to Corlette}
\label{subsubsection;04.1.8.30}

We put $\vartheta=\theta+\theta^{\dagger}$.
We also put $D^{+}=\del+\delbar$.
Let $\vece=(e_i)$ be an orthonormal frame of
$\bigwedge^{1,0}(T^{\lor})$.
Then we have the naturally induced orthonormal frame
$\bar{\vece}$ of $\bigwedge^{0,1}(T^{\lor})$.
Recall that $(1,1)$-form $\sum\varphi_{i\,j}\cdot e_i\wedge\bar{e}_j$
is called primitive if $\sum_i \varphi_{i\,i}=0$.

\begin{lem}[Corlette]
$K$ is harmonic metric
if and only if
the $(1,1)$-form $\delbar\theta$ is primitive.
\end{lem}
\pf
See the page 376 in \cite{corlette}.
We give only a remark on the difference of the notation.
Our $\vartheta$, $\theta$, and $\del$ are denoted by
$\theta$, $\theta^{1,0}$ and $\del^+$ respectively in \cite{corlette}.
\hfill\qed

\vspace{.1in}
Let $\omega$ be the Kahler form of the Kahler manifold $(X,g)$.
Let $e_i$ and $\bar{e}_i$ be as above.
\begin{lem}
There exists a negative constant $C_0$ depending only on $\dim X$
such that the following holds, for any $i\neq j$:
\[
 e_i\wedge \bar{e}_i\wedge e_j\wedge \bar{e}_j
\wedge \omega^{n-2}=C_0\cdot \omega^n.
\]
\end{lem}
\pf
Recall that $\omega$ is described as
$2^{-1}\cdot\sqrt{-1}\sum e_i\wedge \bar{e}_i$.
Then the claim immediately follows.
\hfill\qed

\vspace{.1in}

Let $h$ be a harmonic metric of a flat bundle $(E,\nabla)$.
Then we have the induced hermitian metric 
$\langle\cdot,\cdot\rangle$ of $End(E)$.
It naturally induces the following pairing,
which we also denote by $\langle\cdot,\cdot\rangle$:
\[
 \bigl(
 \End(E)\otimes \Omega^{\cdot,\cdot}_X
 \bigr)
\otimes
\overline{
 \bigl(
 \End(E)\otimes\Omega^{\cdot,\cdot}_X
 \bigr)}
\lrarr
 \Omega^{\cdot,\cdot}_X.
\]
For example, we have the following:
\[
 \big\langle f\cdot dz^I\wedge d\bar{z}^J,\,
 g\cdot dz^K\wedge d\bar{z}^L\big\rangle
=\langle f,g\rangle\cdot
 dz^I\wedge d\bar{z}^J\wedge d\bar{z}^K\wedge dz^L.
\]

\begin{lem}
We have the following formula:
\begin{equation}\label{eq;a12.28.30}
 \del\delbar\big\langle\theta,\,\theta\big\rangle
=-\frac{1}{2}\big\langle[\theta,\theta],[\theta,\theta] \big\rangle
+\big\langle \delbar\theta,\delbar\theta\big\rangle.
\end{equation}
\end{lem}
\pf
We have the following:
\[
 \del\delbar\big\langle\theta,\theta\big\rangle
=\del\big\langle \delbar\theta,\theta\big\rangle
=\big\langle\del\delbar\theta,\theta \big\rangle
+\big\langle\delbar\theta,\delbar\theta\big\rangle.
\]
Due to the vanishings
$\del\delbar+\delbar\del+[\theta,\theta^{\dagger}]=0$
and $\del\theta=0$,
we have the following:
\[
 \del\delbar\theta=
-\delbar\del\theta
-\bigl[[\theta,\theta^{\dagger}],\theta\bigr]
=\frac{1}{2}\bigl[[\theta,\theta],\theta^{\dagger}\bigr]
=-\frac{1}{2}\bigl[\theta^{\dagger},[\theta,\theta]\bigr].
\]
Since we have
$\big\langle\big[\theta^{\dagger},[\theta,\theta]\big],\,
 \theta \big\rangle
=\big\langle[\theta,\theta],\,[\theta,\theta]\big\rangle$,
we obtain (\ref{eq;a12.28.30}).
\hfill\qed

\vspace{.1in}
Thus we obtain the equality
$\del\delbar\big\langle\theta,\,\theta\big\rangle\cdot\omega^{n-2}
=-\frac{1}{2}\big\langle[\theta,\theta],\,[\theta,\theta]\big\rangle
 \cdot\omega^{n-2}
+\big\langle\delbar\theta,\,\delbar\theta\big\rangle\cdot\omega^{n-2}$.

We have the description $\theta=\sum f_i\cdot e_i$
for $f_i\in\End(E)$.
Then we have the following:
\[
 [\theta,\theta]
=2\sum_{i<j}[f_i,f_j]\cdot e_i\wedge e_j.
\]

\begin{lem}\label{lem;a12.28.35}
We have the following formula:
\begin{equation}\label{eq;a12.28.31}
 -\frac{1}{2}\bigl\langle[\theta,\theta],\,[\theta,\theta]\bigr\rangle
 \cdot\omega^{n-2}
=\sum_{i,j}\big|[f_i,f_j]\big|^2_h\cdot C_0\cdot\omega^{n}.
\end{equation}
\end{lem}
\pf
In the case $i<j$ and $k<l$,
we have the following equality:
\[
 \big\langle
 [f_i,f_j],\,[f_k,f_l]
 \big\rangle
 \cdot e_i\wedge e_j\wedge \bar{e}_k\wedge \bar{e}_l\wedge\omega^{n-2}
=\left\{
 \begin{array}{ll}
 -\bigl|
 [f_i,f_j]
 \bigr|^2_h\cdot C_0\cdot \omega^n & (i=k,j=l,i\neq j),\\
 \mbox{{}}\\
 0 & (\mbox{otherwise}).
 \end{array}
 \right.
\]
Then we immediately obtain (\ref{eq;a12.28.31}).
\hfill\qed

\vspace{.1in}
We have the description
$\delbar\theta=\sum\varphi_{i,j}\cdot e_i\wedge \bar{e}_j$.
\begin{lem}\label{lem;a12.28.36}
We have the following formula:
\begin{equation}\label{eq;a12.28.33}
 \bigl\langle\delbar\theta,\,\delbar\theta\bigr\rangle
 \cdot\omega^{n-2}
=C_0\cdot\left(
 \sum_{i,j}\bigl\langle\varphi_{i,j},\varphi_{i,j}\bigr\rangle
 \right)\cdot\omega^{n}.
\end{equation}
\end{lem}
\pf
We have the following:
\begin{multline} \label{eq;a12.28.32}
\bigl\langle
\delbar\theta,\,\delbar\theta
\bigr\rangle\cdot\omega^{n-2}
=\sum_{i,j,k,l}
 \bigl\langle\varphi_{i\,j},\,\varphi_{k,l}\bigr\rangle\cdot
 e_i\wedge \bar{e}_j\wedge \bar{e}_k\wedge e_l\wedge\omega^{n-2}\\
=\sum_{i\neq j}
 \bigl\langle\varphi_{i,j},\,\varphi_{i,j}\bigr\rangle
 e_i\wedge \bar{e}_i\wedge e_j\wedge \bar{e}_j\wedge\omega^{n-2}
-\sum_{i\neq k}
 \bigl\langle \varphi_{i,i},\,\varphi_{k,k}\bigr\rangle
 e_i\wedge \bar{e}_i\wedge e_k\wedge \bar{e}_k\wedge\omega^{n-2}\\
=C_0\cdot
 \Bigl(
 \sum_{i\neq j}\bigl\langle \varphi_{i,j},\,\varphi_{i,j}\bigr\rangle
-\sum_{i\neq k}\bigl\langle \varphi_{i,i},\varphi_{k,k}\bigr\rangle
\Bigr)\cdot\omega^n.
\end{multline}
Recall we have $\sum_{i}\varphi_{i,i}=0$,
for $\delbar\theta$ is primitive.
Therefore we have the following equality:
\[
 \sum_i\langle\varphi_{i,i},\varphi_{i,i}\rangle
+\sum_{i\neq k}\langle\varphi_{i,i},\varphi_{k,k}\rangle=0.
\]
Then we immediately obtain (\ref{eq;a12.28.33}).
\hfill\qed

\begin{prop}\label{prop;04.1.11.1}
There exist negative constants $C_1$ and $C_2$
depending only on $\dim X$
such that the following holds:
\begin{equation}\label{eq;a12.28.37}
 \del\delbar\big\langle\theta,\,\theta\big\rangle\cdot \omega^{n-2}
=\Bigl(
 C_1\cdot\big|[\theta,\theta]\big|^2_h
+C_2\cdot\big|\delbar\theta\big|_h^2
 \Bigr)\cdot\omega^n.
\end{equation}
\end{prop}
\pf
It immediately follows from Lemma \ref{lem;a12.28.35}
and Lemma \ref{lem;a12.28.36}.
\hfill\qed

\vspace{.1in}
The formula (\ref{eq;a12.28.37}) is the Bochner type formula
due to Corlette.
Recall the argument to derive
the pluri-harmonicity of the harmonic metric $h$,
when $X$ is compact Kahler.
If $X$ is compact,
we have the vanishing:
\[
0=\int_X\del\delbar\big\langle\theta,\theta\big\rangle\cdot\omega^{n-2}
=C_1\int_X\big|[\theta,\theta]\big|^2
+C_2\int_X\big|\delbar\theta\big|^2.
\]
It implies $[\theta,\theta]=\delbar\theta=0$ on $X$,
which means the metric $h$ is pluri-harmonic.

\subsubsection{A variation of the Bochner type formula}
\label{subsubsection;04.1.30.350}

Let $f$ and $g$ be sections of $\End(E)$.
We denote the adjoint of them by $f^{\dagger}$ and $g^{\dagger}$
respectively.
\begin{lem}
We have the equality:
$\bigl\langle f\cdot dz_i,\,\,g\cdot dz_j
 \bigr\rangle
=-\bigl\langle g^{\dagger}\cdot d\bar{z}_j,\,\,
 f^{\dagger}\cdot d\bar{z}_i
 \bigr\rangle$.
\end{lem}
\pf
It follows from the equality
$\langle f,g\rangle=\langle g^{\dagger},f^{\dagger}\rangle$.
\hfill\qed

\begin{cor}\label{cor;04.1.11.3}
We have the equality
$\del\delbar\bigl\langle\theta,\,\,\theta\bigr\rangle
=\delbar\del\bigl\langle\theta^{\dagger},\,\,\theta^{\dagger}\bigr\rangle$.
\hfill\qed
\end{cor}

For our argument to derive the pluri-harmonicity
in quasi projective case,
we will also use the following formula.

\begin{prop}\label{prop;04.1.22.40}
We have the following formula:
\begin{equation}\label{eq;04.1.11.4}
 d\Bigl(
 \bigl\langle
 \delbar\theta,\,\theta-\theta^{\dagger}
 \bigr\rangle
\wedge \omega^{n-2}
 \Bigr)
=2\Bigl(
C_1\cdot\big|[\theta,\theta]\big|^2_h
+C_2\cdot\big|\delbar\theta\big|_h^2
 \Bigr)\cdot\omega^{n}.
\end{equation}
Here the negative constants $C_i$ are 
same as those in Proposition {\rm\ref{prop;04.1.11.1}}.
\end{prop}
\pf
We have the following:
\[
 \del\delbar\bigl\langle\theta,\theta\bigr\rangle\wedge\omega^{n-2}
=\del\bigl\langle\delbar\theta,\theta\bigr\rangle\wedge\omega^{n-2}
=\del\bigl\langle\delbar\theta,
  \theta-\theta^{\dagger}\bigr\rangle\wedge\omega^{n-2}.
\]
Similarly, we have the equality
$\delbar\del\langle\theta^{\dagger},\theta^{\dagger}\bigr\rangle
 \wedge\omega^{n-2}
=\delbar\bigl\langle\del\theta^{\dagger},\theta^{\dagger}-\theta
 \bigr\rangle\wedge\omega^{n-2}$,
which can be rewritten as follows:
\[
 \delbar
 \bigl\langle
 -\delbar\theta,\,\,-\bigl(\theta-\theta^{\dagger}\bigr)
 \bigr\rangle\wedge\omega^{n-2}
=\delbar\bigl\langle
 \delbar\theta,\theta-\theta^{\dagger}\bigr\rangle
 \wedge \omega^{n-2}.
\]
Here we have used Corollary \ref{cor;04.1.11.2}.
Then we obtain
the following:
\[
 \del\delbar\bigl\langle\theta,\theta\bigr\rangle\wedge\omega^{n-2}
+\delbar\del\bigl\langle\theta^{\dagger},\theta^{\dagger}
 \bigr\rangle\wedge\omega^{n-2}
=d\Bigl(
 \bigl\langle\delbar\theta,\theta-\theta^{\dagger}
 \bigr\rangle\wedge\omega^{n-2}
 \Bigr).
\]
Then (\ref{eq;04.1.11.4}) follows from Corollary \ref{cor;04.1.11.3}
and Proposition \ref{prop;04.1.11.1}.
\hfill\qed

\section{Tame pure imaginary harmonic bundle}

\label{section;04.1.30.205}

\subsection{Definition}
\label{subsection;04.1.30.206}
%c13.21.tex

Let $X$ be a complex manifold,
and $D=\bigcup_i D_i$ be a normal crossing divisor.
Let $\harmonicbundle$ be a harmonic bundle on $X-D$.
Let $P$ be a point of $X$.
Let us take a neighbourhood $P$ with a coordinate
$(z_1,\ldots,z_n)$ such that
$U\cap D=\bigcup_{i=1}^l\{z_i=0\}$.
We have the description
$\theta=\sum_{i=1}^l f_i\cdot dz_i/z_i+\sum_{j=l+1}^n g_j\cdot dz_j$.
Recall that 
$\harmonicbundle$ is called tame
if the coefficients of the characteristic polynomials
$\det(t-f_i)$ and $\det(t-g_j)$ are holomorphic.

\begin{lem} \label{lem;04.1.23.30}
The harmonic bundle $\harmonicbundle$ is tame
if and only if
there exists a holomorphic vector bundle
$\tilde{E}$ with a regular Higgs field
$\tilde{\theta}\in\End(\tilde{E})\otimes\Omega^{1,0}_X(\log D)$
such that
$(\tilde{E},\tilde{\theta})_{|X-D}=(E,\theta)$.
\end{lem}
\pf
In the subsection 8.6 in \cite{mochi2},
it is proved that
the prolongment $\prolong{E}$ by an increasing order
is locally free and that $\theta$ naturally induces
the regular Higgs field on $\prolong{E}$.
On the contrary,
it is easy to check $\harmonicbundle$ is tame
if there exists a prolongment $(\tilde{E},\tilde{\theta})$.
\hfill\qed

\vspace{.1in}
Let $\harmonicbundle$ be a tame harmonic bundle on $X-D$.
Let $(\tilde{E},\tilde{\theta})$ be any prolongment
of $(E,\theta)$.
Let $\Res_i(\tilde{\theta})$ denote the residue of $\tilde{\theta}$
with respect to the irreducible component $D_i$ of $D$.

\begin{lem}
Let $P$ be a point of $D_i$.
The eigenvalues 
$\Res_i(\tilde{\theta})_{|P}$ is independent of
a choice of a prolongment $(\tilde{E},\tilde{\theta})$.
\end{lem}
\pf
The eigenvalues are solutions of $\det(t-f_i)_{|P}$.
Since $\det(t-f_i)$ is determined independently of
$(\tilde{E},\tilde{\theta})$,
we are done.
\hfill\qed

\vspace{.1in}
Recall that the eigenvalues of $\Res_i(\theta)_{|P}$
is independent of a choice of $P\in D_i$,
which is proved in the subsection 8.1 in \cite{mochi2}.

\begin{df}
If any eigenvalues of the residues $\Res_i(\tilde{\theta})$
is pure imaginary, $\harmonicbundle$ is called
a tame pure imaginary harmonic bundle.
\hfill\qed
\end{df}

\subsection{Tame pure imaginary harmonc bundle on a punctured disc}
\label{subsection;04.1.30.207}
%c13.3.tex

\subsubsection{The estimate of Higgs field}

We use the Poincar\'{e} metric
$g_0:=|z|^{-2}\cdot(-\log|z|^2+A)^{-2}\cdot dz\cdot d\bar{z}$
on $\Deltabarast$.
Here $A$ denotes a positive number.
Let us consider a tame harmonic bundle
$\harmonicbundle$ on a punctured disc $\Delta^{\ast}$.
We have the prolongment $\prolong{E}$
and the description $\theta=f_0\cdot dz/z$
for some $f_0\in \End(\prolong{E})$ on $\Delta$.
We put $\Sp(\theta):=\Sp(f_{0\,|\,O})$.
For any $\alpha\in\Sp(\theta)$,
we denote $\dim\EE(\prolong{E},\alpha)$
by $\multiplicity(\alpha)$.
We put as follows:
\[
 \gminit(\theta)=
 \sum_{\alpha\in\Sp(\theta)}\gminim(\alpha)\cdot|\alpha|^2.
\]

The following lemma can be elementarily shown.
\begin{lem}\label{lem;a12.28.40}
There exist positive constants $R$, $C$ and $\epsilon$,
depending on $(E,\delbar_E,\theta)$,
such that the following holds:
\begin{itemize}
\item
 For any point $P\in\Delta^{\ast}(R)$
 and for any eigenvalue of $f_{0\,|\,P}$,
 there exists the unique element $\alpha\in\Sp(\theta)$
 such that the following holds:
\[
 |a-\alpha|\leq C\cdot |z(P)|^{\epsilon}.
\]
\end{itemize}
\hfill\qed
\end{lem}

\begin{lem}\label{lem;a12.28.42}
There exist positive constants $R',C',\epsilon'$,
depending only on $R,C,\epsilon$ in Lemma {\rm\ref{lem;a12.28.40}},
such that the following claims hold.
\begin{itemize}
\item
 We have the decomposition
 $E=\bigoplus_{\alpha\in\Sp(\theta)}E_{\alpha}$
 on $\Delta^{\ast}(R')$.
 The decomposition is preserved by $f_0$.
 For any point $P\in \Delta^{\ast}(R')$,
 and for any $v\in E_{\alpha_i\,|\,P}$ and $w\in E_{\alpha_j\,|\,P}$
 $(\alpha_i\neq \alpha_j)$,
 the following inequality holds:
\[
 h(v,w)\leq C'\cdot |z(P)|^{\epsilon'}\cdot
 |v|_h\cdot|w|_h.
\]
\item
 The following inequality holds:
\begin{equation}\label{eq;a12.28.41}
 \Bigl|
 |f_0|_h^2-\gminit(\theta)
 \Bigr|
\leq C'\cdot (-\log|z|^2+A)^{-2}.
\end{equation}
\end{itemize}
\end{lem}
\pf
The claims are proved in the subsection 7.1 of \cite{mochi2}.
\hfill\qed

\begin{cor}\label{cor;a12.29.57}
There exists constant $C_0$,
depending only on $R,C,\epsilon$ in Lemma {\rm\ref{lem;a12.28.40}},
such that the following holds:
\begin{equation}\label{eq;04.2.2.1}
 \Bigl|
 |\theta|^2_{h,g_0}
-2\gminit(\theta)\cdot \bigl(-\log|z|^2\bigr)^2
 \Bigr|
\leq C_0.
\end{equation}
\end{cor}
\pf
It immediately follows from Lemma \ref{lem;a12.28.42}.
\hfill\qed

\vspace{.1in}
Lemma \ref{lem;a12.28.40} and Lemma \ref{lem;a12.28.42}
hold for any tame harmonic bundle on $\Delta^{\ast}$.
In the pure imaginary case,
the estimate (\ref{eq;04.2.2.1}) can be refined.
We put $\vartheta=\theta+\theta^{\dagger}$.
We use the polar coordinate $z=r\cdot\exp(\sqrt{-1}\eta)$.
We denote $\del/\del r$ and $\del/\del\eta$
by $\del_r$ and $\del_{\eta}$ respectively.

\begin{lem}\label{lem;a12.29.55}
Assume that $\harmonicbundle$ is pure imaginary.
There exists a positive number $C''$,
depending only on $R,C,\epsilon$ in Lemma {\rm\ref{lem;a12.28.40}},
such that the following inequalities hold:
\[
\Bigl|
 \bigl|
 \vartheta\bigl(\del_{\eta}\bigr)
 \bigr|_h^2
-4\gminit(\theta)
\Bigr|
\leq C''\cdot \bigl(-\log r^2+A\bigr)^{-2}.
\]
\[
 \bigl|
 \vartheta\bigl(\del_{r}\bigr)
 \bigr|_h^2
\leq C''\cdot r^{-2}\cdot\bigl(-\log r^2+A\bigr)^{-2}.
\]
\end{lem}
\pf
We have the following description:
\[
 \vartheta=
 \bigl(f+f^{\dagger}\bigr)\frac{dr}{r}
+\sqrt{-1}\bigl(f-f^{\dagger}\bigr)d\eta.
\]
Then the claims follow from 
Lemma \ref{lem;a12.28.42}.
\hfill\qed

%c13.31.tex

\subsubsection{The estimate of the energy of the associated twisted
   harmonic map}
Let $\harmonicbundle$ be a tame pure imaginary harmonic bundle.
Let $(\nbige^1,\DD^1)$ be the flat bundle
associated to $\harmonicbundle$ on $\Deltabarast$.
\begin{lem}\label{lem;04.1.6.3}
The $KMS$-spectrum of $\nbige^1$ is given as follows:
\[
\begin{array}{l}
 \KMS(\nbige^1)=
 \bigl\{
 (b,2\sqrt{-1}c-b)\,\big|\,
 (b,\sqrt{-1}c)\in \KMS(\nbige^0)
 \bigr\},\\
\mbox{{}}\\
 \KMS^f(\nbige^1)
=\bigl\{
 \bigl(0,\exp(2\pi\sqrt{-1}b+4\pi c)\bigr)\,\big|\,
 (b,\sqrt{-1}c)\in\KMSoverline(\nbige^0)
 \bigr\}.
\end{array}
\]
\end{lem}
\pf
See the the section 5 in \cite{s2} or
the subsection 7.3 and 7.4 in \cite{mochi2}.
We only remark the following equalities:
\[
 \paramap(1,b,\sqrt{-1}c)=b,
\quad
 \eigenmap(1,b,\sqrt{-1}c)=2\sqrt{-1}c-b,
\]
\[
 \paramap^f(1,b,\sqrt{-1}c)
=\Re\bigl(1\cdot\sqrt{-1}c+1\cdot \overline{\sqrt{-1}c}\bigr)=0,
\quad\quad
 \eigenmap^f(1,b,\sqrt{-1}c)
=\exp\bigl(
 2\pi\sqrt{-1}b+4\pi c
 \bigr).
\]
\hfill\qed

\vspace{.1in}

Let $\varphi$ be the monodromy of
the flat bundle $(\nbige^1,\DD^1)$.
\begin{lem}
We have 
$ \rho(\varphi)^2=64\pi^2\cdot\gminit(\theta)$.
(See the subsubsection {\rm\ref{subsubsection;04.1.6.2}}
 for $\rho$.)
\end{lem}
\pf
From Lemma \ref{lem;04.1.6.3},
we obtain the following:
\[
 \rho(\varphi)^2
=\sum_{(b,\sqrt{-1}c)\in \KMS(\nbige^0)}
 \multiplicity(b,\sqrt{-1}c)\cdot(2\cdot 4\pi c_i)^2
=64\pi^2\sum_{\sqrt{-c}\in\Sp(\theta)}
 \multiplicity(\sqrt{-1}c)\cdot c^2
=64\pi^2\cdot\gminit(\theta).
\]
Thus we are done.
\hfill\qed

\vspace{.1in}
Let $\Psi_h$ be the twisted harmonic map
$\Delta^{\ast}\lrarr \PH(r)\big/\langle\varphi\rangle$
associated with $(\nbige^1,\DD^1,h)$.
Recall that we have $e(\Psi_h)=8\cdot|\theta|^2_{h,g_0}$
due to (\ref{eq;a12.28.45}).

\begin{lem}\label{lem;a12.29.56} \label{lem;a12.29.1}
There exists a positive constant $C_1$,
depending only on $R,C,\epsilon$ in Lemma {\rm\ref{lem;a12.28.40}},
such that the following holds:
\[
\Bigl|
 e(\Psi_h)-\frac{\rho(\varphi)^2}{4\pi^2}
 \cdot\bigl(-\log|z|^2+A\bigr)^2
\Bigr|
\leq C_1.
\]
As a result, we obtain the following finiteness:
\[
 \int_{\Deltabarast}
\left|
 e(\Psi_h)
-\frac{\rho(\varphi)^2}{4\pi^2}
 \cdot\bigl(-\log|z|^2+A\bigr)^2
\right|\cdot
 \dvol_{g_0}
<\infty.
\]
\end{lem}
\pf
It immediately follows from
$e(\Psi_h)=8\cdot|\theta|^2_{h,g_0}$
and the estimate of $\theta$.
\hfill\qed

%c13.32.tex

\subsubsection{A characterization of tame pure imaginary harmonic bundle
on a punctured disc}
\label{subsubsection;04.1.29.300}

We put
$T(R_1,R_2):=
 \bigl\{z\in\cnum\, \big|\, R_1\leq -\log|z|\leq R_2
 \bigr\}$,
and $T(R):=T(0,R)$.
We use the Poincar\'{e} metric
$g:=|z|^{-2}\cdot \bigl(-\log|z|^2+A\bigr)^{-2}dz\cdot d\bar{z}$
on $\Deltabarast$.
We use the real coordinate
$z=\exp\bigl(\sqrt{-1}x-y\bigr)$.

\begin{lem}\label{lem;04.1.25.2}
Let $(E,\nabla)$ be a flat bundle on $T(R_1,R_2)$
with the monodromy $\varphi$.
Let $h$ be a hermitian metric of $(E,\nabla)$,
and let $\Psi_h:T(R_1,R_2)\lrarr \PH(r)\big/\langle\varphi\rangle$
be the corresponding twisted map.
Then we have the following a priori lower bound
of the energy:
\begin{equation}\label{eq;04.1.25.1}
 \int_{T(R_1,R_2)}
 e\bigl(\Psi_h\bigr)\dvol_g
\geq
 \int_{T(R_1,R_2)}
 \bigl|\del_x\Psi_h\bigr|^2\cdot\bigl(2y+A\bigr)^2 \dvol_g
\geq
 \int_{T(R_1,R_2)}
 \frac{\rho(\varphi)^2}{4\pi^2}\bigl(2y+A\bigr)^2
 \cdot\dvol_g
=\frac{\rho(\varphi)^2}{2\pi}(R_2-R_1).
\end{equation}
\end{lem}
\pf
We always have the inequality
$e(\Psi_h)\geq
\bigl|\del_x\Psi_h\bigr|^2\cdot\bigl|\del_x\bigr|^{-2}
=\bigl|\del_x\Psi_h\bigr|^2\cdot(2y+A)^2$.
We have the following inequality,
due to Lemma \ref{lem;a12.29.10}:
\[
\int_{0}^{2\pi}
 \left|\frac{\del \Psi_h}{\del x}\right|^2\cdot
 dx
\geq
 \frac{\rho(\varphi)^2}{2\pi}.
\]
Then the inequalities (\ref{eq;04.1.25.1}) immediately follows.
\hfill\qed

\begin{prop}\label{prop;04.1.25.10}
Let $\harmonicbundle$ be a harmonic bundle
on a punctured disc $\Deltabarast$.
Let $\varphi$ be the monodromy of 
the corresponding flat bundle $\bigl(\nbige^1,\DD^1\bigr)$,
and let $\Psi_h:\Deltabarast\lrarr\PH(r)\big/\langle \varphi\rangle$
denote the corresponding twisted harmonic map.
Assume that there exists an integrable function $F$
on $\Deltabarast$ with respect to the measure $\dvol_g$
satisfying the following,
for any sufficiently large $R$:
\begin{equation}\label{eq;04.1.25.3}
 \int_{T(R)}
 e\bigl(\Psi_h\bigr)\dvol_g
\leq
 \int_{T(R)}
 \left(\frac{\rho(\varphi)^2}{4\pi^2}(2y+A)^2+F\right)\dvol_g.
\end{equation}
Then $\harmonicbundle$ is tame and pure imaginary.
\end{prop}
\pf
First we see that the harmonic bundle $\harmonicbundle$ is tame.
Due to Lemma \ref{lem;04.1.25.2},
we have the lower bound of the energy
for any $R$:
\begin{equation}\label{eq;04.1.25.4}
\int_{T(R)}
 e\bigl(\Psi_h\bigr)\cdot\dvol_g
\geq
 \int_{T(R)}\frac{\rho(\varphi)^2}{4\pi^2}
\cdot\bigl(2y+A\bigr)^2\cdot\dvol_g.
\end{equation}
From (\ref{eq;04.1.25.3}) and (\ref{eq;04.1.25.4}),
we obtain the following inequality:
\begin{equation}\label{eq;04.1.25.5}
\int_{T(R_1,R_2)}
 e\bigl(\Psi_h\bigr)\cdot\dvol_g
\leq
 \int_{T(R_1,R_2)}
 \frac{\rho(\varphi)^2}{4\pi^2}\cdot\bigl(2y+A\bigr)^2
\cdot\dvol_g
+\int_{T(0,R_2)}F\cdot\dvol_g
\leq
 \frac{\rho(\varphi)^2}{2\pi}\cdot(R_2-R_1)
+C.
\end{equation}
Here we put $C=\int_{\Deltabarast}F\cdot\dvol$.

Let us consider the universal covering
$\hyperh\lrarr \Delta^{\ast}$
given by $x+\sqrt{-1}y\longmapsto z=\exp\bigl(\sqrt{-1}x-y\bigr)$.
We have the induced Poincar\'{e} metric
$g=(2y+A)^{-2}\cdot\bigl(dx\cdot dx+dy\cdot dy\bigr)$
on $\hyperh$.
We have the induced map
$\Psi_h:\hyperh\lrarr \PH(r)$.
We put as follows:
\[
 \tilde{S}(x_0,y_0):=
 \bigl\{
 x+\sqrt{-1}y\,
 \big|\,y_0-1\leq y\leq y_0+1,\,\,x_0-\pi\leq x\leq x_0+\pi
 \bigr\}.
\]
From (\ref{eq;04.1.25.5}),
we have the following inequality,
for the energy of $\Psi_h$:
\begin{equation}\label{eq;04.1.25.6}
 \int_{\tilde{S}(x_0,y_0)}
 e(\Psi_h)\cdot\dvol_g
\leq
 \frac{\rho(\varphi)^2}{\pi}+C.
\end{equation}

Let $g_1:=dx\cdot dx+dy\cdot dy$ be the Euclidean metric.
Let $e_{g_1}(\Psi_h)$ denote the energy function of $\Psi_h$
with respect to the metric $g_1$.
From (\ref{eq;04.1.25.6}),
we have the following inequality:
\begin{equation}\label{eq;04.1.25.7}
 \int_{\tilde{S}(x_0,y_0)}
 e_{g_1}(\Psi_h)\cdot\dvol_{g_1}
\leq
 \frac{\rho(\varphi)^2}{\pi}+C.
\end{equation}

\begin{lem}
We have the estimate
$e(\Psi_h)=O\bigl((2y+A)^2\bigr)$.
\end{lem}
\pf
Since the right hand side of (\ref{eq;04.1.25.7})
is independent of a choice of $(x_0,y_0)$,
we obtain the constant $C_1$ such that
$e_{g_1}\bigl(\Psi_h\bigr)\leq C_1$ on $\hyperh$,
due to Lemma \ref{lem;a12.29.4}.
Since we have the relation
$e(\Psi_h)=e_{g_1}(\Psi_h)\cdot \bigl(2y+A\bigr)^2$,
we obtain the estimate
$e\bigl(\Psi_h\bigr)\leq C_1\cdot (2y+A)^2$.
\hfill\qed

\vspace{.1in}

Recall that we have the relation
$8\cdot\bigl|\theta\bigr|^2=e\bigl(\Psi_h\bigr)$.
Let us describe $\theta=f\cdot dz/z$,
and then we have
$\bigl|\theta\bigr|^2=
2\cdot |f|_h^2\cdot (2y+A)^2$.
Hence we obtain the boundedness of $|f|_h$ on $\Deltabarast$.
Then it is easy to derive that 
the coefficients of $\det(t-f)$ are holomorphic 
on $\Deltabar$,
namely, the harmonic bundle $\harmonicbundle$ is tame.

\vspace{.1in}
Let us show that the harmonic bundle $\harmonicbundle$ is pure imaginary.

\begin{lem}
$\bigl|\del_y\Psi_h\bigr|^2\cdot (2y+A)^2$ is integrable
on $\Deltabarast$
with respect to the measure $\dvol_g$.
\end{lem}
\pf
From (\ref{eq;04.1.25.1}),
we have the following inequality:
\begin{equation}\label{eq;04.1.25.8}
\int_{T(R)}
 \bigl|\del_x\Psi_h\bigr|^2\cdot (2y+A)^2\cdot\dvol_g
\geq
\int_{T(R)}\frac{\rho(\varphi)^2}{4\pi^2}\cdot
 \bigl(2y+A\bigr)^2 \cdot\dvol_g.
\end{equation}
From (\ref{eq;04.1.25.3}) and (\ref{eq;04.1.25.8}),
we obtain the following inequality for any $R$:
\[
 \int_{T(R)}
 \bigl|\del_y\Psi_h\bigr|^2\cdot (2y+A)^2\cdot\dvol_g
\leq
 \int_{T(R)}F\cdot\dvol_g.
\]
It implies the integrability of
$\bigl|\del_y\Psi_h\bigr|^2\cdot(2y+A)^2$.
\hfill\qed

\vspace{.1in}
We have $\theta=f\cdot dz/z=f\cdot \bigl(\sqrt{-1}dx-dy\bigr)$.
Let $f^{\dagger}$ denote the adjoint of $f$,
and then we have
$\theta^{\dagger}=f^{\dagger}\cdot\bigl(-\sqrt{-1}dx-dy\bigr)$.
We have the following equalities:
\[
 \bigl|\del_{y}\Psi_h\bigr|^2\cdot(2y+A)^2
=4\cdot\bigl|
 \theta(\del_y)+\theta^{\dagger}(\del_y)
 \bigr|^2\cdot (2y+A)^2
=4\cdot \bigl|f+f^{\dagger}\bigr|^2\cdot(2y+A)^2.
\]
Since we have already known that the harmonic bundle
$\harmonicbundle$ is tame,
we have the following estimate, due to Lemma \ref{lem;a12.28.42}:
\begin{equation}\label{eq;a12.29.21}
 \bigl|
 f+f^{\dagger}
 \bigr|^2\cdot (2y+A)^2
=\sum_{a\in \Sp(\theta)}
 \bigl| 2\Re(a)
 \bigr|^2\cdot (2y+A)^2
+O\bigl(1\bigr).
\end{equation}
Then we obtain the following vanishing,
from the integrability of $\bigl|\del_y\Psi_h\bigr|^2\cdot(2y+A)^2$
on $\Deltabarast$
with respect to the measure $\dvol_g$:
\[
 \sum_{a\in \Sp(\theta)}
 \bigl| 2\Re(a)
 \bigr|^2=0.
\]
Namely the harmonic bundle $\harmonicbundle$ is pure imaginary.
Therefore the proof of Proposition \ref{prop;04.1.25.10} is accomplished.
\hfill\qed

\subsection{Semisimplicity}
\label{subsection;04.1.30.208}
%c13.23.tex

\subsubsection{Statement}

\begin{prop} \label{prop;04.1.2.2}
Assume that $X$ is a smooth projective variety, and $D$ is a normal
crossing divisor of $X$.
Let $\harmonicbundle$ be a tame pure imaginary harmonic bundle on $X-D$.
Then the corresponding flat bundle
$(\nbige^1,\DD^1)$ is semisimple.
\hfill\qed
\end{prop}
We will prove the proposition in the next subsubsections
\ref{subsubsection;04.1.5.1}--\ref{subsubsection;04.1.5.3}.
We will also prove the reverse of the proposition
in the next sections.

\subsubsection{Stability and semistability}
\label{subsubsection;04.1.5.1}

Let $C$ be a smooth quasi projective curve over $\cnum$,
and $\overline{C}$ be a smooth projective completion.
We put $D=\overline{C}-C=\{P_1,\ldots,P_l\}$.
Let $E$ be a holomorphic bundle over $\overline{C}$.
Let us consider a neighbourhood $U$ of $P$
with a coordinate $z$ such that $z(P_i)=0$.

\begin{lem}\label{lem;04.1.5.4}
The following data are equivalent:
\begin{itemize}
\item
 A filtration $F'$ of $E_{|P_i}$ indexed by $\openclosed{-1}{0}$.
\item
 A filtration $F$ of $\bigcup_{h}E(h\cdot P_i)$ indexed by $\real$
 such that $F_0=E$ and $F_a\cdot z^{-1}=F_{a+1}$.
\end{itemize}
\end{lem}
\pf
Let $F'$ be a filtration of $E_{|P}$ indexed by
$\openclosed{-1}{0}$.
We put
$F_a:=\bigl\{f\in E\,|\,f_{|P}\in F'_a\bigr\}$
for any number $a\in\openclosed{-1}{0}$.
For any real number $a\in\real$,
we take the number $a_0\in\openclosed{-1}{0}$
and the integer $a_1$ which are uniquely determined
by $a_0+a_1=a$.
Then we put $F_a:=F_{a_0}\cdot z^{-a_1}$.
Thus we obtain the filtration of $\bigcup E(h\cdot P_i)$
satisfying the condition.

The claim in the reverse direction can be shown similarly.
\hfill\qed

\vspace{.1in}
We will not distinguish two kind of data
in Lemma \ref{lem;04.1.5.4}.
They are called the parabolic structure of $E$.

Let $F$ be a filtration of $\bigcup_hE(h\cdot D)$ as above,
i.e.,
we are given filtrations $F'$ of $E_{|P_i}$ $(P_i\in D)$.
For any $a\in\openclosed{-1}{0}$,
we put $\multiplicity(a)=\sum_{P_i}\dim\Gr^{F'}_a(E_{|P_i})$.
We put as follows:
\begin{equation}\label{eq;04.1.5.5}
 \deg(E,F):=\deg(E)-\sum_{a\in\openclosed{-1}{0}} a\cdot
 \multiplicity(a),
\quad\quad
 \mu(E,F):=\frac{\deg(E,F)}{\rank E}.
\end{equation}
\begin{rem}
The formula {\rm(\ref{eq;04.1.5.5})} is same as
that given in the section {\rm 6} of {\rm\cite{s2}}.
Note our parabolic filtration is increasing.
\hfill\qed
\end{rem}

A connection of $E_{|C}$ is called regular,
if we have $\nabla f\in F_a(E)\otimes\Omega^{1,0}(\log D)$
for any $f\in F_a(E)$.

Let $E'\subset E$ be a subsheaf,
then $F$ induces the filtration of $E'$
by $F_a(E')=F_a(E)\cap E'$.
We denote it also by $F$.
Recall that the filtered regular connection
$(E,F,\nabla)$ is called stable,
if the inequality
$\mu(E',F)<\mu(E,F)$ holds for any sub-connection $(E',F,\nabla)$.
We recall only a part of
the Kobayashi-Hitchin correspondence for harmonic metric
(Theorem 5 in \cite{s2}):

\begin{prop}\mbox{{}}\label{prop;04.1.2.6}
 Let $(E,\nabla,h)$ be a tame harmonic bundle on $C$.
 We obtain the filtration $F$ of $E$ by an increasing order.
 Then the regular filtered connection $(E,\nabla,F)$
 is polystable.
\hfill\qed
\end{prop}

\subsubsection{Quasi canonical prolongment and the canonical filtration}
\label{subsubsection;04.1.5.2}

Let $X$ be a complex manifold,
and $D$ be a normal crossing divisor of $X$.
Let $(E,\nabla)$ be a flat bundle on $X-D$.
Then we have the quasi canonical prolongment $QC(E)$ of $E$:
\begin{itemize}
\item $QC(E)$ is a holomorphic vector bundle on $X$.
\item For any $f\in QC(E)$,
 we have $\nabla(f)\in QC(E)\otimes\Omega^{1,0}(\log D)$.
\item
 Let $\alpha$ be any eigenvalue of $\Res(\nabla)$ for any 
 irreducible component of $D$.
Then it satisfies 
$0\leq \Re(\alpha)<1$.
\end{itemize}

For the quasi canonical filtration,
we have the canonical filtration 
of $QC(E)_{|D_i}$ for any irreducible component $D_i$ of $D$.
For simplicity, we only consider the case 
that $X=\overline{C}$ is a smooth projective curve.
We put $C=\overline{C}-D$.
Let $P_i$ be a point of $D$.
Then we have the generalized eigen decomposition
with respect to the residue $\Res_{P_i}(\nabla)$:
\[
 QC(E)_{|P_i}=\bigoplus_{\alpha\in\cnum}
 \EE\bigl(\Res_{P_i}(\nabla),\,\alpha\bigr).
\]
Here $\EE\bigl(\Res_{P_i}(\nabla),\alpha\bigr)$ denotes
$\Ker\bigl(\Res_{P_i}(\nabla)-\alpha\bigr)^N$ 
for any sufficiently large integer $N$.
Then we put as follows, for any $-1< a\leq 0$:
\[
 F_a\bigl(QC(E)_{|P_i}\bigr)
=\bigoplus_{-\Re(\alpha)\leq a}\EE(\Res_{P_i}(\nabla),\,\alpha).
\]
It also induces the filtration $F$ of
$\bigcup_h QC(E)(h\cdot D)$ (Lemma \ref{lem;04.1.5.4}).
The filtration is called the canonical filtration.
It is well known that $\mu(QC(E),F)=0$ holds.

\begin{lem}[Sabbah]\label{lem;04.1.2.1}
Let $(E,\nabla)$ be a flat connection on $C$.
Then $(E,\nabla)$ is simple if and only if
$(QC(E),F,\nabla)$ is stable.
\end{lem}
\pf
It is easy to see that the simplicity of $(E,\nabla)$
implies the stability of $(QC(E),F,\nabla)$.
Let us assume that $(E,\nabla)$ is not simple.
Then we have a sub-connection $(E',\nabla)\subset (E,\nabla)$.
Then we obtain the filtered subbundle
$\bigl(QC(E'),\nabla\bigr)\subset \bigl(QC(E),\nabla\bigr)$.
Since we have $\mu(QC(E'),F)=0=\mu(QC(E),F)$,
$\bigl(QC(E),F,\nabla\bigr)$ is not stable.
Thus we are done.
\hfill\qed

\begin{cor}\label{cor;04.1.2.7}
Let $X$ be a projective variety, and $D$ be a normal crossing divisor
of $X$. Let $(E,\nabla)$ be a flat connection on $X-D$.
Let $C$ be a smooth projective curve in $X$,
which is transversal with $D$
such that $\pi_1(C\setminus D)\lrarr \pi_1(X-D)$ is surjective
(see Lemma {\rm\ref{lem;04.1.22.10}}).

Assume that $\bigl(QC(E_{|C}),F,\nabla\bigr)$ is poly-stable,
then $(E,\nabla)$ is semisimple.
\end{cor}
\pf
Due to Lemma \ref{lem;04.1.2.1},
we know that $(E,\nabla)_{|C}$ is semisimple.
Since $\pi_1(C\setminus D)\lrarr\pi_1(X-D)$ is surjective,
we obtain that $(E,\nabla)$ is also semisimple.
\hfill\qed

\subsubsection{The end of the proof of Proposition \ref{prop;04.1.2.2}}
\label{subsubsection;04.1.5.3}

Let $\harmonicbundle$ be a tame  harmonic bundle
on a quasi projective curve $C$ with the completion $\overline{C}$.
We denote the corresponding flat bundle 
by $(\nbige^1,\DD^1)$.
Then we have the two kind of prolongment of $\nbige^1$
with the filtration.
One prolongation is $QC(\nbige^1)$ with the canonical filtration.
The other is the prolongment
$\prolong{\nbige^1}$ by an increasing order
with the filtration with respect to the metric $h$.

\begin{lem}\label{lem;04.1.2.5}
If $(E,\nabla,h)$ is pure imaginary,
we have the canonical isomorphism
$QC(\nbige^1)\simeq \prolong{\nbige^1}$
preserving the filtrations.
\end{lem}
\pf
It follows from Lemma \ref{lem;04.1.6.3}
and the uniqueness of $QC(\nbige^1)$ (\cite{d}).
\hfill\qed

\vspace{.1in}
Let us show Proposition \ref{prop;04.1.2.2}.
Let $\harmonicbundle$ be a tame pure imaginary harmonic bundle
on $X-D$.
It is well known that
we can take a smooth projective curve $C$ in $X$
such that $\pi_1(C\setminus D)\lrarr\pi_1(X-D)$ is surjective.
Let us consider the restriction 
$\harmonicbundle_{|C\setminus D}$.
Due to Lemma \ref{lem;04.1.2.5} and Proposition \ref{prop;04.1.2.6},
we know that $QC(E)$ with the canonical filtration
is polystable.
Thus we obtain the semisimplicity of $(E,\nabla)$
due to Corollary \ref{cor;04.1.2.7}.
Thus the proof of Proposition \ref{prop;04.1.2.2} is accomplished.
\hfill\qed

\subsection{The maximum principle}
\label{subsection;04.1.30.209}
%c13.4.tex

Let $X$ be a compact Riemannian surface
with the continuous boundary $\del X$.
Let $U$ denote the interior part of $X$.
Let us take $P_1,\ldots,P_l\in U$.
We put $X^{\ast}:=X-\{P_1,\ldots,P_l\}$
and $U^{\ast}:=U-\{P_1,\ldots,P_l\}$.
Let $(E,\nabla)$ be a flat bundle on $X^{\ast}$.
A continuous hermitian metric $h$ 
of $(E,\nabla)$ is called a tame pure imaginary harmonic bundle,
if  $(E,\nabla,h)_{|U^{\ast}}$ is tame pure imaginary harmonic bundle.

Let $h_i$ $(i=1,2)$ be tame pure imaginary harmonic bundle
of $(E,\nabla)$.
The identity of $E$ induces the flat morphisms
$\Phi:(E,\nabla,h_1)\lrarr (E,\nabla,h_2)$.
We obtain the norms $|\Phi|$ and $|\Phi^{-1}|$
obtained from $h_1$ and $h_2$.

\begin{lem}
The morphism $\Phi$ is bounded.
\end{lem}
\pf
Let $\prolong{(E,h_i)}$ denote the prolongment of $E$
by an increasing order with respect to $h_i$.
Due to Lemma \ref{lem;04.1.2.5},
the morphism $\Phi$ is prolonged to the morphism
$\prolong{(E,h_1)}\lrarr\prolong{(E,h_2)}$,
which preserves the parabolic filtrations.
Since the weight filtrations of
$\prolong{(E,h_i)}_{|P_i}$ are determined by
the residue $\Res_{P_i}(\nabla)$,
the morphism $\Phi$ also preserves the weight filtrations.
Then it follows from the norm estimate of a tame harmonic bundle
on a punctured disc.
(See \cite{s2} or \cite{mochi2}).
\hfill\qed

\begin{lem}\label{lem;04.1.6.4}
We have the following inequalities of the distributions
on $U$:
\begin{equation}\label{eq;a12.28.50}
 \Delta\log|\Phi|^2\leq 0,
\quad\quad
 \Delta\log|\Phi^{-1}|^2\leq 0.
\end{equation}
\end{lem}
\pf
If follows from the Simpson-Weitzenbeck formula
and the boundedness of $\Phi$ and $\Phi^{-1}$
(see Lemma 4.1 and Corollary 4.2 in \cite{s2}).
\hfill\qed

\begin{lem}\mbox{{}}
\begin{itemize}
\item
$|\Phi|$ and $|\Phi^{-1}|$ take the maximum value
at points in $\del X$.
\item
We have the following inequality:
\[
 |\Phi_{|P}|^2+|\Phi^{-1}_{|P}|^2-2r
\leq
 \max\Bigl\{
 |\Phi_{|Q}|^2+|\Phi^{-1}_{|Q}|^2-2r
 \,\Big|\,
 Q\in \del X
 \Bigr\}.
\]
\end{itemize}
\hfill\qed
\end{lem}

\begin{lem}\label{lem;a12.29.5}
Let $R$ be a real number such that
$d_{\PH(r)}(h_{1\,|\,Q},h_{2\,|\,Q})\leq R$
for any point $Q\in\del X$.
Then the following inequality holds for any point $P\in X^{\ast}$:
\[
 d_{\PH(r)}\bigl(
 h_{1\,|\,P},h_{2\,|\,P}
 \bigr)
\leq 
 \left(
 \frac{e^R-e^{-R}}{2R}
 \right)\cdot
 \max\Bigl\{
 d_{\PH(r)}(h_{1\,|\,Q},h_{2\,|\,Q})\,\Big|\,
 Q\in \del X
 \Bigr\}.
\]
\end{lem}
\pf
We always have the following:
\[
 d_{\PH(r)}(h_{1\,|\,P},h_{2\,|\,P})
\leq
 \left(
 \frac{|\Phi_{|P}|^2+|\Phi^{-1}_{|P}|^2-2r}{2}
 \right)^{1/2}.
\]
For any point $Q\in \del X$,
we have the following:
\[
 \left(
 \frac{|\Phi_{|Q}|^2+|\Phi^{-1}_{|Q}|^2-2r}{2}
 \right)^{1/2}
\leq
 \frac{e^R-e^{-R}}{2R}
\cdot d_{\PH(r)}\bigl(h_{1\,|\,Q},h_{2\,|\,Q}\bigr).
\]
Thus we are done.
\hfill\qed

\begin{cor}\label{cor;a12.29.25}
If we have $h_{1\,|\,\del X}=h_{2\,|\,\del X}$,
we have $h_1=h_2$.
\hfill\qed
\end{cor}

\subsection{The uniqueness of tame pure imaginary pluri-harmonic metric}
\label{subsection;04.1.30.210}
%c18.1.tex

\subsubsection{The statement and the reduction to the one dimensional case}

Let $X$ be a smooth projective variety over $\cnum$
and $D$ be a normal crossing divisor of $X$.
Let $(E,\nabla)$ be a flat bundle over $X-D$.
Let $h_1$ and $h_2$ be tame pure imaginary pluri-harmonic metric
of $(E,\nabla)$.
\begin{prop}\label{prop;04.1.3.4}
Assume that $\dim X\geq 1$.
There exists a positive constant $a$
such that $h_0=a\cdot h_1$.
\end{prop}
If $X$ is compact,
the claim is proved by Corlette \cite{corlette}.
We essentially follow his argument.
Since we have to care the infinite energy,
we need some modification of the argument.

We will prove the claim in the case $\dim X=1$ later.
Here we give an argument to reduce the higher dimensional case
to the one dimensional case.
We use an induction on $\dim X$.
We assume the claim holds in the case $\dim X\leq n-1$,
and we show that the claim holds in the case $\dim X=n$.

Let $P$ be any point of $X$.
We take a smooth hypersurface $Y$ of $X$
such that $P\in Y$ and $Y\cap D$ is normal crossing.
Note $\dim (Y)=n-1\geq 1$.
Then there exists a positive constant $a_Y$
such that $h_{0\,|\,Y}=a_Y\cdot h_{1\,|\,Y}$.
In particular, there exists a positive constant
 $a(P)$ such that $h_{0\,|\,P}=a(P)\cdot h_{1\,|\,P}$.

For any point $P$ and $Q$,
we can take a smooth hypersurface $Y$
such that $P,Q\in Y$ and $Y\cap D$ is normal crossing.
Then we have $a(P)=a_Y=a(Q)$,
i.e.,
$h_{0}=h_{1}$.

Thus we have reduced the higher dimensional case
to the one dimensional case.
In the following subsubsections
\ref{subsubsection;04.1.5.15}--\ref{subsubsection;04.1.5.16},
we only consider a simple flat bundle $(E,\nabla)$
on a smooth quasi projective curve
with the smooth projective completion $\overline{C}$.

%c18.tex

\subsubsection{Constantness of the eigenvalues}
\label{subsubsection;04.1.5.15}

The hermitian metric $h_1$ induces the self adjoint morphism $H$
of $E$ with respect to the metric $h_0$.
Let $\alpha_1,\ldots,\alpha_r$ be the eigenvalues of $H$,
and then we have $|\Phi|^2=\sum \alpha_i^2$.

\begin{lem}\label{lem;04.1.3.1}
We have the constantness of $\sum\alpha_i^2$ on $C$.
\end{lem}
\pf
The identity map of $E$ induces the flat morphism
$\Phi:(E,\nabla,h_0)\lrarr(E,\nabla,h_1)$.
Due to Lemma \ref{lem;04.1.6.4},
we obtain the constantness
of $|\Phi|$,
i.e., the constantness of $\sum\alpha_i^2$.
\hfill\qed

\vspace{.1in}

The tame pure imaginary harmonic metrics $h_0$ and $h_1$
induce those on $(\bigwedge^l E,\nabla)$ for any $l$.
By applying Lemma \ref{lem;04.1.3.1},
we obtain the constantness of
any symmetric functions of $\alpha_i^2$ $(i=1,\ldots,r)$,
which implies the constantness of $\alpha_i$.
Thus we obtain the following:
\begin{lem}
There exist the mutually different real numbers $\beta_1,\ldots,\beta_s$
and the decomposition $E=\bigoplus_i E_i$
satisfying the following:
\begin{itemize}
\item 
 $E_i$ are mutually orthogonal with respect to both of
 the metrics $h_0$ and $h_1$.
\item
 On $E_i$, we have $h_1=e^{2\beta_i}\cdot h_0$.
\hfill\qed
\end{itemize}
\end{lem}

We put $\nbigl=\bigoplus e^{-\beta_i}\cdot id_{E_i}$.
Then we have
$h_1(x,y)=h_0\bigl(\nbigl^{-1} x,\nbigl^{-1}y\bigr)$.
We put as follows:
\[
 \nbigl_t:=\bigoplus e^{-t\beta_i}\cdot id_{E_i},
\quad\quad
 h_t(x,y):=h_0\bigl(\nbigl_t^{-1}x,\nbigl_t^{-1}y\bigr).
\]
Let $\theta_t$ denote the $(1,0)$-form for
$(E,\nabla,h_t)$ (see the subsubsection \ref{subsubsection;04.1.5.10}).

\subsubsection{The description by connection forms}

Let $\vece$ be the orthonormal frame of $E$ with respect to
$h_0$ on some open subset of $C$ with a coordinate $z$.
The $(1,0)$-form $A^{1,0}dz$ and
the $(0,1)$-form $A^{0,1}d\bar{z}$ are determined by
$\nabla\vece=\vece\cdot \bigl(A^{1,0}dz+A^{0,1}d\bar{z}\bigr)$.
We have the following:
\[
 d''\vece=\vece\cdot A^{0,1}d\bar{z},
\quad
 d'\vece=\vece\cdot A^{1,0}dz,
\quad
 \delta''\vece=\vece\cdot
 \bigl(-\overline{\lefttop{t}A^{1,0}}\bigr)d\bar{z},
\quad
 \delta'\vece=\vece\cdot
 \bigl(-\overline{\lefttop{t}A^{0,1}} \bigr)dz.
\]
Thus we have the following:
\begin{equation}\label{eq;04.1.5.11}
 \theta_0\vece=\vece\cdot\frac{A^{1,0}+\overline{\lefttop{t}A^{0,1}}}{2}dz,
\quad
 \del\vece=\vece\cdot\frac{A^{1,0}-\overline{\lefttop{t}A^{0,1}}}{2}dz.
\end{equation}
Assume that $\vece$ is compatible with the decomposition
$E=\bigoplus E_i$.
Then we have $\nbigl_t\vece=\vece\cdot L_t$
for the constant diagonal matrices $L_t$.
The frame $\vece\cdot L_t$ is the orthonormal frame
with respect to $h_t$.
Via the frame $\vece\cdot L_t$,
we identify $End(E,E)$ with $M(r)$.
We have the following:
\[
 \nabla\bigl(\vece\cdot L_t\bigr)
=\vece\cdot L_t
\cdot\bigl(
 L_t^{-1}A^{1,0}L_t\cdot dz
+L_t^{-1}A^{0,1}L_t\cdot d\bar{z}
 \bigr).
\]
Hence we have the following:
\[
 \theta_t\bigl(e\cdot L_t\bigr)
=\vece\cdot L_t\cdot
\frac{\bigl(L_t^{-1}A^{1,0}L_t
   +L_t\lefttop{t}\overline{A^{0,1}}L_t^{-1}\bigr)}{2}dz.
\]

The metric $h_t$ induces the hermitian metric 
of $End(E)$.
It induces the skew linear pairing
\[
 \langle\cdot,\cdot\rangle:
 \bigl(
 End(E)\otimes\Omega^{1,0}\bigr)
 \otimes
\overline{
 \bigl(
 End(E)\otimes\Omega^{1,0}\bigr)}\lrarr\Omega^{1,1}.
\]
We have the following formula:
\[
 \langle\theta_t,\theta_t\rangle
=\frac{1}{4}
 \bigl|\bigl|
 L_t^{-1}A^{1,0}L_t+L_t\lefttop{t}\overline{A^{0,1}}L_t^{-1}
 \bigr|\bigr|^2\cdot dz d\bar{z}.
\]
Here $||\cdot||$ denote the norm of matrices.
We have the decompositions:
\[
 A^{1,0}
=\sum A^{1,0}_{i\,j},\quad
 \lefttop{t}\overline{A^{0,1}}
=\sum\bigl( \lefttop{t}\overline{A^{0,1}}\bigr)_{i\,j},\quad\quad
\Bigl(
 A^{1,0}_{i\,j},\,\,
 \bigl(\lefttop{t}\overline{A^{0,1}}\bigr)_{i\,j}
 \in Hom(E_i,E_j)
\Bigr).
\]
We also have the following:
\[
 L_t^{-1}A^{1,0}L_t=\sum_{i\,j}A_{i\,j}^{1,0}\cdot e^{(\beta_i-\beta_j)t}.
\]
Thus we obtain the following:
\[
 \bigl|\bigl|
 L_t^{-1}A^{1,0}L_t
+L_t\lefttop{t}\overline{A^{0,1}}L_t^{-1}
 \bigr|\bigr|^2
=\sum_{i\,j}
 \bigl|\bigl|
 A^{1,0}_{i\,j}\cdot e^{(\beta_i-\beta_j)t}
+\bigl( \lefttop{t}\overline{A^{0,1}}\bigr)_{i\,j}
 \cdot e^{-(\beta_i-\beta_j)t}
 \bigr|\bigr|^2.
\]

\subsubsection{The convexity}

Note the following positivity, or convexity:
\begin{multline}
 \left(\frac{d}{dt}\right)^2
 \bigl|\bigl|
 A^{1,0}_{i\,j}\cdot e^{(\beta_i-\beta_j)t}
+\lefttop{t}\overline{A^{0,1}}_{i\,j}\cdot e^{-(\beta_i-\beta_j)t}
 \bigr|\bigr|^2 \\
=2(\beta_i-\beta_j)^2\cdot
 \left(
 \bigl|\bigl|
 A^{1,0}_{i\,j}\cdot e^{(\beta_i-\beta_j)t}
+\bigl(
 \lefttop{t}\overline{A^{0,1}}\bigr)_{i\,j}
\cdot e^{-(\beta_i-\beta_j)t}
 \bigr|\bigr|^2
+
  \bigl|\bigl|
 A^{1,0}_{i\,j}\cdot e^{(\beta_i-\beta_j)t}
-\bigl(\lefttop{t}\overline{A^{0,1}}_{i\,j}\bigr)
 \cdot e^{-(\beta_i-\beta_j)t}
 \bigr|\bigr|^2
 \right)\geq 0.
\end{multline}
Thus we obtain the following.
\begin{lem}
We have the following:
\[
 \sqrt{-1}\left(\frac{d}{dt}\right)^2
 \langle\theta_t,\theta_t\rangle\geq 0.
\]
The equality holds if and only if
$A_{i\,j}^{1,0}=A_{i\,j}^{0,1}=0$ for any pair $i\neq j$.
\hfill\qed
\end{lem}

\subsubsection{The end of the proof of Proposition \ref{prop;04.1.3.4}}
\label{subsubsection;04.1.5.16}

We have the line bundle $\nbigo(D)$ on $\overline{C}$
with the canonical section $s:\nbigo\lrarr\nbigo(D)$.
Let us take a $C^{\infty}$-hermitian metric $h_D$ of $\nbigo(D)$.
Then we obtain the $C^{\infty}$-function $-\log|s|$ on $C$.
We put $C(R):=\bigl\{P\in C\,\big|\,-\log|s(P)|\leq R\bigr\}$.
Let us consider the following functions:
\[
 F_R(t):=
   \int_{C(R)}\sqrt{-1}\cdot\frac{d}{dt}\langle\theta_t,\theta_t\rangle,
\quad\quad
 F_{R_1,R_2}(t):=
   \int_{C(R_1,R_2)}\sqrt{-1}\cdot\frac{d}{dt}\langle\theta_t,\theta_t\rangle.
\]

\begin{lem}\label{lem;04.1.5.12}
We have $F_{R_1,R_2}(0)\leq F_{R_1,R_2}(1)$ for any $R_1$ and $R_2$.
If $F_{R_1,R_2}(0)=F_{R_1,R_2}(1)$,
then we obtain $(\frac{d}{dt})^2\langle\theta_t,\theta_t\rangle=0$
on $C(R_1,R_2)$.
\end{lem}
\pf
It follows from 
$\sqrt{-1}\cdot\bigl(\frac{d}{dt}\bigr)^2
 \langle\theta_t,\theta_t\rangle\geq 0$.
\hfill\qed

\begin{cor} \label{cor;04.1.5.13}
For any pair $R\geq R'$,
we have $F_R(1)-F_R(0)\geq F_{R'}(1)-F_{R'}(0)\geq 0$.
\hfill\qed
\end{cor}

Let us show that $F_R(0)$ converges to $0$ when $R\to \infty$
by using the assumption 
that $(E,\nabla,h_0)$ is tame pure imaginary harmonic.
We put $\xi:=\frac{d L_t}{dt}$,
which is self adjoint with respect to $h_0$.
We have the following formula:
\[
 \frac{d}{dt}\left(
 L_t^{-1}A^{1,0}L_t+L_t\lefttop{t}\overline{A^{0,1}}L_t^{-1}
 \right)_{|t=0}
=\left[
 \frac{A^{1,0}-\lefttop{t}\overline{A^{0,1}}}{2},\,\xi
 \right]
=\del\xi.
\]
Here we have used the formula (\ref{eq;04.1.5.11}).
Thus we have the following:
\[
 \frac{d}{dt}\langle\theta_t,\theta_t\rangle_{|t=0}
=\langle\del\xi,\theta_0\rangle+\langle\theta_0,\del\xi\rangle
=-\langle\xi,\delbar\theta_0\rangle
+\del\langle\xi,\theta_0\rangle
+\langle\delbar\theta_0,\xi\rangle
-\delbar\langle\theta_0,\xi\rangle
=\del\langle\xi,\theta_0\rangle-\delbar\langle\theta_0,\xi\rangle.
\]
Hence we obtain the following:
\[
 F_R(0)
=\sqrt{-1}\int_{\del C(R)}
\Bigl(
 \langle\xi,\theta_0\rangle
-\langle\theta_0,\xi\rangle
\Bigr).
\]
Let $P$ be a point of $\overline{C}-C$.
We take a coordinate $z=re^{\sqrt{-1}\eta}$ around $P$.
Then we have the description:
\[
 \theta_0=g\cdot\frac{dz}{z}
=g\left(\frac{dr}{r}+\sqrt{-1}d\eta\right).
\]
Here $g$ is an endomorphism of $E$.
We recall that the eigenvalues of $g_{|P}$ is pure imaginary.
Due to Simpson's Main estimate,
we have the decomposition $g=g_0+g_1$
(see the subsection 7.1 in \cite{mochi2}):
\begin{itemize}
\item There is a decomposition
 $E=\bigoplus_{\alpha\in\sqrt{-1}\real} E'_{\alpha}$ which is orthogonal
 with respect to $h_0$.
 The endomorphism $g_0$ is given by
 $\bigoplus \alpha\id_{E'_{\alpha}}$.
\item
 We have the estimate
 $|g_1|\leq C\cdot (-\log r)^{-1}$ for some positive constant $C$.
\end{itemize}
Since $\sqrt{-1}\cdot g_0$ is self adjoint,
we have the cancellation:
$\bigl\langle\xi,\,\sqrt{-1}g_0\cdot d\eta\bigr\rangle
-\bigl\langle\sqrt{-1}g_0\cdot d\eta,\,\xi\bigr\rangle=0$.
We also have the estimate
$\bigl|
 \langle\xi,\sqrt{-1}g_1\rangle\bigr|\leq C\cdot(\log R)^{-1}$
for some positive constant $C$.
Then we obtain that
$\lim_{R\to\infty}F_R(0)=0$.

\begin{lem}\label{lem;04.1.5.14}
We have the convergence
$\lim_{R\to\infty}\bigl(F_R(1)-F_R(0)\bigr)=0$.
\hfill\qed
\end{lem}
\pf
We have only used the property that
$(E,\nabla,h_0)$ is tame pure imaginary harmonic
to show $\lim_{R\to\infty}F_R(0)=0$.
Hence we also obtain 
$\lim_{R\to\infty}F_R(1)=0$.
Then the lemma immediately follows.
\hfill\qed

\vspace{.1in}

We obtain $F_R(1)-F_R(0)=0$ for any $R$,
from Corollary \ref{cor;04.1.5.13} and Lemma \ref{lem;04.1.5.14}.
It implies
$(\frac{d}{dt})^2\langle\theta_t,\theta_t\rangle=0$
for any $t\in[0,1]$,
and thus
$A^{0,1}_{i\,j}=A^{1,0}_{i\,j}=0$
for $i\neq j$,
due to Lemma \ref{lem;04.1.5.12}.
It implies that
the decomposition $E=\bigoplus_i E_i$ is flat
with respect to the connection $\nabla$.

Since $E$ is simple,
we have $E=E_i$ for some $i$,
and thus we obtain
$h_1=e^{2\beta_i}\cdot h_0$.
Therefore the proof of Proposition \ref{prop;04.1.3.4} is accomplished.
\hfill\qed

\section{Depencence on boundary value in the case of a punctured disc}
\label{section;04.1.30.211}

\subsection{The Dirichlet problem and sequence of the boundary values}
\label{subsection;04.1.6.10}
%c13.51.tex

\subsubsection{The Dirichlet problem}

Let $\varphi$ be an element of $GL(r)$,
and $(E,\nabla)$ be a flat bundle on $\Deltabar^{\ast}$
whose monodromy is conjugate with $\varphi$.
Let $h_{\del\Deltabar}$ be a $C^{\infty}$-hermitian metric
on $(E,\nabla)_{|\del\Deltabar}$.
We denote the corresponding $C^{\infty}$-map
$\del\Deltabar\lrarr \PH(r)\big/\langle \varphi\rangle$
by $\psi$.

\begin{prop}\label{prop;a12.29.24}
There exists the tame pure imaginary harmonic metric $h$
on $(E,\nabla)$ such that $h_{|\del\Deltabar}=h_{\del\Deltabar}$.
It is unique up to constant multiplication.
\end{prop}
\pf
The uniqueness immediately follows from Corollary \ref{cor;a12.29.25}.
Hence we have only to show the existence.
The main idea is clearly explained in \cite{lohkamp} and \cite{JZ2}.

We put $T(R_1,R_2)
:=\bigl\{z\in\cnum\,|\,R_1\leq -\log|z|\leq R_2\bigr\}$,
and $T(R):=T(0,R)$.
We use the Poincar\'{e} metric
$g:=|z|^{-2}\cdot\bigl(-\log|z|^2+A\bigr)^{-2}\cdot dz\cdot d\bar{z}$
of $\Deltabar^{\ast}$ for some positive constant $A$.
We also use the real coordinate
$z=\exp\bigl(\sqrt{-1}x-y\bigr)$.

\begin{lem}\label{lem;a12.30.13}
We have a $C^{\infty}$-hermitian metric $h_0$ on $\Deltabar^{\ast}$
satisfying the following:
\[
 h_{0\,|\,\del\Deltabar}=h_{\del\Deltabar},
\quad\quad
 \int_{\Deltabarast}
 \left|
 e(\Psi_{h_0})-\frac{\rho(\varphi)^2}{4\pi^2}\bigl(2y+A\bigr)^2
 \right|\cdot
  \dvol_g
 <\infty.
\]
\end{lem}
\pf
Recall that we have a tame pure imaginary harmonic metric
of $(E,\nabla)$ on $\Deltabarast$,
by using the model bundle for example
(see \cite{mochi2}).
Hence we can take a $C^{\infty}$-hermitian metric $h_0$ of $(E,\nabla)$
satisfying the following:
\begin{itemize}
\item
$h_{0\,|\,\del\Deltabar}=h_{\del\Deltabar}$.
\item
The restriction of $h_0$ to $\Deltabar^{\ast}-T(1)$ is
a tame pure imaginary harmonic metric.
\end{itemize}
Due to Lemma \ref{lem;a12.29.1},
we are done.
\hfill\qed

\begin{rem}
We can avoid to use the model bundle
as in {\rm\cite{JZ2}}.
See also Lemma {\rm \ref{lem;04.1.7.21}},
for example.
\hfill\qed
\end{rem}

The function $F$ is given as follows,
which is integrable with respect to $\dvol_g$:
\begin{equation}\label{eq;04.1.25.15}
 F:=\left|
 e\bigl(\Psi_{h_0}\bigr)
-\frac{\rho(\varphi)^2}{4\pi^2}\cdot(2y+A)^2
 \right|.
\end{equation}

\begin{lem}[Hamilton-Schoen-Corlette \cite{corlette3}]
There exists the twisted harmonic map
$\Phi_n:T(n)\lrarr \PH(r)\big/\langle \varphi\rangle$
satisfying the following:
\[
 \Phi_{n\,|\,\del\Deltabar}=\psi,
\quad\quad
 \Phi_{n\,|\,|z|=e^{-n}}=\Psi_{h_0\,|\,|z|=e^{-n}}.
\]
\[
 \int_{T(n)}e(\Phi_n)\dvol_g
\leq
 \int_{T(n)}e(\Psi_{h_0})\dvol_g.
\]
\end{lem}
\pf
See the proof of Theorem 2.1 of \cite{corlette3}.
\hfill\qed

\begin{lem}\label{lem;a12.29.3}
Let $R$ be a positive number.
For any $n>R$, we have the following inequality:
\[
 \int_{T(R)}e(\Phi_n)\cdot\dvol_g
\leq
 R\cdot\frac{\rho(\varphi)^2}{2\pi}
+\int_{\Deltabarast}F\cdot\dvol_g.
\]
\end{lem}
\pf
We have the following inequalities:
\begin{multline}\label{eq;04.1.25.20}
 \int_{T(n)}e\bigl(\Psi_{h_0}\bigr)\cdot\dvol_g
\geq
 \int_{T(n)}e\bigl(\Phi_{n}\bigr)\cdot\dvol_g
=\int_{T(R)}e\bigl(\Phi_n\bigr)\cdot\dvol_g
+\int_{T(R,n)}e\bigl(\Phi_n\bigr)\cdot\dvol_g \\
\geq
 \int_{T(R)}e\bigl(\Phi_n\bigr)\cdot\dvol_g
+\frac{\rho(\varphi)^2}{2\pi}\cdot(n-R).
\end{multline}
Here we have used Lemma \ref{lem;04.1.25.2}.
On the other hand, we have the following inequality
by our choice of the integrable function $F$ (\ref{eq;04.1.25.15}):
\begin{equation}\label{eq;04.1.25.21}
 \int_{T(n)}e(\Psi_{h_0})\cdot\dvol_g
\leq
 \frac{\rho(\varphi)^2}{2\pi}\cdot n
+\int_{T(n)}F\cdot\dvol_g.
\end{equation}
From (\ref{eq;04.1.25.20}) and (\ref{eq;04.1.25.21}),
we obtain the following:
\[
 \int_{T(R)}e(\Phi_n)\cdot\dvol_g
\leq
 R\cdot\frac{\rho(\varphi)^2}{2\pi}
+\int_{T(n)}F\cdot\dvol_g.
\]
Thus we are done.
\hfill\qed

\vspace{.1in}

\vspace{.1in}
We have the projection $p:\Deltabarast\lrarr \del\Deltabar$
given by $p(e^{\sqrt{-1}x-y})=e^{\sqrt{-1}x}$.
We have the following,
for any point $P=e^{\sqrt{-1}x_0-y_0}\in T(R)$:
\begin{multline}
 d_{\PH(r)}\bigl(\Phi_{n}(P),\psi(p(P))\bigr)
=d_{\PH(r)}\bigl(\Phi_{n}(P),\Phi_{n}(p(P))\bigr)
\leq
 \int_{0}^{y_0}\left|\del_y\Phi_{n}\right|\cdot dy
\leq
 \int_{0}^{R}\left|\del_y \Phi_{n}\right|\cdot dy\\
\leq
 C_R\cdot\left(
 \int_{0}^R\left|\del_y \Phi_{n}\right|^2\cdot dy
 \right)^{1/2}
\leq
 C_R\cdot\left(
 \int_{0}^R e(\Phi_n)\cdot \frac{dy}{(2y+A)^2}
 \right)^{1/2}.
\end{multline}
Hence for any subregion $B\subset T(R)$,
there exists a positive constant $C_{R,B}$,
which is independent of a choice of $n$,
such that the following inequalities hold:
\begin{equation}\label{eq;a12.29.2}
 \int_B d_{\PH(r)}\bigl(\Phi_{n}(P),\psi(p(P))\bigr)^2
  \cdot d\eta \cdot dr
\leq C_{R,B}\cdot
 \int_{T(R)}e(\Phi_n)\cdot\dvol_g.
\end{equation}

\begin{lem}\label{lem;a12.29.30}
Let $\vecn_0$ be any infinite subset of $\nnum$.
Then there exist an infinite subset $\vecn_1\subset\vecn_0$
such that
the sequence
$\bigl\{\Phi_{n}\,\big|\,n\in\vecn_1\bigr\}$ is $C^{\infty}$-convergent
on any compact subset $K\subset\Delta^{\ast}$.
\end{lem}
\pf
From (\ref{eq;a12.29.2}),
we can take a point $P\in\Delta^{\ast}$
and a subsequence $\vecn_2$
such that $\bigl\{\Phi_n(P)\,|\,n\in\vecn_2\bigr\}$ is convergent.
We have the boundedness of 
$e(\Phi_n)$ on any compact subset $K\subset\Delta^{\ast}$,
which can be derived from the boundedness of
the energy $E(\Phi_n)$ (Lemma \ref{lem;a12.29.3})
and Lemma \ref{lem;a12.29.4}.
Then it is standard to show the result by using the boot strapping
argument.
(See \cite{lohkamp} and \cite{schoen-yau}).
\hfill\qed

\vspace{.1in}
Let us fix the subsequence $\vecn_1$ 
for which $\{\Phi_n\,|\,n\in\vecn_1\}$ is convergent
on any compact subset $K\subset\Delta^{\ast}$.
Let $\Phi_{\infty}$ denote the limit,
and let $h$ denote the corresponding harmonic metric
of $(E,\nabla)$.
We also use the notation $\Psi_h$ to denote $\Phi_{\infty}$.

\begin{lem}\label{lem;a12.29.41}
The sequence $\bigl\{\Phi_n\,\big|\,n\in\vecn_1\bigr\}$
is $C^0$-convergent on $T(R)$ for any $R$.
In particular, $\Phi_{\infty\,|\,\del\Delta}=\psi$.
\end{lem}
\pf
Due to our choice of $\vecn_1$,
the sequence $\bigl\{\Phi_{n\,|\,|z|=e^{-k}}\,|\,n\in\vecn_1\bigr\}$
is $C^{\infty}$-convergent.
By our choice of $\Phi_{n}$,
we have $\Phi_{n\,|\,\del\Delta}=\psi$.
Then we obtain the $C^0$-convergence
due to Lemma \ref{lem;a12.29.5}.
\hfill\qed

\vspace{.1in}
Due to Lemma \ref{lem;a12.29.3},
we obtain the following inequality:
\[
  \int_{T(R)} e\bigl(\Psi_h\bigr)\cdot\dvol_g
\leq
 R\cdot \frac{\rho(\varphi)^2}{2\pi}+C.
\] 
Here $C$ denotes a positive constant,
which is independent of $R$.
Then it is easy to obtain an integrable function $\tilde{F}$
on $\Deltabarast$,
satisfying the following inequalities
for any sufficiently large real number $R$:
\[
 \int_{T(R)} e\bigl(\Psi_h\bigr)\cdot\dvol_g
\leq
 \int_{T(R)}\left(
 \frac{\rho(\varphi)^2}{4\pi^2}\cdot\bigl(2y+A\bigr)^2
+\tilde{F}
\right)\cdot\dvol_g.
\]
Therefore the harmonic bundle
$(E,\nabla,h)$ is tame pure imaginary 
due to Proposition \ref{prop;04.1.25.10}.
Thus the proof of Proposition \ref{prop;a12.29.24} is accomplished.
\hfill\qed

\vspace{.1in}
We give a remark on the convergency of the sequence
$\bigl\{\Phi_n\bigr\}$.
\begin{lem}\label{lem;a12.29.40}
The sequence $\bigl\{\Phi_n\bigr\}$ is $C^{\infty}$-convergent
on any compact subset $K\subset\Delta^{\ast}$,
and it is $C^0$-convergent on $T(R)$ for any $R$.
\end{lem}
\pf
Let $\vecn_0$ and $\vecn_1$ be as in Lemma \ref{lem;a12.29.30}.
We have the limit $\Phi_{\infty}$.
It is tame pure imaginary harmonic.
It also satisfies $\Phi_{\infty\,|\,\del\Delta}=\psi$.
Then it follows that $\Phi_{\infty}$ is independent
of a choice of $\vecn_1$, due to Corollary \ref{cor;a12.29.25}.
It implies the desired convergence properties
of the sequence $\bigl\{\Phi_n\,\big|\,n\in\nnum\bigr\}$.
\hfill\qed

%c13.61.tex

\subsubsection{Dependence of a convergent sequence of boundary values}

Let $\{\varphi_i\}$ be a sequence in $GL(r)$
converging to $A$.
Let $\Phi_i:\Deltabar^{\ast}\lrarr
            \PH(r)\big/\big\langle \varphi_i\big\rangle$
and $\Phi:\Deltabar^{\ast}\lrarr\PH(r)\big/\big\langle \varphi\bigr\rangle$
be maps,
which are corresponding to tame pure imaginary harmonic bundles
on $\Deltabar^{\ast}$.
\begin{prop}\label{prop;a12.29.50}
Assume the following:
\begin{itemize}
\item
 The sequence of the boundary values
 $\bigl\{\Phi_{i\,|\,\del\Deltabar}\bigr\}$
 converges to $\Phi_{|\del\Deltabar}$ in $C^{\infty}$-sense.
\item
 There exists a positive constant $C$, which is independent of $i$ and $k$,
 satisfying the following:
\[
 \int_{T(k)}e(\Phi_{i\,|\,T_k})\cdot\dvol_g
\leq
 k\cdot \frac{\rho(\varphi_i)^2}{2\pi}+C.
\]
\end{itemize}
Then the sequence $\{\Phi_i\}$ is $C^{\infty}$-convergent
on any compact subset $K\subset \Delta^{\ast}$,
and $C^0$-convergent on any $T_k$.
\end{prop}
\pf
The argument is essentially same as the proof of
Proposition \ref{prop;a12.29.24}
and Lemma \ref{lem;a12.29.40}.
Hence we only indicate an outline.

Let $\vecn_0$ be an infinite subset of $\nnum$.
Then we can take an infinite subset $\vecn_1\subset\vecn_0$
such that
the sequence $\{\Phi_i\,|\,n\in\vecn_1\}$ is 
$C^{\infty}$-convergent on any compact subset $K\subset\Delta^{\ast}$,
by an argument similar to the proof of Lemma \ref{lem;a12.29.30}.
We denote the limit by $\Phi_{\infty}$.
We can show that the sequence
$\bigl\{\Phi_i\,\big|\,n\in\vecn_1\bigr\}$ is
$C^0$-convergent on any $T_k$
by an argument similar to the proof of Lemma \ref{lem;a12.29.41}.
In particular, we obtain
$\Phi_{\infty\,|\,\del\Deltabar}=\Phi_{|\del\Deltabar}$.
By using the estimate of the energy and 
Proposition \ref{prop;04.1.25.10},
we can show that $\Phi_{\infty}$ corresponds
to the tame pure imaginary harmonic bundle.
Thus we obtain $\Phi_{\infty}=\Phi$
from the maximum principle (Lemma \ref{lem;a12.29.5}).
Since $\Phi_{\infty}$ is independent of a choice of $\vecn_1$,
we obtain the desired convergency of $\bigl\{\Phi_i\bigr\}$.
\hfill\qed

\subsection{Family version}
\label{subsection;04.1.30.212}

%c13.7.tex

\subsubsection{Estimate for a flat family}

\label{subsubsection;04.1.7.6}

Let $X$ be a $C^{\infty}$-manifold.
Let $(E,\nabla)$ be a flat bundle
on $\Deltabar^{\ast}\times X$.
Let $h_{1}$ be a $C^{\infty}$ hermitian metric
of $(E,\nabla)_{|\del\Deltabar\times X}$.
We have the hermitian metric $h$ on $(E,\nabla)$
such that
$(E,\nabla,h)_{|\Deltabarast\times P}$ are
tame pure imaginary harmonic bundles on $\Deltabar^{\ast}$
for any $P\in X$.

Let us pick any point $P_0$ of $X$,
and take an appropriate neighbourhood $U$ of $P_0$.
For any point $P\in U$,
we may assume to have the natural identification
$(E,\nabla)_{|\Deltabar^{\ast}\times P}\simeq
 (E,\nabla)_{|\Deltabar^{\ast}\times P_0}$.
Hence we obtain the family of the maps
$\phi_P:\Deltabar^{\ast}\times\{P_0\}\lrarr
 \PH(r)\big/\langle\varphi\rangle$.
Here $\varphi$ denotes the monodromy.

\begin{lem} \label{lem;a12.29.51}
We regard $h$ as the family of tame pure imaginary harmonic metrics
$\bigl\{h_{\Deltabar^{\ast}\times P}\,\big|\,P\in X\bigr\}$.
Then the assumption in Proposition {\rm\ref{prop;a12.29.50}}
is satisfied.
\end{lem}
\pf
It is easy to see that
we can take a $C^{\infty}$-hermitian metric $h_2$
of $(E,\nabla)_{|\Deltabarast\times U}$
such that the following holds:
\begin{itemize}
\item
$h_{2\,|\,\del\Deltabar\times U}=h_{1\,|\,\del\Deltabar\times U}$,
\item
The restriction of $h_2$ to $\Deltabarast-S_1$
is the fixed tame pure imaginary harmonic metric,
i.e.,
$h_{2\,|\,(\Deltabarast-S_1)\times\{P\}}
=h_{2\,|\,(\Deltabarast-S_1)\times\{P_0\}}$
under the isomorphism
$(E,\nabla)_{|\Deltabarast\times\{P\}}\simeq
 (E,\nabla)_{|\Deltabarast\times\{P_0\}}$.
\end{itemize}

Let $\Psi_{h_2\,P}:
 \Deltabar^{\ast}\times \{P_0\}\lrarr \PH(r)\big/\langle\varphi\rangle$
denote the map corresponding to
$h_{2\,|\,\Deltabarast\times P}$.
There exists a positive constant $C$,
which is independent of a choice of $P$,
such that the following holds:
\[
 \int_{\Deltabar^{\ast}}
 \Big|
 e(\Psi_{h_2\,P})-\frac{\rho(\varphi)^2}{4\pi^2}
 \cdot\bigl(-\log |z|^2+A\bigr)^2
 \Big|\cdot
 \dvol<C.
\]
Then the claim of Lemma \ref{lem;a12.29.51}
is clear from the proof of
Proposition \ref{prop;a12.29.24}.
\hfill\qed

\begin{lem}\label{lem;a12.29.52}
The corresponding map $\Psi_h$ is continuous.
We also have the continuity of
$\del_z^l\delbar_z^m\Psi_h$ for any $l$ and $m$.
\end{lem}
\pf
It follows from Proposition \ref{prop;a12.29.50}
and Lemma \ref{lem;a12.29.51}.
\hfill\qed

\vspace{.1in}
For any point $P\in X$,
we put as follows:
\[
 (E^P,\nabla^P,h^P):=
 (E,\nabla,h)_{|\Deltabar^{\ast}\times P}.
\]
We obtain $\delbar_{E^P}$ and $\theta^P$.
\begin{lem}
$\bigl\{\delbar_{E^P}\,\big|\,P\in X\bigr\}$
and
$\bigl\{\theta^P\,\big|\,P\in X\bigr\}$ are continuous
with respect to $P$.
\end{lem}
\pf
It follows from Lemma \ref{lem;a12.29.52}.
\hfill\qed

\vspace{.1in}
We have the description $\theta^P=f^P_0\cdot dz/z$.
We obtain the family of the characteristic polynomial:
\[
 \det(t-f_0^P)=\sum t^i\cdot A_i(z,P).
\]
Here $A_i(z,P)$ is continuous on $\Deltabar^{\ast}\times X$.
The restrictions of $A_i$
to $\Deltabar^{\ast}\times P$ are holomorphic,
and it is naturally extended to the holomorphic functions
on $\Deltabar\times P$.
Then it is easy to see that $A_i$ is naturally extended
to the continuous function on $\Deltabar\times X$,
by using the Cauchy's integral formula.

For any point $P\in X$,
we can take positive constants $R(P),C(P),\epsilon(P)$
as in Lemma \ref{lem;a12.28.40}.

\begin{lem}\label{lem;a12.30.30}
We can take the constants $R(P),C(P),\epsilon(P)$
independently of $P$.
\end{lem}
\pf
It follows from the continuity of $A_i$
on $\Deltabar\times X$.
\hfill\qed

\begin{lem}\label{lem;a12.30.31}
Positive constants
$R',C',\epsilon'$
in Lemma {\rm\ref{lem;a12.28.42}},
a positive constant $C''$ in Lemma {\rm\ref{lem;a12.29.55}},
a positive constant $C_1$ in Lemma {\rm\ref{lem;a12.29.56}}
and a positive constant $C_0$ in Corollary {\rm\ref{cor;a12.29.57}}
for $(E^P,\delbar_{E^P},h^P,\theta^P)$
can be taken independently of $P$.
\hfill\qed
\end{lem}

Let $(x_1,\ldots,x_l)$ be a coordinate of $X$.
Recall that $\del\Psi_h/\del x_i$ are defined almost everywhere.

\begin{lem}\label{lem;04.1.7.7}
For any point $T\in X$,
we have the following:
\[
 \left|
 \frac{\del \Psi_h}{\del x_i}
 \right|(P,T)
\leq
 \max\left\{
 \left|\frac{\del \Psi_h}{\del x_i}\right|(Q,T)
 \,\Big|\,Q\in\del\Deltabar
 \right\}
\]
\end{lem}
\pf
It follows from Lemma \ref{lem;a12.29.5}.
\hfill\qed
%c13.8.tex

\subsubsection{Another family}
\label{subsubsection;04.1.7.20}

Recall we have the continuous family of $C^{\infty}$-twisted maps
$\bigl\{\tilde{\Phi}_{(\xi_1,\kappa)}\,
 \big|\,(\xi_1,\kappa)\in \real\times[0,1]\bigr\}$
(see the subsubsection \ref{subsubsection;a12.30.11}).

\begin{lem}\label{lem;04.1.7.21}
For any $(\xi_1,\kappa)\in\real\times[0,1]$,
we have twisted harmonic maps
$\Psi_{h\,(\xi_1,\kappa)}:
\Deltabar^{\ast}\lrarr \PH(r)
   \big/\big\langle \tilde{f}_2(\kappa)\big\rangle$
which satisfy the following:
\begin{itemize}
\item
$\Psi_{h\,(\xi_1,\kappa)\,|\,\del\Deltabar}
=\tilde{\Phi}_{(\xi_1,\kappa)\,|\,\del\Deltabar}$.
\item
There exists a positive constant $C$,
which is independent of a choice of $(\xi_1,\kappa)$,
such that the following holds:
\begin{equation}\label{eq;a12.30.14}
 \int_{\Deltabarast}
 \left|
 e\bigl(\tilde{\Psi}_{h\,(\xi_1,\kappa)}\bigr)
 -\frac{\rho(f_2)^2}{4\pi^2}
 \cdot(2\eta_2+b)^2
 \right|
 \dvol<C.
\end{equation}
\end{itemize}
\end{lem}
\pf
We have the finiteness as in Lemma \ref{lem;a12.30.12}.
We can take $\Phi_{(\xi_1,\kappa)}$ as the hermitian metric
in Lemma \ref{lem;a12.30.13},
and we construct the twisted harmonic maps
$\Psi_{h\,(\xi_1,\kappa)}$ as in
the subsection \ref{subsection;04.1.6.10}.
Then we obtain the estimate of the integral (\ref{eq;a12.30.14}),
independently of $(\xi_1,\kappa)$.
\hfill\qed

\begin{cor}
For any non-negative integers $l$ and $m$,
the family
$\bigl\{\del^m_{\xi_2}\del^l_{\eta_2}\Psi_{h,(\xi_1,\kappa)}
 \,\big|\,(\xi_1,\kappa)\in\real\times[0,1]\bigr\}$
is continuous with respect to $(\xi_1,\kappa)$.
\hfill\qed
\end{cor}

From the family
$\bigl\{\Psi_{h\,(\xi_1,\kappa)}\,
 \big|\,(\xi_1,\kappa)\in \real\times[0,1]\bigr\}$,
we obtain the continuous family of the flat bundles
$(E,\nabla)_{(\xi_1,\kappa)}$
with the tame pure imaginary harmonic metrics $h_{(\xi_1,\kappa)}$
on $\Deltabar^{\ast}$.
We obtain the Higgs fields $\theta_{(\xi_1,\kappa)}$.
Similarly to Lemma \ref{lem;a12.30.30}
and Lemma \ref{lem;a12.30.31},
we obtain the following.

\begin{lem}
We can take the constants
$R(\xi_1,\kappa),C(\xi_1,\kappa),\epsilon(\xi_1,\kappa)$
in Lemma {\rm\ref{lem;a12.28.40}}
independently of $(\xi_1,\kappa)$.
\hfill\qed
\end{lem}

\begin{lem}\label{lem;04.1.7.22}
Positive constants
$R',C',\epsilon'$
in Lemma {\rm\ref{lem;a12.28.42}},
a positive constant $C''$ in Lemma {\rm\ref{lem;a12.29.55}},
a positive constant $C_1$ in Lemma {\rm\ref{lem;a12.29.56}}
and a positive constant $C_0$ in Corollary {\rm\ref{cor;a12.29.57}}
for $(E,\delbar_E,h,\theta)_{(\xi_1,\kappa)}$
can be taken independently of $(\xi_1,\kappa)$.
\hfill\qed
\end{lem}

\section{Control of the energy of twisted map on a Kahler surface}
\label{section;04.1.30.410}

\subsection{Around smooth points of divisors}
\label{subsection;04.1.8.20}
%c14.tex

\subsubsection{Metric}
\label{subsubsection;04.1.28.1}

Let $X$ be a closed region of $\cnum^l$
whose boundary is continuous.
We consider the set $X\times\Delta(R)\subset\cnum^l\times\cnum$.
We use the real coordinate $(x_1,\ldots,x_{2l})$ of $\cnum^l$.
We denote $\del/\del x_i$ also by $\del_i$.

We use the real coordinate $z=x+\sqrt{-1}y$ of $\cnum$.
We also use the polar coordinate $z=r\cdot e^{\sqrt{-1}\eta}$ of
$\cnum$.
We denote $\del/\del r$ and $\del/\del \eta$
also by $\del_r$ and $\del_{\eta}$ respectively.
We denote the origin of $\cnum$ by $O$.

We have the natural complex structure $J_X$ of $X$,
which is induced by the complex structure of $\cnum^l$.
For any point $P\in X$,
we denote $\{P\}\times\Delta(R)$ by $P\times\Delta(R)$
for simplicity.
We have the natural complex structure $J_{P\times\Delta(R)}$
on $P\times\Delta(R)$,
which is induced by the complex structure of $\cnum$.

Let $J$ be a complex structure of $X\times\Delta(R)$
satisfying the following conditions:
\begin{condition}\label{condition;04.1.30.1}\mbox{{}}
\begin{itemize}
\item
 The natural inclusion $X\times\{O\}\lrarr X\times\Delta(R)$
 is a holomorphic embedding.
\item
 Let $P$ be any point of $X$.
 We have the natural inclusion
 $T_{(P,O)}(P\times\Delta(R))\lrarr T_{(P,O)}(X\times\Delta(R))$
 of the tangent spaces.
 Then $T_{(P,O)}(P\times\Delta(R))$ is a subspace 
 of $\cnum$-vector space $T_{(P,O)}(X\times\Delta(R))$,
 and the restriction of $J$ to $T_{(P,O)}(P\times\Delta(R))$
 is same as $J_{P\times\Delta(R)}$.
\end{itemize}
\end{condition}

\begin{rem}
The complex structure $J$ is not necessarily
same as the natural complex structure of $X\times\Delta(R)$.
\hfill\qed
\end{rem}

Let $P$ be any point of $X$.
Let $\zeta$ be a holomorphic function of $(X\times\Delta(R),J)$
defined on a neighbourhood $U$ of $(P,O)$
satisfying the following:
\[
 d\zeta_{|(P,O)}\neq 0,
\quad
 \zeta^{-1}(0)=U\cap (X\times \{O\}).
\]
Let $z$ be the holomorphic coordinate of $\cnum$.
It induces a $C^{\infty}$-function on $X\times\Delta(R)$,
which is not necessarily holomorphic 
with respect to the complex structure $J$.
\begin{lem}
There exists the complex number $a(P)$
such that $dz_{|(P,O)}=a(P)\cdot d\zeta_{|(P,O)}$.
\end{lem}
\pf
Both of $d\zeta_{|(P,O)}$ and $dz_{|(P,O)}$ vanish on $T(X\times \{O\})$.
Hence we have only to compare them on 
the one dimensional subspace $T_{(P,O)}(P\times\Delta(R))$.
Since both of $d\zeta_{|(P,O)}$ and $dz_{|(P,O)}$ are
$\cnum$-linear, we obtain the complex number $a(P)$
as claimed.
\hfill\qed

\vspace{.1in}
For simplicity,
we consider the case $\zeta$ is defined over $X\times\Deltabar$.
We can take a $C^{\infty}$-function
$a:X\times\Deltabar\lrarr \cnum^{\ast}$
such that $(a\cdot d\zeta)_{|X\times\{O\}}=dz_{|X\times\{O\}}$.
We have the following estimates:
\begin{equation}\label{eq;04.1.7.2}
  |a\zeta|^2=|z|^2+O(|z|^3),
\quad\quad
 \left|
 d\zeta\left(\frac{\del}{\del x_i}\right)
 \right|
=O\bigl(|z|\bigr).
\end{equation}

Let $g_1$ be a Kahler metric of $(X\times\Deltabar,J)$,
and $\omega_1$ be the associated Kahler form.
We have the following estimates
with respect to the metric $g_1$:
\begin{equation}\label{eq;04.1.7.1}
 \bigl|
 dz-d(a\zeta)
 \bigr|_{g_1}=O(|z|),
\quad
 \frac{dz}{z}-\frac{d\zeta}{\zeta}
=O(1).
\end{equation}

Let $b$ be any $C^{\infty}$-function on
$X\times\Deltabar$
such that $-\log|\zeta|^2+b>0$ on $X\times\Deltabarast$
and $\bigl(b+\log|a|^2\bigr)_{|X\times\{O\}}>0$ on $X\times\{O\}$.
We put $A:=\bigl(b+\log|a|^2\bigr)_{X\times\{O\}}$.
We have the estimate:
\begin{equation}\label{eq;04.1.14.1}
 \Bigl|
 (-\log|z|^2+A)
-(-\log|\zeta|^2+b)
 \Bigr|=O(|z|).
\end{equation}
We have the following formula:
\begin{multline}\label{eq;04.1.7.3}
 \del\delbar \log\Bigl(
 -\log |\zeta|^2+b
 \Bigr)
=\del\left(
 \frac{\delbar\bigl( -\log|\zeta|^2+b\bigr)}
 {-\log|\zeta|^2+b}
 \right)\\
=
-\frac{1}{\bigl(-\log|\zeta|^2+b\bigr)^2}
 \left(\frac{d\zeta}{\zeta}-\del b\right)
\wedge
 \left(\frac{d\bar{\zeta}}{\bar{\zeta}}-\delbar b\right)
+\frac{-\del\delbar b}
 {(-\log|\zeta|^2+b)}.
\end{multline}

Assume that the following 
gives the Kahler form:
\[
 \omega=\omega_1
-\sqrt{-1}\del\delbar
 \log\Bigl(-\log|\zeta|^2+b\Bigr).
\]
Let $g$ denote the associated Kahler metric.
Let $g_2$ (resp. $g_3$) denote the hermitian form
corresponding to the first (resp. second) term 
in the right hand side of (\ref{eq;04.1.7.3}).
Thus we have the decomposition $g=g_1+g_2+g_3$.

Recall we use the polar coordinate $z=r\cdot e^{\sqrt{-1}\eta}$.
\begin{lem}\label{lem;a12.31.3}
We have the following estimates:
\[
 \bigl|
 \del_r
 \bigr|^2_g=(1+F_1)\cdot\frac{1}{r^2\cdot(-\log r^2+A)^2},
\quad\quad
 |F_1|=O\bigl(r\cdot(-\log r^2+A)\bigr).
\]
\[
 \bigl|
 \del_{\eta}
 \bigr|^2_g=(1+F_2)\cdot\frac{1}{(-\log r^2+A)^2},
\quad\quad
 |F_2|=O\bigl(r\cdot(-\log r^2+A)\bigr).
\]
\end{lem}
\pf
It follows from (\ref{eq;04.1.7.1}), (\ref{eq;04.1.14.1})
and (\ref{eq;04.1.7.3}).
Note $|\del_{\eta}|_{g_1}=O(r)$.
\hfill\qed

\begin{lem}\label{lem;04.1.7.5}
We have the following estimate:
\[
 \bigl|
 g(\del_r,\del_{\eta})
 \bigr|=O\bigl(
 (-\log r^2+A)^{-2}
 \bigr),
\quad
 \left|
 \frac{g(\del_r,\del_{\eta})}{|\del_{\eta}|_g^2}
 \right|
=O(1).
\]
\end{lem}
\pf
We have the following estimates: 
\[
\bigl|g_1\bigl(\del_r,\del_{\eta}\bigr) \bigr|=O(r),
\quad\quad
 \left|
g_3\bigl(\del_r,\del_{\eta}\bigr)
 \right|
=O\bigl(
 r\cdot(-\log r^2+A)^{-1}
 \bigr).
\]
Since
$\del_r$ and $\del_{\eta}$ are orthogonal
with respect to the hermitian form $dz\cdot d\bar{z}$,
we obtain the following estimate
from (\ref{eq;04.1.7.1}).
\[
  \left|
 g_2\bigl( \del_r,\del_{\eta}\bigr)
 \right|
=O\Bigl(
 \bigl(
 -\log r^2+A
 \bigr)^{-2}
 \Bigr).
\]
Then we obtain the result
from (\ref{eq;04.1.7.3}).
\hfill\qed

\vspace{.1in}

We put as follows:
\[
 e_1:=\del_{\eta},
\quad\quad
 e_2:=\del_r
  -\frac{g(\del_r,\del_{\eta})}{|\del_{\eta}|^2_g}\cdot\del_{\eta}.
\]
\begin{lem}\label{lem;a12.31.4}
We have the following estimate:
\[
 |e_2|_g^2=
 (1+F_3)\cdot\frac{1}{r^2(-\log r^2+A)^2},
\quad\quad
|F_3|=O\bigl(r\cdot(-\log r^2+A)\bigr).
\]
\end{lem}
\pf
We have 
$|e_2|_g^2
=|\del_r|^2_g-|g(\del_r,\del_{\eta})|^2\cdot|\del_{\eta}|_g^{-2}$.
Then it is easy to check the claim.
\hfill\qed

\vspace{.1in}

We use the notation $\del_i$ to denote $\del/\del x_i$.
\begin{lem}
We have the following estimate:
\[
 g\bigl(\del_i,e_1\bigr)
=O\bigl(
 (-\log r^2+A)^{-2}
 \bigr),
\quad\quad
 g\bigl(\del_i,e_2\bigr)
=O\bigl(r^{-1}\cdot (-\log r^2+A)^{-2}\bigr).
\]
\end{lem}
\pf
We have estimates
$ (g_1+g_3)\bigl(\del_i,e_1\bigr)
=O(r)$.
We also have
$ g_2\bigl(\del_i,e_1\bigr)
=O\bigl((-\log r^2+A)^{-2}\bigr)$,
which follows from (\ref{eq;04.1.7.1}).
Thus we obtain the first estimate.
Similarly, we can obtain the second estimate.
\hfill\qed

\vspace{.1in}
We put as follows:
\[
 f_i:=
 \del_i-\frac{g(\del_i,e_1)}{|e_1|_g^2}\cdot e_1
-\frac{g(\del_i,e_2)}{|e_2|^2_g}\cdot e_2.
\]
\begin{lem}\label{lem;a12.31.2}
We have the following estimate:
\[
 \bigl|
 g\bigl(f_i,f_j\bigr)-g_1(\del_i,\del_j)
 \bigr|
=O\bigl((-\log r^2+A)^{-1}\bigr).
\]
\end{lem}
\pf
We have the following:
\[
 g(f_i,f_j)
=g(\del_i,\del_j)
-\sum_{a=1,2}\frac{g(\del_i,e_a)\cdot g(\del_j,e_a)}{|e_a|_g^2}.
\]
The second term is $O\bigl((-\log r^2+A)^{-2}\bigr)$,
which follows from  the previous lemmas.
It is easy to obtain the following estimate,
from (\ref{eq;04.1.7.3}):
\[
 g(\del_i,\del_j)=g_1(\del_i,\del_j)+O\bigl((-\log r^2+A)^{-1}\bigr).
\]
Thus we are done.
\hfill\qed

%c14.1.tex

\subsubsection{The volume form}

The restriction of $g_{1}$ to $X\times\{0\}$
induces the volume form of $X$,
which we denote by $\dvol_X$.

\begin{lem}\label{lem;04.1.7.16}
Let $\pi:X\times\Deltabar\lrarr X$ be the projection.
We have the following estimate:
\begin{equation}\label{eq;a12.31.1}
 \dvol_g=(1+F_4)\cdot
 \pi^{\ast}\dvol_X\cdot\frac{dr\cdot d\eta}{r\cdot(-\log r^2+A)^2},
\quad\quad
 |F_4|=O\bigl((-\log r^2+A)^{-1}\bigr).
\end{equation}
\end{lem}
\pf
We have the naturally defined subbundle
$T\Deltabar\subset T(X\times\Deltabar)$.
Let $H$ be the orthogonal complement of $T\Deltabar$
in $T(X\times\Deltabar)$.
The bundle $H$ is generated by $f_1,\ldots,f_{2l}$.

The metric $g$ induces the volume form of
$T\Deltabar$, which is as follows, due to Lemma \ref{lem;a12.31.3}
and Lemma \ref{lem;a12.31.4}:
\[
 (1+G_1)\cdot \frac{dr\cdot d\eta}{r\cdot(-\log r^2+A)^{2}},
\quad\quad
 |G_1|= O\bigl(r\cdot(-\log r^2+A)\bigr).
\]
On the other hand, we have the volume form of $H$ induced by $g$,
which is as follows, due to Lemma \ref{lem;a12.31.2}:
\[
 (1+G_2)\cdot\pi^{\ast}\dvol_X,
\quad\quad
 |G_2|=O\bigl((-\log r^2+A)^{-1}\bigr).
\]
Then (\ref{eq;a12.31.1}) follows.
\hfill\qed

\subsubsection{An inequality for energy}
\label{subsubsection;04.1.30.10}

Let $Y$ be a Riemannian manifold
on which $\pi_1(X\times\Deltabarast)$ acts.
Let $F:X\times\Deltabarast\lrarr Y\big/\pi_1(X\times\Deltabarast)$
be a twisted map.
Let $e_g(F)$ denote the energy function of $F$
with respect to the metric $g$ of $X$.

\begin{lem}\label{lem;04.1.7.8}
There exists a positive number $C$
such that the following holds for any twisted map $F$:
\begin{equation}\label{eq;a12.31.10}
\left|
 e_g(F)-(-\log r^2+A)^2\cdot \left|\frac{\del F}{\del \eta}\right|^2
\right|
\leq
 C\cdot\left(
 \sum_{i=1}^l\left|\frac{\del F}{\del x_i}\right|^2
+r^2\cdot(-\log r^2+A)^2\cdot
 \left|\frac{\del F}{\del r}\right|^2
+\left|\frac{\del F}{\del \eta}\right|^2
 \right).
\end{equation}
\end{lem}
\pf
We have the following estimate from 
Lemma \ref{lem;a12.31.3},
which is independent of $F$:
\[
|e_1|_g^{-2}\cdot
 \bigl|dF(e_1)\bigr|^2
=\bigl((-\log r^2+A)^2+G_1\bigr)
\cdot
 \left|
 \frac{\del F}{\del\eta}
 \right|^2,
\quad\quad
 |G_1|=O\bigl(r\cdot (-\log r^2+A)\bigr).
\]
From Lemma \ref{lem;04.1.7.5},
there exists a positive constant,
which is independent of $F$, such that the following holds:
\[
 \bigl|dF(e_2)\bigr|^2
\leq
 C\cdot
 \bigl(\bigl|
 \del_r F\bigr|^2
+\bigl|
 \del_{\eta}F\bigr|^2
\bigr).
\]
Thus we have the following estimate
from Lemma \ref{lem;a12.31.4},
which is independent of $F$:
\[
 |e_2|_g^{-2}\cdot
 \bigl|dF(e_2) \bigr|^2
=O\Bigl(
 r^2\cdot(-\log r^2+A)^2
\cdot \bigl(
 |\del_r F|^2+|\del_{\eta}F|^2\bigr)
 \Bigr).
\]
We take the orthogonal frame 
$\tilde{f}_1,\ldots,\tilde{f}_{2l}$ of the orthogonal complement 
of $T\Deltabar$,
by applying the orthogonalization of Schmidt to $f_1,\ldots,f_{2l}$.
Then it is easy to obtain the following estimate,
which is independent of $F$:
\[
 |dF(\tilde{f}_i)|^2
=O\Bigl( 
 \sum_i |\del_iF|^2
+r^2\cdot|\del_r F|^2+|\del_{\eta}F|^2
 \Bigr).
\]
Then (\ref{eq;a12.31.10}) immediately follows.
\hfill\qed

%c14.2.tex

\subsubsection{An estimate for some twisted map}

Let $(E,\nabla)$ be a flat connection on $X\times\Deltabarast$.
Let $h'$ be a $C^{\infty}$-hermitian metric of 
$E_{|X\times\del\Delta}$.
\begin{lem}
We have the unique continuous hermitian metric $h$ of $E$
satisfying the following:
\begin{itemize}
\item
 The restriction $h_{|P\times\Deltabarast}$ is
 tame pure imaginary harmonic metric of
 $(E,\nabla)_{|P\times\Deltabarast}$.
 Here the metric of $P\times\Deltabarast$ is
 given by $|z|^{-2}\cdot(-\log|z|^2+A)^{-2}\cdot dz\cdot d\bar{z}$.
\item
 $h_{|X\times\del\Deltabar}=h'$.
\end{itemize}
\end{lem}
\pf
See the subsubsection \ref{subsubsection;04.1.7.6}.
\hfill\qed

\vspace{.1in}

For any point $P$,
let $\varphi^P$ denote the monodromy
of the harmonic bundle $(E,h,\nabla)_{|P\times\Deltabarast}$.
The conjugacy classes are independent of a choice of $P$.
Thus we denote $\varphi^P$ simply by $\varphi$.
Let $\Psi_h$ denote the twisted map corresponding to $h$.

\begin{lem}\label{lem;04.1.7.9}
There exists a positive constant $C_0$,
which is independent of $P$,
such that the following holds:
\[
 \left|
 \frac{\del\Psi_h}{\del r}
 \right|^2\cdot r^2\cdot(-\log r^2+A)^2
\leq C_0,
\quad\quad
 \left|
 \left|\frac{\del \Psi_h}{\del\eta}\right|^2
 \cdot (-\log r^2+A)^2
-\frac{\rho(\varphi)^2}{4\pi^2}
\cdot (-\log r^2+A)^2
 \right|
\leq C_0.
\]
\end{lem}
\pf
It follows from Lemma \ref{lem;a12.30.31}.
\hfill\qed

\begin{lem}\label{lem;04.1.7.10}
$\del_i\Psi_h$ is defined almost everywhere,
and we have the following,
for any point $P\in X$:
\[
 \bigl|\del_i\Psi_h\bigr|^2(T,P)
\leq
 \max\Bigl\{
 |\del_i\Psi_h|^2(T,Q)\,\Big|\,
 Q\in \del\Deltabar
 \Bigr\}.
\]
\end{lem}
\pf
It follows from Lemma \ref{lem;04.1.7.7}.
\hfill\qed

\begin{prop}
The following function is bounded on $X\times\Deltabarast$:
\[
 \left|
 e(\Psi_h)-\frac{\rho(\varphi)^2}{4\pi^2}\cdot(-\log|z|^2+A)^2
 \right|.
\]
\end{prop}
\pf
We obtain the boundedness of the following,
due to Lemma \ref{lem;04.1.7.8}, Lemma \ref{lem;04.1.7.9}
and Lemma \ref{lem;04.1.7.10}:
\[
 \left|
 e(\Psi_h)-\left|\frac{\del\Psi_h}{\del \eta}\right|^2\cdot(-\log|z|^2+A)^2
 \right|.
\]
Then the desired boundedness immediately follows
from Lemma \ref{lem;04.1.7.9}.
\hfill\qed

\begin{cor}\label{cor;04.1.8.25}
We have the boundedness of the following
function on $X\times\Deltabarast$:
\begin{equation}\label{eq;04.1.7.11}
 \left|
 e(\Psi_h)
-\frac{\rho(\varphi)^2}{4\pi^2}\bigl(-\log |\zeta|^2+b\bigr)^2
 \right|.
\end{equation}
In particular, we obtain the following finiteness:
\begin{equation}\label{eq;04.1.7.12}
 \int_{X\times\Deltabarast}
 \left|
 e(\Psi_h)
-\frac{\rho(\varphi)^2}{4\pi^2}
 \bigl(-\log |\zeta|^2+b\bigr)^2
 \right|\cdot\dvol_g<\infty.
\end{equation}
\end{cor}
\pf
The function
$\bigl| \bigl(-\log|\zeta|^2+b\bigr)^2-\bigl(-\log|z|^2+A\bigr)^2\bigr|$
is bounded on $X\times\Deltabar$,
due to (\ref{eq;04.1.14.1}).
Thus we obtain the boundedness (\ref{eq;04.1.7.11}).
The finiteness (\ref{eq;04.1.7.12}) immediately follows from
(\ref{eq;04.1.7.11}).
\hfill\qed

%c14.3.tex

\subsubsection{A priori lower estimate of the energy for arbitrary map}
\label{subsubsection;04.1.15.1}

Let $b_1$ be a $C^{\infty}$-function on $X\times \Deltabar$.
The functions $w$, $s$ and $\phi$ are given by
the relations $e^{-(b+b_1)/2}\cdot \zeta=w=s\cdot e^{\sqrt{-1}\phi}$.
We have $dw/w=d\zeta/\zeta-d(b+b_1)/2$.
We use the coordinate change on $X\times\Deltabarast$
given as follows:
\[
 (x_1,\ldots,x_{2l},r,\eta)\longleftrightarrow
 (x_1,\ldots,x_{2l},s,\phi).
\]
We have the equality
$-\log|\zeta|^2+b=-\log s^2-b_1>0$.

A $1$-form $\sum_i A_i\cdot dx_i+B\cdot ds+C\cdot d\phi$
on $X\times\Deltabarast$ is called independent of $\phi$,
if $A_i$, $B$ and $C$ are independent of $\phi$.
It is similarly defined that an $i$-form is
independent of $\phi$.

\begin{lem}\label{lem;04.1.7.15}
Let $\omega$ be a $C^{\infty}$-form on $X\times\Deltabar$.
Then we have a decomposition
$\omega=\omega_1+\omega_2$, where $\omega_1$ is
independent of $\phi$ and $|\omega_2|_{g_1}=O(s)$.
\end{lem}
\pf We can check it elementarily by using Taylor development.
\hfill\qed

\begin{lem}\label{lem;04.1.14.2}
We have the following estimate with respect to
the Kahler form $g_1$ on $X\times\Deltabar$:
\[
 d\eta-d\phi=O(1),\quad\quad
 \frac{dr}{r}-\frac{ds}{s}=O(1).
\]
\end{lem}
\pf
The first claim follows from the following:
\[
  2\sqrt{-1}d\eta=\frac{dz}{z}-\frac{d\bar{z}}{\bar{z}}
=\frac{dw}{w}-\frac{d\bar{w}}{\bar{w}}+O(1)
=2\sqrt{-1}d\phi+O(1).
\]
The second claim can be shown similarly.
\hfill\qed

\begin{lem}\label{lem;04.1.14.3}
We have the following estimate of the volume form:
\[
 \dvol_g=(H_1+H_2)\cdot dx_1\cdots dx_{2l}\cdot
 \frac{ds\cdot d\phi}{s\cdot(-\log s^2-b_1)^2}.
\]
Here $H_1$ is independent of $\phi$
and $0<C_1<H_1<C_2$ for some positive constants $C_i$.
We also have the estimate $|H_2|=O\bigl((-\log s^2-b_1)^{-2}\bigr)$.
\end{lem}
\pf
We apply Lemma \ref{lem;04.1.7.15} to the $C^{\infty}$-form
$\del\delbar b$,
and then we have
$g(\del_i,\del_j)=
 g_1(\del_i,\del_j)+K_1\cdot\bigl(-\log s^2-b_1\bigr)^{-1}+K_2$,
where $K_1$ is independent of $\phi$
and $K_2=O\bigl((-\log s^2-b_1)^{-2}\bigr)$.
Then the claim can be shown
by an argument similar to the proof of Lemma \ref{lem;04.1.7.16}
using Lemma \ref{lem;04.1.14.2}.
\hfill\qed

\begin{lem}\label{lem;a12.31.20}
We have the following estimate:
\[
 \frac{\del \eta}{\del \phi}=1+O(s),
\quad\quad
 \frac{\del r}{\del \phi}=O(s^2).
\]
\end{lem}
\pf
It follows from Lemma \ref{lem;04.1.14.2}.
\hfill\qed

\begin{lem}
We have the following estimate:
\[
 \left|
 \del_{\phi}
 \right|^2_g=(1+I_1)\cdot (-\log s^2-b_1)^{-2},
\quad
 |I_1|=O\bigl(s\cdot(-\log s^2-b_1)\bigr).
\]
\end{lem}
\pf
We obtain the following
from Lemma \ref{lem;a12.31.20}
and Lemma \ref{lem;a12.31.3}:
\[
 \left|\del_{\phi}\right|^2_g
=\bigl(1+O(s)\bigr)\cdot\left|\del_{\eta}\right|_g^2
+O(s^4)\cdot\bigl|\del_r\bigr|_g^2
+O(s^2)\cdot g\bigl(\del_{\eta},\del_r\bigr)
=(1+I_1)\cdot(-\log s^2-b_1)^2.
\]
Thus we are done.
\hfill\qed

\vspace{.1in}
Let $F:X\times\Deltabarast\lrarr \PH(r)\big/\langle\varphi\rangle$ 
be any twisted map.
Then we obtain the following:
\begin{multline}
 e_g(F)\cdot \dvol_g
\geq
 \bigl|\del_{\phi}\bigr|_g^{-2}
\cdot
 \bigl|\del_{\phi}F
 \bigr|^2\cdot\dvol_g\\
=(1+I_1)\cdot(-\log s^2-b_1)^2\cdot
 \bigl|\del_{\phi}F\bigr|^2
\times (H_1+H_2)\cdot dx_1\cdots dx_{2l}
 \cdot \frac{ds\cdot d\phi}{s\cdot (-\log s^2-b_1)^2} \\
=:(H_1+H_3)\cdot
 \bigl|\del_{\phi}F\bigr|^2\cdot dx_1\cdots dx_{2l}
 \cdot \frac{ds\cdot d\phi}{s}.
\end{multline}
Here $H_1$ and $H_2$ are as in Lemma \ref{lem;04.1.14.3},
and we put $H_3:=(1+I_1)\cdot (H_1+H_2)-H_1$.
We remark $|H_3|=O\bigl((-\log s^2-b_1)^{-2}\bigr)$.

\begin{lem}\label{lem;04.1.8.10}
There is a function $H_4$, which is independent of $F$,
satisfying  $|H_4|=O\bigl((-\log s^2-b_1)^{-2}\bigr)$
and the following inequality:
\[
 \int_{0}^{2\pi}
 \bigl|
 \del_{\phi}F
 \bigr|^2\cdot (H_1+H_3)\cdot s^{-1}\cdot d\phi
\geq
 \int_{0}^{2\pi}\Bigl(\frac{\rho(\varphi)^2}{4\pi^2}-H_4\Bigr)
 \cdot (H_1+H_2)\cdot s^{-1}\cdot d\phi.
\]
\end{lem}
\pf
Due to Lemma \ref{lem;a12.29.10}, we have the following:
\[
\int_{0}^{2\pi}
 \bigl|
 \del_{\phi}F
 \bigr|^2\cdot \bigl(1+H_1^{-1}\cdot H_3\bigr)\cdot d\phi
\geq
 \rho(\varphi)^2\cdot
 \left(
\int_0^{2\pi}(1+H_1^{-1}\cdot H_3)^{-1}\cdot
 d\phi
 \right)^{-1}.
\]
Since we have $|H_1^{-1}\cdot H_3|=O\bigl((-\log s^2-b_1)^{-2}\bigr)$,
we obtain the following,
for some $H_4'$ such that $|H_4'|=O\bigl((-\log s^2-b_1)^{-2}\bigr)$:
\[
\int_{0}^{2\pi}
 \bigl|
 \del_{\phi}F
 \bigr|^2\cdot \bigl(1+H_1^{-1}\cdot H_3\bigr)\cdot d\phi
\geq
 \frac{\rho(\varphi)^2}{2\pi}-H_4'.
\]
Hence we obtain the following:
\begin{multline}
 \int_0^{2\pi}
 \bigl|\del_{\phi}F\bigr|^2\cdot (H_1+H_3)\cdot s^{-1}\cdot d\phi
\geq
 \int_0^{2\pi}
 \left(\frac{\rho(\varphi)^2}{4\pi^2}-\frac{H_4'}{2\pi}\right)
 \cdot H_1\cdot s^{-1}\cdot d\phi \\
\geq
 \int_0^{2\pi}
 \left(
 \frac{\rho(\varphi)^2}{4\pi^2}-H_4
 \right)
\cdot (H_1+H_2)\cdot s^{-1}\cdot d\phi.
\end{multline}
Here we take $H_4$ appropriately by using
$|H_2|=O\bigl((-\log s^2-b_1)^{-2}\bigr)$.
Thus we are done.
\hfill\qed

\vspace{.1in}
Let us consider the region $\nbigs$ of the following form:
\[
 \nbigs:=\bigl\{
 (P,Q)\in X\times\Deltabarast\,\big|\,
 R_1\leq -\log s(Q)\leq R_2
 \bigr\}.
\]
\begin{cor}\label{cor;04.1.12.1}
There is a function $H_4$,
which is independent of $F$,
satisfying $|H_4|=O\bigl((-\log|\zeta|^2+b)^{-2}\bigr)$
and the following inequality:
\[
 \int_{\nbigs}e(F)\cdot \dvol_g
\geq
 \int_{\nbigs}
 \bigl|\del_{\phi}F\bigr|^2\cdot
 \bigl|\del_{\phi}\bigr|_g^{-2}\cdot \dvol_g
\geq
 \int_{\nbigs}
\left(
 \frac{\rho(\varphi)^2}{4\pi^2}
 -H_4
\right)\cdot(-\log |\zeta|^2+b)^2\cdot \dvol_g.
\]
Note the finiteness
$\int_{X\times\Deltabarast} |H_4|\cdot(-\log |\zeta|^2+b)^2\dvol_g<\infty$.
\hfill\qed
\end{cor}

\subsection{Around the intersection}
\label{subsection;04.1.1.2}
%c15.2.tex

\subsubsection{Preliminary}
\label{subsubsection;04.1.7.40}

We use the real coordinate $(x_1,y_1,x_2,y_2)$
given by $z_i=\exp(\sqrt{-1}x_i-y_i)$ on $(\Deltabarast)^2$.
For the moment, we use the Poincar\'{e} metric $g_0$ on $(\Deltabarast)^2$:
\[
 g_0:=\sum_{a=1,2}
 \frac{dx_a\cdot dx_a+dy_a\cdot dy_a}
 {\bigl(2y_a+A\bigr)^2}.
\]
Let $(E,\nabla)$ be a flat bundle on $(\Deltabarast)^2$.
Let $f_i$ $(i=1,2)$ be monodromies
around $z_i=0$.
We put as follows:
\[
 (\Deltabarast)^2
\supset
 Y:=\bigl\{
 (z_1,z_2)\in(\Deltabarast)^2\,\big|\,
 |z_1|=|z_2| \bigr\}
\simeq \openclosed{0}{1}\times T^2.
\]
We have the twisted map $\psi':Y\lrarr \PH(r)\big/\langle f_1,f_2\rangle$
given as follows:
\[
 \psi'(x_1,x_2,y_1,y_2)=
F\bigl(B\cdot\log(y_1+y_2+A),x_1,x_2\bigr)
=F\bigl(B\cdot \log(2y_1+A),x_1,x_2\bigr).
\]
(See the subsubsections \ref{subsubsection;04.1.27.100}
and \ref{subsubsection;a12.30.10}.)
Here $B$ denotes sufficiently large positive constant.
We remark that the morphism $\Phi$ in the subsubsection
\ref{subsubsection;a12.30.10} is just twisted by the isomorphisms.

Let $\pi:(\Deltabarast_w)^2\lrarr(\Deltabarast_z)^2$
be the map given by
$\pi(w_1,w_2)=(w_1,|w_1|\cdot w_2)$.
Then we obtain the isomorphism
$\Deltabarast_{w_1}\times \del\Delta_{w_2}\simeq Y$.

We use the real coordinate $(\xi_1,\xi_2,\eta_1,\eta_2)$
of $(\Deltabarast_w)^2$,
given by $w_i=\exp\bigl(\sqrt{-1}\xi_i-\eta_i\bigr)$.
Then we have the relation
$x_i=\xi_i$ $(i=1,2)$, $y_1=\eta_1$ and $y_2=\eta_1+\eta_2$.

We have the $C^{\infty}$-metric $h'$
of $\pi^{-1}(E)_{|\Deltabarast\times\del\Deltabar}$,
which corresponds to a twisted map $\psi'$.
We take the unique continuous metric $h$
of $\pi^{-1}(E)$ on $(\Deltabarast)^2$
satisfying the following conditions:
\begin{itemize}
\item
 $h_{|\Deltabarast\times \del\Deltabar}=h'$.
\item
 The restrictions
 $(\pi^{-1}(E),\pi^{-1}\nabla,h)_{|P\times\Deltabarast_{w_2}}$
 is tame pure imaginary harmonic bundle,
 for any point $P\in\Deltabarast$.
\end{itemize}
Let $\Psi_h:(\Deltabarast)^2\lrarr\PH(r)\big/\langle f_1,f_2\rangle$
denote the corresponding twisted map.

We have the following estimate, due to Lemma \ref{lem;a12.30.5}:
\[
 \left|
 \left|\frac{\del \psi'}{\del\xi_1}\right|^2
-\frac{\rho(\varphi)^2}{4\pi^2}
 \right|
\leq C_1\cdot e^{-C_2\cdot B\log(2\eta_1+A)}
\leq \frac{C_3}{(2\eta_1+A)^{C_2B}}.
\]
Hence we obtain the following, due to Lemma \ref{lem;04.1.7.7}:
\begin{equation}\label{eq;04.1.7.35}
 \left|
 \frac{\del \Psi_h}{\del\xi_1}
 \right|^2
\leq \frac{\rho(\varphi)^2}{4\pi^2}+\frac{C_3}{(2\eta_1+A)^{C_2B}}.
\end{equation}

We have the following estimate, due to Lemma \ref{lem;a12.30.6}:
\[
 \left|\frac{\del\psi'}{\del\eta_1}
 \right|
=\frac{B\cdot \sqrt{\sum \alpha_i^2}}{2\eta_1+A}.
\]
Therefore we obtain the following, due to Lemma \ref{lem;04.1.7.7}:
\begin{equation}\label{eq;04.1.7.30}
 \left|
 \frac{\del \Psi_h}{\del \eta_1}
 \right|
\leq
 \frac{B\cdot \sqrt{\sum \alpha_i^2}}{2\eta_1+A}.
\end{equation}

\begin{lem}\label{lem;04.1.7.25}
There exists a positive constant $C$ such that the following holds:
\[
 \left|
 \frac{\del \Psi_h}{\del \eta_2}
 \right|^2\leq \frac{C}{(2\eta_2+A)^2},
\quad\quad
 \left|
 \left|\frac{\del\Psi_h}{\del\xi_2}\right|^2
-\frac{\rho(f_2)^2}{4\pi^2}
 \right|
\leq \frac{C}{(2\eta_2+A)^2}.
\]
\end{lem}
\pf
It follows from Lemma \ref{lem;04.1.7.22}.
We remark that
the boundary value $\tilde{\Phi}$ in the subsubsection
\ref{subsubsection;04.1.7.20} is obtained from $\Phi$
(see the subsubsection \ref{subsubsection;a12.30.10}),
and $\Phi$ is just twisted by isomorphisms from $F$,
which is the boundary value we consider here.
Hence we can use Lemma \ref{lem;04.1.7.22} here.
\hfill\qed

\subsubsection{The twisted maps and the estimates of their energy
 on $Z_1$ and $Z_2$}
\label{subsubsection;04.1.29.5}

We reformulate the result in the subsubsection
\ref{subsubsection;04.1.7.40}.
We put
$Z_1:=\bigl\{(z_1,z_2)\,\big|\,|z_1|\geq |z_2|\bigr\}$.
We have the naturally defined projection
$\pi_1:Z_1\lrarr\Deltabarast_{z_1}$.
For any point $P\in\Deltabarast_{z_1}$,
we have
$\pi_1^{-1}(P)\simeq
  \bigl\{z_2\,\big|\,0<|z_2|\leq |z_1(P)|\bigr\}$.
We take the continuous hermitian metric $h_{Z_1}$
of $E_{|Z_1}$ satisfying the following:
\begin{itemize}
\item
 $h_{Z_1\,|\,Y}=h'$.
\item
 For any point $P\in\Deltabarast_{z_1}$,
 the restriction $h_{Z_1\,|\,\pi_1^{-1}(P)}$ is
 a tame pure imaginary harmonic metric.
\end{itemize}

\begin{lem}\label{lem;04.1.8.1}
We have the following estimate:
\[
 \left|
 \frac{\del \Psi_{h_{Z_1}}}{\del x_1}
 \right|^2
\leq\frac{\rho(f_1)^2}{4\pi^2}
+O\bigl((2y_1+A)^{-C_2B}\bigr),
\]
\[
 \left|
 \frac{\del \Psi_{h_{Z_1}}}{\del x_2}
 \right|^2
\leq \frac{\rho(f_2)^2}{4\pi^2}
+O\bigl(\bigl(2(y_2-y_1)+A\bigr)^{-2}\bigr),
\]
\[
 \left|
 \frac{\del\Psi_{h_{Z_1}}}{\del y_1}
 \right|^2
=O\Bigl(
 (2y_1+A)^{-2}+\bigl(2(y_1-y_2)+A\bigr)^{-2}
 \Bigr),
\]
\[
 \left|
 \frac{\del \Psi_{h_{Z_1}}}{\del y_2}
 \right|^2
=O\bigl(\bigl(2(y_2-y_1)+A\bigr)^{-2}\bigr).
\]
\end{lem}
\pf
It follows from (\ref{eq;04.1.7.35}),
(\ref{eq;04.1.7.30}),
Lemma \ref{lem;04.1.7.25}
and the relations
$x_i=\xi_i$ $(i=1,2)$, $y_1=\eta_1$ and $y_2=\eta_1+\eta_2$.
\hfill\qed

\vspace{.1in}

\begin{lem}\label{lem;04.1.8.2}
We have the following estimate on $Z_1$:
\[
 \left|
 \del_{z_1}
 \right|^{-2}_{g_0}\cdot
 \left|
 \del_{z_1}\Psi_{h_{Z_1}}
 \right|^2
\leq\frac{(2y_1+A)^2}{2}\cdot
 \left(
 \frac{\rho(f_1)^2}{4\pi^2}
+O\Bigl(
 (2y_1+A)^{-C_2B}
+(2y_1+A)^{-2}
+\bigl(2(y_2-y_1)+A\bigr)^{-2}
 \Bigr)
 \right).
\]
\[
  \left|
 \del_{z_2}
 \right|^{-2}_{g_0}\cdot
 \left|
 \del_{z_2}\Psi_{h_{Z_1}}
 \right|^2
\leq\frac{(2y_2+A)^{2}}{2}
 \left(
  \frac{\rho(f_2)^2}{4\pi^2}
+O\Bigl(
 \bigl(2(y_2-y_1)+A\bigr)^{-2}
 \Bigr)
 \right).
\]
\end{lem}
\pf
It follows from Lemma \ref{lem;04.1.8.1}.
\hfill\qed

\begin{cor}\label{cor;04.1.8.3}
We have the following estimate of the energy
$e_{g_0}\bigl(\Psi_{h_{Z_1}}\bigr)$ with respect to the metric $g_0$
on $Z_1$:
\[
  e_{g_0}\bigl(\Psi_{h_{Z_1}}\bigr)
\leq\sum_{i=1,2}\frac{\rho(f_i)^2}{4\pi^2}(2y_i+A)^2
+O\left(
 1
+\frac{(2y_1+A)^2}{\bigl(2(y_2-y_1)+A\bigr)^2}
+\frac{(2y_2+A)^2}{\bigl(2(y_2-y_1)+A\bigr)^2}
 \right).
\]
The last term in the right hand side is 
integrable on $Z_1$ with respect to the measure $\dvol_{g_0}$
induced by the metric $g_0$.
\hfill\qed
\end{cor}

Similarly we put
$Z_2:=\bigl\{(z_1,z_2)\,\big|\,|z_1|\geq |z_2|\bigr\}$.
Let $\pi_2:Z_2\lrarr \Delta_{z_2}^{\ast}$ be the naturally defined projection.
We have the continuous hermitian metric $h_{Z_2}$
of $E_{|Z_2}$ satisfying the following:
\begin{itemize}
\item
 $h_{Z_2\,|\,Y}=h'$.
\item
 For any point $P\in\Deltabarast_{z_2}$,
 the restriction $h_{Z_2\,|\,\pi_2^{-1}(P)}$ is
 tame pure imaginary harmonic metric.
\end{itemize}

We obtain the following lemma.
\begin{lem}\label{lem;04.1.8.4}
Lemma {\rm\ref{lem;04.1.8.1}}, 
Lemma {\rm\ref{lem;04.1.8.2}}
and Corollary {\rm\ref{cor;04.1.8.3}}
for $h_{Z_2}$ hold, when $(x_1,y_1)$ and $(x_2,y_2)$ are exchanged
and $Z_1$ is replaced with $Z_2$.
\hfill\qed
\end{lem}

\vspace{.1in}

Since $h_{Z_1}$ and $h_{Z_2}$ coincides on $Y$,
they give the continuous and locally $L_1^2$
hermitian metric $h$ on $(\Deltabar^{\ast})^2$.
\begin{lem}\label{lem;04.1.7.60}
There exists an integrable function $J_i$ $(i=1,2)$
with respect to the measure $\dvol_{g_0}$,
such that the following holds for $h$:
\[
 \left|
 \del_{z_1}
 \right|^{-2}_{g_0}\cdot
 \left|
 \del_{z_1}\Psi_{h}
 \right|^2
=\frac{(2y_1+A)^2}{2}\cdot
 \frac{\rho(f_1)^2}{4\pi^2}
+J_1
\]
\[
  \left|
 \del_{z_2}
 \right|^{-2}_{g_0}\cdot
 \left|
 \del_{z_2}\Psi_{h}
 \right|^2
=\frac{(2y_2+A)^{2}}{2}
  \frac{\rho(f_2)^2}{4\pi^2}
+J_2
\]
\[
 e_{g_0}\bigl(\Psi_{h}\bigr)
\leq\sum_{i=1,2}\frac{\rho(f_i)^2}{4\pi^2}(-2y_i+A)^2
+2(J_1+J_2).
\]
\end{lem}
\pf
It follows from Corollary \ref{cor;04.1.8.3}
and Lemma \ref{lem;04.1.8.4}.
\hfill\qed
%c15.3.tex

\subsubsection{Perturbation of the metric and the estimate of
 the energy function}
\label{subsubsection;04.1.8.8}

Let $g_1$ be a $C^{\infty}$-Kahler metric of $\Deltabar^2$.
Let us consider the Kahler metric $g=g_0+g_1$.
We put as follows:
\[
 e_1:=\del_{z_1},
\quad
 e_2:=\del_{z_2}
-\frac{g(\del_{z_2},\del_{z_1})}{|\del_{z_1}|_g^2}\cdot\del_{z_1}.
\]
We have the following estimate:
\begin{equation}\label{eq;04.1.29.2}
 \bigl|e_2\bigr|_g^2
=\bigl|\del_{z_2}\bigr|^2_g\cdot
\left(
1-\frac{\bigl|g(\del_{z_1},\del_{z_2})\bigr|^2}{\bigl|\del_{z_1}\bigr|^2}
\right)
=\bigl|\del_{z_2}\bigr|^2_g\cdot
\left(
1+O\Bigl(
 |z_1|^{2}\cdot |z_2|^{2}\cdot
 \bigl(-\log|z_1|^2+A\bigr)^2\cdot
 \bigl(-\log|z_2|^2+A\bigr)^2
 \Bigr)
\right).
\end{equation}

\begin{lem}\label{lem;04.1.29.3}
Let $\Phi:(\Deltabarast)^2\lrarr \PH(r)\big/\langle f_1,f_2\rangle$
be any twisted map.
The energy function $e_g(\Phi)$ can be described as follows:
\[
 e_g\bigl(\Phi\bigr)
=2\sum_{i=1,2}\bigl(1+G_i\bigr)\cdot
 \bigl|\del_{z_1}\bigr|_g^{-2}\cdot \bigl|\del_{z_i}\Phi\bigr|^2.
\]
We have the estimate of $G_i$:
\[
 |G_i|=
 O\Bigl(
 |z_1|\cdot |z_2|\cdot
 \bigl(-\log|z_1|^2+A\bigr)\cdot
 \bigl(-\log|z_2|^2+A\bigr)
 \Bigr).
\]
The estimate is independent of a choice of $\Phi$.
\end{lem}
\pf
We have the following equality:
\[
 \bigl|d\Phi(e_2)\bigr|^2
=\bigl|\del_{z_2}\Phi\bigr|^2
-2\Re\left(
 \frac{g(\del_{z_1},\del_{z_2})}{\bigl|\del_{z_1}\bigr|^2_g}
  \bigl(\del_{z_2}\Phi,\del_{z_1}\Phi\bigr)
 \right)
+
 \frac{\bigl|g(\del_{z_1},\del_{z_2})\bigr|^2}
 {\bigl|\del_{z_1}\bigr|_g^4}
\cdot
 \bigl|\del_{z_1}\Phi\bigr|^2.
\]
Hence we have the following inequality:
\[
 \Bigl|
 \bigl|d\Phi(e_2)\bigr|^2
-\bigl|\del_{z_2}\Phi\bigr|^2
 \Bigr|
\leq
 2\cdot
 \frac{\bigl|g(\del_{z_1},\del_{z_1})\bigr|}
 {\bigl|\del_{z_1}\bigr|_g^2}
\cdot \bigl|\del_{z_1}\Phi\bigr|\cdot\bigl|\del_{z_2}\Phi\bigr|
+
 \frac{\bigl| g(\del_{z_1},\del_{z_2})\bigr|^2}{\bigl|\del_{z_1}\bigr|_g^4}
 \cdot
\bigl|\del_{z_1}\Phi\bigr|^2.
\]
Hence we obtain the following inequality:
\begin{multline}\label{eq;04.1.29.1}
 \Bigl|
 |e_2|_g^{-2}\cdot
 \bigl|d\Phi(e_2)\bigr|^2
-\bigl|\del_{z_2}\bigr|_g^{-2}\cdot
 \bigl|\del_{z_2}\Phi\bigr|^2
 \Bigr|
\leq
 |e_2|_g^{-2}\cdot
 \Bigl| \bigl|d\Phi(e_2)\bigr|^2-\bigl|\del_{z_2}\Phi\bigr|^2
 \Bigr|
+\bigl|
 |e_2|_g^{-2}-\bigl|\del_{z_2}\bigr|^{-2}_g
 \bigr|
\cdot \bigl|\del_{z_2}\Phi\bigr|^2 \\
\leq
 C_1\frac{\bigl|g(\del_{z_1},\del_{z_2})\bigr|}
{\bigl|\del_{z_1}\bigr|^2_g\cdot \bigl|\del_{z_2}\bigr|_g^{2}}
\cdot
 \bigl|\del_{z_2}\Phi\bigr|\cdot \bigl|\del_{z_1}\Phi\bigr|
+
C_1\frac{\bigl|g\bigl(\del_{z_1},\del_{z_2}\bigr)\bigr|^2}
  {\bigl|\del_{z_1}\bigr|_g^2\cdot \bigl|\del_{z_2}\bigr|_g^2}
 \cdot \bigl|\del_{z_1}\bigr|^{-2}_g\cdot
 \bigl|\del_{z_1}\Phi\bigr|^2
+\Bigl|
 |e_2|^{-2}_g-\bigl|\del_{z_2}\bigr|_g^{-2}
 \Bigr|
\cdot
 \bigl|\del_{z_2}\Phi\bigr|^2.
\end{multline}
Here $C_1$ is a positive constant,
which depends only on $g_1$.
The first term in the right hand side
of (\ref{eq;04.1.29.1}) is dominated as follows:
\[
 C_1\cdot \frac{\bigl|g(\del_{z_1},\del_{z_2})\bigr|}
 {\bigl|\del_{z_1}\bigr|_g\cdot \bigl|\del_{z_2}\bigr|_g}
\cdot
 \Bigl(
 \bigl|\del_{z_1}\bigr|_g^{-1}\cdot \bigl|\del_{z_1}\Phi\bigr|
 \Bigr)
\cdot
 \Bigl(
 \bigl|\del_{z_2}\bigr|_g^{-1}\cdot \bigl|\del_{z_2}\Phi\bigr|
 \Bigr)
\leq
 C_1\cdot
 \frac{\bigl|g(\del_{z_1},\del_{z_2})\bigr|}
 {\bigl|\del_{z_1}\bigr|_g\cdot \bigl|\del_{z_2}\bigr|_g}
\cdot
 \Bigl(
 \bigl|\del_{z_1}\bigr|^{-2}_g
\cdot \bigl|\Phi\bigr|^2
+\bigl|\del_{z_2}\bigr|^{-2}_g
\cdot \bigl|\Phi\bigr|^2
 \Bigr).
\]
We can control the third term in the right hand side 
of (\ref{eq;04.1.29.1})
by using (\ref{eq;04.1.29.2}).
Then it is easy to derive the claim of the lemma,
by using the formula
$e_g\bigl(\Phi\bigr)=2\sum_{i=1,2} |e_i|_g^{-2}\cdot \bigl|d\Phi(e_i)\bigr|^2$.
\hfill\qed

\begin{lem}\label{lem;04.1.7.61}
Let $\Psi_h$ be the twisted map given in 
the subsubsection {\rm\ref{subsubsection;04.1.29.5}}.
We have the following estimate:
\begin{equation}\label{eq;04.1.8.6}
 \bigl|
 e_g(\Psi_h)-e_{g_0}(\Psi_h)
 \bigr|
=
 O\left(
 \sum_{i=1,2}
 \bigl(-\log|z_i|^2+A\bigr)\cdot|z_i|
 \cdot\bigl|\del_{z_i}\bigr|_{g_0}^{-2}
 \cdot\bigl|\del_{z_i}\Psi_h\bigr|^2
\right).
\end{equation}
The right hand side of {\rm(\ref{eq;04.1.8.6})}
is integrable with respect to
the measure $\dvol_{g_0}$.
\end{lem}
\pf
The estimate (\ref{eq;04.1.8.6}) follows from 
Lemma \ref{lem;04.1.29.3} and
the equality
$\bigl|\del_{z_i}\bigr|_g^2=
 \bigl|\del_{z_i}\bigr|^2_{g_0}
+\bigl|\del_{z_i}\bigr|^2_{g_1}$.
The integrability of the right hand side of (\ref{eq;04.1.8.6})
follows from Lemma \ref{lem;04.1.7.60}.
\hfill\qed

\subsubsection{The volume form and the estimate of the energy}

The volume form for the metric $g$ is given as follows:
\[
 \dvol_{g}=
 \left|
 |\del_{z_1}|_{g}^2\cdot|\del_{z_2}|_g^2
-g(\del_{z_1},\del_{z_2})
 \right|\cdot dz_1d\bar{z}_1dz_2d\bar{z}_2.
\]
We have
$|\del_{z_i}|_g^2=|z_i|^{-2}(-\log|z_i|^2+A)^2+g_1(\del_{z_i},\del_{z_i})$,
and $g(\del_{z_1},\del_{z_2})=g_1(\del_{z_1},\del_{z_2})$.
Here $g_1(\del_{z_i},\del_{z_j})$ is $C^{\infty}$ on $\Deltabar^2$.

\begin{lem}\label{lem;a12.31.30}
$\dvol_g$ is of the following form:
\[
 \dvol_g=\Bigl(
 (1+F_1)\cdot(1+F_2)+F_3
 \Bigr)\cdot\dvol_{g_0}.
\]
Here $F_i$ $(i=1,2)$ are of the form
 $\tilde{F}_i\cdot |z_i|^2\cdot(-\log|z_i|^2+A)^2$
for $C^{\infty}$-functions $\tilde{F}_i$ on $\Deltabar^2$,
and $F_3$ is of the form
$\tilde{F}_3\cdot \prod_{i=1,2} |z_i|^2\cdot(-\log|z_i|^2+A)^2$
for $C^{\infty}$-function $\tilde{F}_3$.
\hfill\qed
\end{lem}

\begin{lem}
There exists an integrable function $J_3$ with respect to $\dvol_g$
such that the following holds:
\begin{equation}\label{eq;04.1.8.7}
 e_g\bigl(\Psi_h\bigr)
\leq \sum_{i=1,2}\frac{\rho(f_i)^2}{4\pi^2}\bigl(-\log|z_i|^2+A\bigr)^2
+J_3.
\end{equation}
\end{lem}
\pf
From Lemma \ref{lem;04.1.7.60} and Lemma \ref{lem;04.1.7.61},
there exists an integrable function $J_3$ with respect to $\dvol_{g_0}$
such that the estimate (\ref{eq;04.1.8.7}) holds.
The integrability of $J_3$ with respect to $\dvol_{g}$
follows from Lemma \ref{lem;a12.31.30}.
\hfill\qed

\begin{cor}\label{cor;04.1.8.26}
There exists an integrable function $J_3$ with respect to $\dvol_g$,
such that the following holds for any compact region $K$
of $(\Deltabarast)^2$:
\[
 \int_{K}e_g(\Psi_h)\dvol_g
\leq
 \int_K
\left(
 \frac{\rho(f_1)^2}{4\pi^2}(-\log |z_1|^2+A)^2
+\frac{\rho(f_2)^2}{4\pi^2}(-\log|z_2|^2+A)^2
+J_3
 \right)\cdot\dvol_g.
\]
\hfill\qed
\end{cor}

%c15.4.tex

\subsubsection{A priori lower estimate for arbitrary map}
\label{subsubsection;04.1.22.22}

Let $\Phi:(\Deltabarast)^2\lrarr\PH(r)\big/\langle f_1,f_2\rangle$
be any twisted map.
We would like to obtain the lower bound of the energy
with respect to the metric $g$ given in
the subsubsection \ref{subsubsection;04.1.8.8}.

\begin{lem}
We have the following:
\[
 \dvol_g=(1+F_3')(1+F_1)(1+F_2)\dvol_{g_0}.
\]
Here $F_1$ and $F_2$ are as in Lemma {\rm\ref{lem;a12.31.30}},
and we have $|F_3'|=O\bigl(|z_1|\cdot |z_2|\bigr)$.
\end{lem}
\pf
It immediately follows from Lemma \ref{lem;a12.31.30}.
\hfill\qed

\vspace{.1in}

We use the real coordinate $z_i=\exp\bigl(\sqrt{-1}x_i-y_i\bigr)$
as before.
Recall Lemma \ref{lem;04.1.29.3}.
We have the following:
\begin{multline}
 2\bigl(1+G_i\bigr)\cdot
 \bigl|\del_{z_i}\bigr|^{-2}_g\cdot\bigl|\del_{z_i}\Phi\bigr|^2
=\bigl(1+G_i\bigr)\cdot
 \Bigl(
 \bigl|\del_{x_i}\bigr|^{-2}_g\cdot
 \bigl|\del_{x_i}\Phi\bigr|^2
+\bigl|\del_{y_i}\bigr|^{-2}_g\cdot
 \bigl|\del_{y_i}\Phi\bigr|^2
 \Bigr) \\
\geq
 \bigl(1+G_i\bigr)\cdot
 \Bigl(
 \bigl|\del_{x_i}\bigr|^{-2}_g\cdot
 \bigl|\del_{x_i}\Phi\bigr|^2
 \Bigr)
=(1+H_i)\cdot
 \bigl|\del_{x_i}\Phi\bigr|^2\cdot
 \bigl(-\log|z_i|^2+A\bigr)^2.
\end{multline}
Here we have the estimate
$\bigl|H_i\bigr|=O\bigl(|z_i|^{1/2}\bigr)$.
The estimate is independent of a choice of $\Phi$.

Let us consider the case $i=1$.
We can decompose $(1+F_2)=(1+F_2'')\cdot(1+F_2')$
such that  $F_2'$ is independent of $z_2$,
and the estimate $|F_2''|=O(|z_1|)$ holds.
(See Lemma \ref{lem;a12.31.30}.)
We put $I_1:=
 (1+H_1)\cdot (1+F_3')\cdot (1+F_1)-1$,
and then we have the estimate
$|I_1|=O(|z_1|^{1/2})$.
The estimate is independent of a choice of $\Phi$.
Then we have the following:
\begin{equation}\label{eq;04.1.8.16}
 (1+G_1)\cdot\bigl| \del_{x_1}\Phi \bigr|^2
\cdot \bigl|\del_{x_1}\bigr|^{-2}_g
\cdot\dvol_{g}\\
=\Bigl(
 (1+I_1)\cdot\bigl|\del_{x_1}\Phi\bigr|^2
 \cdot dy_1\cdot dx_1\Bigr)
\cdot(1+F_2')\cdot \frac{dy_2\cdot dx_2}{(2y_2+A)^2}
\end{equation}

\begin{lem}\label{lem;04.1.8.17}
There is a function $J_5$
satisfying $|J_5|=O(|z_1|^{1/2})$
and the following inequality:
\[
 \int_0^{2\pi}
 (1+I_1)\cdot \bigl|\del_{x_1}\Phi\bigr|^2\cdot dx_1
\geq
 \int_0^{2\pi}
 \Bigl(
 \frac{\rho(f_1)^2}{4\pi^2}-J_5
 \Bigr)\cdot dx_1.
\]
\end{lem}
\pf
It follows from Lemma \ref{lem;a12.29.10}.
See the proof of Lemma \ref{lem;04.1.8.10}.
\hfill\qed

\vspace{.1in}

We have the natural projection
of $p:(\Deltabarast)^2\lrarr
 \openclosed{0}{1}^2$.
For any compact subset $K$ of $\openclosed{0}{1}^2$,
we put $\tilde{K}=p^{-1}(K)$.

\begin{lem}\label{lem;04.1.1.1}
There exists an integrable function $J_6$ with respect to
the measure $\dvol_g$ such that the following holds
for any compact subset $K\subset\openclosed{0}{1}^2$:
\[
 \int_{\tilde{K}} (1+I_1)\cdot\bigl|\del_{x_1}\Phi\bigr|^2
 \cdot dx_1\cdot dy_1\cdot
(1+F_2')\cdot\frac{dx_2\cdot dy_2}{(2y_2+A)^2}
\geq
 \int_{\tilde{K}}
 \left(
 \frac{\rho(f_1)^2}{4\pi^2}
\cdot\bigl(-\log|z_1|^2+A\bigr)^2
-J_6
 \right)\cdot
 \dvol_g.
\]
In other words, we have the following inequality:
\begin{multline}\label{eq;04.1.28.10}
\int_{\tilde{K}} 2\cdot (1+G_1)\cdot
 \bigl|\del_{z_1}\Phi\bigr|^2
\cdot\bigl|\del_{z_1}\bigr|_g^{-2}
\cdot\vol_g
\geq
 \int_{\tilde{K}}(1+G_1)\cdot\bigl|\del_{x_1}\Phi\bigr|^2
 \cdot \bigl|\del_{x_1}\bigr|_g^{-2}
 \cdot\dvol_g \\
\geq
  \int_{\tilde{K}}
 \left(
 \frac{\rho(f_1)^2}{4\pi^2}
\cdot\bigl(-\log|z_1|^2+A\bigr)^2
-J_6
 \right)\cdot
 \dvol_g.
\end{multline}
\end{lem}
\pf
We have only to put
$J_6:=
 J_5\cdot (1+F_3')^{-1}\cdot (1+F_1)^{-1}\cdot (1+F_2'')^{-1}
 \cdot(-\log|z_1|+A)^2$.
\hfill\qed

\vspace{.1in}

Similarly, there exists an integrable function $J_7$
with respect to the measure $\dvol_g$ such that
the following holds
for any compact subset $K\subset\openclosed{0}{1}^2$:
\begin{multline}\label{eq;04.1.8.18}
  \int_{\tilde{K}}2\cdot(1+G_2)
 \cdot\bigl|\del_{z_2}\Phi\bigr|^2
 \cdot \bigl|\del_{z_2}\bigr|_g^{-2}
 \cdot\dvol_g
\geq
 \int_{\tilde{K}}(1+G_2)
 \cdot\bigl|\del_{x_2}\Phi\bigr|^2
 \cdot \bigl|\del_{x_2}\bigr|_g^{-2}
 \cdot\dvol_g \\
\geq
  \int_{\tilde{K}}
 \left(
 \frac{\rho(f_2)^2}{4\pi^2}\cdot\bigl(-\log|z_2|^2+A\bigr)^2
-J_7
 \right)\cdot
 \dvol_g.
\end{multline}

\begin{cor}\label{cor;04.1.8.27}
There exists an integrable function $J_8$
with respect to $\dvol_g$
such that the following inequality holds
for any twisted map $\Phi$ and
for any compact subset $K\subset\openclosed{0}{1}$:
\[
 \int_{\tilde{K}}
 e(\Phi)\dvol_g
\geq
 \int_{\tilde{K}}
 \left(
 \sum_{i=1,2}\frac{\rho(f_i)^2}{4\pi^2}\cdot(-\log|z_i|^2+A)^2
  -J_8
 \right)
\cdot \dvol_g.
\]
\end{cor}
\pf
It follows from (\ref{eq;04.1.28.10}),
(\ref{eq;04.1.8.18})
and Lemma \ref{lem;04.1.29.3}.
\hfill\qed

\subsection{On $X-D$}
\label{subsection;04.1.11.5}

%c16.tex

\subsubsection{hermitian metrics of line bundles and neighbourhoods of
 divisors}
\label{subsubsection;04.1.29.350}

Let $X$ be a compact complex surface.
Let $D_i$ $(i=1,\ldots,l)$ be a normal crossing divisor of $X$
such that $D=\bigcup D_i$ is normal crossing.
We have the canonical section $s_i:\nbigo\lrarr\nbigo(D_i)$.

Let $P$ be a point of $D_i\cap D_j$.
We take a sufficiently small neighbourhood $U_P$ of $P$.
We may assume $U_P\cap D_k=\emptyset$ unless $k=i,j$.
We may also assume that
$U_P\cap D_i$ and $U_P\cap D_j$ are holomorphically 
isomorphic to an open disc $\Delta$.
We may assume that we can take a holomorphic trivialization
$e_i$ and $e_j$ of $\nbigo(D_i)$ and $\nbigo(D_j)$ respectively.
The holomorphic functions
$z_i$ and $z_j$ are determined by
$s_i=z_i\cdot e_i$ and $s_j=z_j\cdot e_j$.
Then we may assume that
$\varphi_P=(z_i,z_j):U_P\lrarr\cnum^2$ gives a holomorphic embedding.
We may also assume that $\varphi_P(U_P)\supset\Deltabar^2$.

We have the hermitian metrics $h_{i\,U_P}$ (resp. $h_{j\,U_P}$)
of $\nbigo(D_i)_{|U_P}$ (resp. $\nbigo(D_j)_{|U_P}$)
given by $|e_i|=1$ (resp. $|e_j|=1$).
By shrinking $U_P$ appropriately,
we can take a $C^{\infty}$-hermitian metric $h_i$
of $L_i$ such that $h_{i\,|\,U_P}=h_{i\,U_P}$
for any $P\in D_i\cap D_j$.

We have the metric $dz_i\cdot d\bar{z}_i+dz_j\cdot d\bar{z}_j$ on $U_P$.
Take a hermitian metric of
the tangent bundle $TX$, which is not necessarily Kahler,
such that the following holds:
\begin{itemize}
\item $g_{U_P}=dz_i\cdot d\bar{z}_i+dz_j\cdot d\bar{z}_j$.
 We shrink $U_P$ if it is necessary.
\item
 The hermitian metric $g$ induces the orthogonal decomposition
 $TX_{|D_i}=TD_i\oplus N_{D_i}(X)$.
 We have the natural isomorphism
 $N_{D_i}(X)\simeq \nbigo(D_i)_{|D_i}$.
 Then the restriction of $g$ to $N_{D_i}(X)$
 is same as the restriction of $h_i$ to $\nbigo(D_i)_{|D_i}$.
\end{itemize}

We have the exponential map
$T(X)_{|D_i}\lrarr X$.
Let us consider the restriction
$\exp_i:N_{D_i}(X)\lrarr X$.
For any $R\in\real_{>0}$,
we put as follows:
\[
 N'_{i\,R}:=\bigl\{v\in N_{D_i}(X)\,\big|\,
 h_i(v,v)<R^2\bigr\},
\quad\quad
 N_{i\,R}:=\exp_i\Bigl(N_{i\,R}'\Bigr)\subset X.
\]
It is well known that there exists a positive number $R_0$
such that
$N_{i\,R}'$ and $N_{i\,R}$ are diffeomorphic
for any $R\leq R_0$.

We have the naturally defined projection
$\pi_i:N_{i\,R}'\lrarr D_i$.
In the case $R\leq R_0$,
it induces the projection
$\pi_i:N_{i\,R}\lrarr D_i$.
For any point $P\in D_i$ and for any $R\leq R_0$,
we put $N'_{i\,R\,|\,P}:=\pi_i^{-1}(P)$
and $N_{i\,R\,|\,P}=\exp_i\bigl(N'_{i\,R\,|\,P}\bigr)$.
We have the natural complex structure
of $N'_{i\,R\,|\,P}\simeq\Delta(R)$,
which induces the complex structure of $N_{i\,R\,|\,P}$.
The inclusion $N_{i\,R\,|\,P}\lrarr X$ is not necessarily
holomorphic embedding.
However the following lemma is obtained from our construction.
(Compare the lemma with the condition \ref{condition;04.1.30.1}.)
\begin{lem}\label{lem;04.1.30.2}
The inclusion of the tangent spaces
$T_{(P,O)}N_{i\,R\,|\,P}\lrarr T_{(P,O)}X$
is compatible with their complex structures.
\hfill\qed
\end{lem}

Let $P$ be a point of $D_i\cap D_j$.
The holomorphic function
$z_j$ gives a holomorphic coordinate of $D_i\cap U_P$.
The bundle $N_{D_i}(X)_{U_P\cap D_i}$ is trivialized by $e_i$.

\begin{lem}
There exists a positive constant $R_1$ such that
the map $\exp_i:N_{i\,R}'\lrarr X$
 is given by $(z_j,z_i\cdot e_i)\longmapsto(z_i,z_j)$,
for any $R\leq R_1$.
\end{lem}
\pf
Since the metric on $U_P$ is given by
$dz_i\cdot d\bar{z}_i+dz_j\cdot d\bar{z}_j$,
the claim is clear.
\hfill\qed

\begin{lem}
There exists a positive constant $R_2$
such that the following holds, for any $R\leq R_2$:
\begin{itemize}
\item
 The set of the connected components of $N_{i\,R}\cap N_{j\,R}$
 corresponds bijectively to $D_i\cap D_j$.
\item
 Let $(N_{i\,R}\cap N_{j\,R})_P$ denote the connected components
 of $N_{i\,R}\cap N_{j\,R}$ corresponding to $P$.
 Then $(N_{i\,R}\cap N_{j\,R})_P\subset U_P$.
\item
 $\varphi_P\bigl((N_{i\,R}\cap N_{j\,R})_P\bigr)=
 \Delta(R)^2$.
\end{itemize}
\end{lem}
\pf
It is clear from our construction.
\hfill\qed

\vspace{.1in}
For a positive constant $C$,
we put $\tilde{h}_i=C^{-2}\cdot h_i$,
$\tilde{e}_i=C\cdot e_i$ and $\tilde{z}_i=C^{-1}\cdot z_i$.
If $C$ is taken appropriately,
we may assume that $R_0$, $R_1$ and $R_2$ are larger than $1$.
Hence we may assume $R_i>1$ from the beginning.
In the following,
we use the notation $N_i$ and $N_i'$
to denote $N_{i\,1}$ and $N_{i\,1}'$ respectively.

\subsubsection{The Kahler metric and the decomposition}
\label{subsubsection;04.1.28.20}

We assume $X$ is a Kahler surface with the Kahler metric $g_1$.
\begin{rem}
We do not assume that 
the Kahler metric $g_1$ is not necessarily same as
the hermitian metric of $TX$ used in the subsubsection
{\rm\ref{subsubsection;04.1.29.350}}.
\hfill\qed
\end{rem}

Let $\omega_1$ be the associated Kahler form.
We take a positive constant $A$
such that $|s_i|^2<e^A$ for any $i$.
If a positive constant $C$ is sufficiently large,
the following form $\omega$  also gives a Kahler form
(\cite{cg}):
\[
 \omega=
 C\cdot \omega_1-\sum_i\sqrt{-1}\del\delbar\log\bigl(-\log|s_i|^2+A\bigr).
\]

We put $X^{\circ}=X-\bigcup_i N_i$.
It is $C^{\infty}$ compact submanifold of $X$,
which possibly has a boundary with the corner.

We denote the closure of $N_i$ by $\overline{N}_i$.
It is $C^{\infty}$-submanifold of $X$
with the boundary.
We put $M_P:=\overline{N}_i\cap \overline{N}_j$,
which is holomorphically isomorphic to $\Deltabar^2$.
It is a compact submanifold of $X$
with the boundary and the corner.
We may assume that the norm of the canonical section $s_k$
is constant on $M_P$ unless $k=i,j$.
We can identify $s_k$ and $z_k$ for $k=i,j$ on $M_P$.
We would like to apply the result in the subsection
\ref{subsection;04.1.1.2}.

We put $D_i^{\circ}:=D_i-\bigcup_{j\neq i}(N_j\cap D_i)$,
which is a $C^{\infty}$ compact submanifold of $D_i$
with the boundary.
Let $\del_PD_i$ denote the component of the boundary of $D_i$
contained in $U_P$.

$\overline{N}_i$ can be regarded as $\Deltabar$-bundle over $D_i$
in the $C^{\infty}$-category.
We put $\overline{N}_i^{\circ}:=D_i^{\circ}\times_{D_i}\overline{N}_i$.
We also put $N_i^{\circ}:=D_i^{\circ}\times_{D_i}N_i$.
We put
$\del_0\overline{N}_i^{\circ}
 :=\overline{N}_i^{\circ}-N_i^{\circ}$.
We put $\del_P\overline{N}_i
 :=\del_PD_i^{\circ}\times_{D_i^{\circ}}\overline{N_i}^{\circ}$.

Let $D_i^{\circ}=\coprod U_{i\,j}$ be disjoint unions,
where $U_{i\,j}$ denote subregions of $D_i^{\circ}$,
which are isomorphic to subregions of $\cnum$.
We put $N_{i\,j}^{\circ}:=N_i^{\circ}\times_{D_i^{\circ}}U_{i\,j}$.
The fibration $N_{i\,j}^{\circ}\lrarr U_{i\,j}$
can be trivialized in the $C^{\infty}$-category.
We would like to apply the result in the subsection \ref{subsection;04.1.8.20}
to $N_{i\,j}$,
by putting $-\log|\zeta|^2-b=\log|s_i|^2+A$
(the subsubsection \ref{subsubsection;04.1.28.1})
and $-\log s^2=-\sum_{i=1}^l\bigl(-\log |s_i|^2+A\bigr)$
(the subsubsection \ref{subsubsection;04.1.15.1}).
We remark that the constant $A$ here is different from
the function $A$ in the subsubsection \ref{subsubsection;04.1.28.1}.
We also remark Lemma \ref{lem;04.1.30.2}.

%c16.1.tex

\subsubsection{The construction of maps with the controlled energy}

Let $(E,\nabla)$ be a flat bundle of rank $r$ on $X-D$.
Then we have a homotopy class of twisted maps $X-D\lrarr \PH(r)/\pi_1(X)$.
We take any $C^{\infty}$-twisted map
$F^{\circ}:X^{\circ}\lrarr \PH(r)/\pi_1(X)$,
which corresponds to a hermitian metric $h^{\circ}$
of $(E,\nabla)_{|X^{\circ}}$.

We take a continuous hermitian metric $h_{i}^{\circ}$
of $(E,\nabla)_{|\overline{N}_i^{\circ}-D_i^{\circ}}$,
satisfying the following:
\begin{itemize}
\item
The restrictions of $h_i^{\circ}$ and $h^{\circ}$
to $\del_0\overline{N}_i$ are same.
\item
 The restriction of $h_{i}^{\circ}$
 to $(E,\nabla)_{|N_{i\,|\,P}^{\ast}}$ is
 a tame pure imaginary harmonic metric.
 Here we put $N_{i\,|\,P}^{\ast}:=N_{i\,|\,P}-\{(P,O)\}$,
 and the conformal structure is induced by that of
 $N'_{i\,|\,P}$.
\end{itemize}
(See the subsection \ref{subsection;04.1.8.20}.)
The corresponding twisted map is denoted by $F_i^{\circ}$.

Let $P$ be any point of $D_i\cap D_j$.
We have the continuous hermitian metric  $h_{M_P}$
of $(E,\nabla)_{|M_P\setminus D}$
as in the subsection \ref{subsection;04.1.1.2}.
We may assume that
the restrictions of $h_{M_P}$ and $h^{\circ}$
to $\del\Deltabar\times\del\Deltabar$ are same,
if we modify $h^{\circ}$ appropriately.
Then we also obtain
$h_{M_P\,|\,\del_P\overline{N}_i^{\circ}}
 =h^{\circ}_{i\,|\,\del_P\overline{N}_i^{\circ}}$
due to our construction.

Hence we obtain the continuous hermitian metric $h_0$
of $(E,\nabla)$ on $X-D$,
such that $h_{0\,|\,X^{\circ}}=h^{\circ}$,
$h_{0\,|\,N_i^{\circ}}=h_i^{\circ}$
and $h_{0\,|\,M_P}=h_{M_P}$.
Let $\Psi_{h_0}$ denote the corresponding twisted map.

Let us take a small loop $\gamma_i$ around $D_i$.
Then we obtain the monodromy $\varphi_i$ with respect to $\gamma_i$.
It is easy to see that
the number $\rho(\varphi_i)$ is independent of a choice of $\gamma_i$.
We denote the number by $\rho_i$.

Let $K$ be any compact subset of $X-D$.
\begin{lem}\label{lem;04.1.30.6}
There exist the integrable functions $J_{10}$
on $X-D$ with respect to the measure $\dvol_g$ such that
the following holds, for any compact subset $K\subset X-D$:
\[
 \int_K
 e(\Psi_{h_0})\cdot\dvol_g
\leq
 \int_K\left(
 \sum_i\frac{\rho_i^2}{4\pi^2}\cdot\bigl(-\log|s_i|^2+A\bigr)^2
 +J_{10}
 \right)\cdot\dvol_g.
\]
\end{lem}
\pf
It follows from our construction.
See Corollary \ref{cor;04.1.8.25}
and Corollary \ref{cor;04.1.8.26}.
\hfill\qed

\vspace{.1in}

For any real numbers $R,R_1,R_2$, we put as follows:
\[
 X(R):=\Bigl\{P\in X\,\Big|\,
 \sum -\log|s_i|(P)\leq R\Bigr\},
\quad
 X(R_1,R_2):=\Bigl\{ P\in X\,\Big|\,
 R_1\leq  \sum -\log|s_i|(P)\leq R_2
 \Bigr\}.
\]
When we consider $X(R)$ (resp. $X(R_1,R_2)$),
the number $R$ (resp. $R_1$ and $R_2$) is chosen
such that the boundary of $X(R)$ (resp. $X(R_1,R_2)$)
is $C^{\infty}$.

We can take a $C^{\infty}$-twisted map
$F_N:X(N)\lrarr \PH(r)/\pi_1(X)$,
which approximates $\Psi_{h_0\,|\,X(N)}$
in $L_1^2$-sense sufficiently closely,
such that there exists an integrable function $J$ on $X-N$
such that the following holds for any compact subset $K\subset X(N)$:
\[
 \int_{K}\bigl|e(\Psi_{h_0})-e(F_N)\bigr|\cdot\dvol<
 \int_KJ\cdot\dvol.
\]

\subsubsection{A priori lower bound for energy of any map}

Let $P$ be a point of $D_i\cap D_j$.
We may assume that $|s_k|$ is constant on $M_P$
unless $k=i,j$.

\begin{lem}\label{lem;04.1.30.5}
There exists an integrable function $J_{11}$ on $X-D$
with respect to the measure $\dvol_g$,
such that the following holds:
\begin{itemize}
\item
 For any twisted harmonic map $F:X(R_1,R_2)\lrarr \PH(r)/\pi_1(X)$,
 the following holds:
\[
 \int_{X(R_1,R_2)}e(F)\cdot\dvol_g
\geq
 \int_{X(R_1,R_2)}
 \left(
 \sum_i\frac{\rho_i^2}{4\pi^2}\bigl(-\log|s_i|^2+A\bigr)^2
-J_{11}\right)
 \cdot\dvol_g.
\]
\end{itemize}
\end{lem}
\pf
It follows from Lemma \ref{lem;04.1.8.10}
and Corollary \ref{cor;04.1.8.27}.
\hfill\qed

\section{The existence of tame pure imaginary pluri-harmonic metric}

\label{section;04.1.30.220}

\subsection{Harmonic metric of
 a semisimple flat bundle on a quasi compact Kahler surface}

\label{subsection;04.1.30.221}

%c17.tex

\subsubsection{The existence of harmonic metric}

Let $X$ be a compact Kahler surface.
Let $D$ be a normal crossing divisor of $X$.
Let $(E,\nabla)$ be a semisimple flat bundle on $X-D$.
Let us see the existence of harmonic metric of $(E,\nabla)$.
We use a notation in the subsection \ref{subsection;04.1.11.5}.

Due to the theorem of Hamilton-Schoen-Corlette
(see the proof of Theorem 2.1 of \cite{corlette}),
we can take a twisted harmonic map $\Psi_N:X(N)\lrarr \PH(r)/\pi_1(X)$
satisfying the following:
\[
 \Psi_{N\,|\,\del X(N)}=F_{N\,|\,\del X(N)},
\quad\quad
 \int_{X(N)}e(\Psi_N)\cdot\dvol_g
\leq
 \int_{X(N)}e(F_N)\cdot\dvol_g.
\]
In the case $N>k$,
we have the following inequalities,
due to Lemma \ref{lem;04.1.30.5}:
\begin{multline}
 \int_{X(N)}e(\Psi_N)\cdot\dvol_g
=\int_{X(k)}e(\Psi_N)\cdot\dvol_g+\int_{X(k,N)}e(\Psi_N)\cdot\dvol_g \\
\geq
 \int_{X(k)}e(\Psi_N)\cdot\dvol_g
+\int_{X(k,N)}\left(
 \sum \frac{\rho_i^2}{4\pi^2}\cdot\bigl(-\log|s_i|^2+A\bigr)^2
 -J_{11}
 \right)\cdot\dvol_g.
\end{multline}
We have the following inequality due to Lemma \ref{lem;04.1.30.6}:
\[
 \int_{X(k,N)}\left(
 \sum\frac{\rho_i^2}{4\pi^2}\bigl(-\log|s_i|^2+A\bigr)^2
 \right)\cdot\dvol_g
\geq
 \int_{X(k,N)} \bigl(
 e(\Psi_{h_0})-J_{10}\bigr)\cdot\dvol_{g}
\]
Due to our choice of $F_N$,
we have
$\int_{X(k,N)}e(\Psi_{h_0})\dvol_g\geq \int_{X(k,N)}(e(F_N)-J)\cdot\dvol_g$
for some integrable function $J$,
which is independent of $N$.
Thus we obtain the following inequality:
\[
 \int_{X(k,N)}\left(
 \sum\frac{\rho_i^2}{4\pi^2}\bigl(-\log|s_i|^2+A\bigr)
 \right)\cdot\dvol_g
\geq
 \int_{X(k,N)}(e(F_N)-J_{10}-J)\cdot\dvol_g
\]
Then we obtain the following:
\[
 \int_{X(N)}e(F_N)\dvol_g
\geq
 \int_{X(N)}e(\Psi_N)\cdot\dvol_g
\geq
 \int_{X(k)}e(\Psi_N)\cdot\dvol_g
+\int_{X(k,N)}e(F_N)\cdot\dvol_g
-\int_{X(k,N)}(J_{10}+J_{11}+J)\cdot\dvol_g
\]
Thus we obtain the following:
\begin{multline} \label{eq;04.1.15.2}
 \int_{X(k)}e(\Psi_N)\cdot\dvol_g
\leq
 \int_{X(k)}e(F_N)\cdot\dvol_g
+\int_{X(k,N)}(J_{10}+J_{11}+J)\cdot\dvol_g \\
\leq
 \int_{X(k)}e(\Psi_{h_0})\cdot\dvol_g
+\int_{X(k,N)}(J_{10}+J_{11}+2J)\cdot\dvol_g,\\
\leq
 \int_{X(k)}e(\Psi_{h_0})\cdot\dvol_g
+\int_{X(k,\infty)}(J_{10}+J_{11}+2J)\cdot\dvol_g.
\end{multline}

\begin{lem}
Assume that $(E,\nabla)$ is semisimple.
Then there exists an infinite subset $\vecn_1$ of $\nnum$
such that
the sequence $\bigl\{\Psi_N\,\big|\,N\in\vecn_1\bigr\}$
is $C^{\infty}$-convergent
to a twisted harmonic map $\Psi_{\infty}:X-D\lrarr \PH(r)/\pi_1(X)$.
\end{lem}
\pf
Since we have the estimate of the energy (\ref{eq;04.1.15.2}),
we have only to apply the arguments in
the section 2 in Jost-Yau \cite{jy4} (using semisimplicity)
and \cite{schoen-yau}.
\hfill\qed

\vspace{.1in}

Let $h$ denote the harmonic metric corresponding to
$\Psi_{\infty}$,
and $\theta$ and $\theta^{\dagger}$
denote the corresponding $(1,0)$-form and $(0,1)$-form respectively.
We denote $\Psi_{\infty}$ by $\Psi_h$.

Thus we obtain the harmonic metric $h$
for any semisimple flat bundle $(E,\nabla)$ on $X-D$.
We will show that the harmonic metric $h$ constructed
is, in fact, pluri harmonic (Proposition \ref{prop;04.1.12.30})
and tame pure imaginary (Theorem \ref{thm;04.1.29.305}).

%c17.51.tex

\subsubsection{The decomposition and the energy}
\label{subsubsection;04.1.15.15}

We put $N_{i\,j}^{\circ}(R_1,R_2):=N_{i\,j}^{\circ}\cap X(R_1,R_2)$
and $N_{i\,j}^{\circ}(R):=N_{i\,j}^{\circ}\cap X(R)$.
Due to Corollary \ref{cor;04.1.12.1},
there exist integrable functions $J_{i,j}$
on $N_{i\,j}^{\circ}\setminus D$
with respect to $\dvol_g$,
such that the following holds
for any $R_1<R_2$:
\begin{multline}\label{eq;04.1.15.4}
 \int_{N_{i\,j}^{\circ}(R_1,R_2)}
 e(\Psi_h)\cdot\dvol_g
\geq
 \int_{N_{i\,j}^{\circ}(R_1,R_2)}
 \bigl|
 \del_{\phi}\Psi_h
 \bigr|^2\cdot\bigl|\del_{\phi}\bigr|_g^{-2}\cdot\dvol_g \\
\geq
 \int_{N_{i\,j}^{\circ}(R_1,R_2)}
 \left( 
 \frac{\rho_i^2}{4\pi^2}\cdot\bigl(-\log|s_i|^2+A\bigr)^2
-J_{i,j}
 \right)\cdot\dvol_g.
\end{multline}

We put $M_{P}(R_1,R_2):=M_P\cap X(R_1,R_2)$ and
$M_P(R):=M_P\cap X(R)$.
Due to Corollary \ref{cor;04.1.8.27},
there exists an integrable function $J_P$ on $M_P\setminus D$
on $\dvol_g$,
such that the following inequality for any $R_1<R_2$:
\begin{equation}\label{eq;04.1.15.5}
 \int_{M_P(R_1,R_2)}
 e\bigl(\Psi_h\bigr)\cdot\dvol_g
\geq
 \int_{M_P(R_1,R_2)}\left(
 \sum_{i}\frac{\rho_i^2}{4\pi^2}\cdot\bigl(-\log|s_i|^2+A\bigr)^2
-J_{P}
 \right)\cdot\dvol_g.
\end{equation}

On the other hand,
we obtain the following inequality due to our construction
(see (\ref{eq;04.1.15.2}) and Lemma \ref{lem;04.1.30.6}):
\[
 \int_{X(R)}e(\Psi_h)\cdot\dvol_g
\leq
 \int_{X(R)}\left(
 \sum_i\frac{\rho_i^2}{4\pi^2}\cdot\bigl(-\log|s_i|^2+A\bigr)^2
+J_{10}
 \right)\cdot\dvol_g
+\int_{X(R,\infty)}
 \left(
 J_{10}+J_{11}+2J
 \right)\cdot\dvol_g.
\]
The second term in the right hand side converges to $0$
when $R\to\infty$.
Hence there exists a positive constant $\tilde{C}$,
such that the following holds for any $R$:
\begin{equation}\label{eq;04.1.15.3}
  \int_{X(R)}e(\Psi_h)\cdot\dvol_g
\leq
 \int_{X(R)}
 \sum_i\frac{\rho_i^2}{4\pi^2}\cdot\bigl(-\log|s_i|^2+A\bigr)^2
 \cdot\dvol_g
+\tilde{C}.
\end{equation}
From (\ref{eq;04.1.15.4}), (\ref{eq;04.1.15.5}) and
(\ref{eq;04.1.15.3}),
there exist positive  constants $C_{i\,j}$,
such that the following holds for any $R>0$:
\begin{equation}\label{eq;04.1.15.6}
 \int_{N_{i\,j}^{\circ}(R)} e(\Psi_h)\cdot\dvol_g
\leq
 \int_{N_{i\,j}^{\circ}(R)}
  \frac{\rho_i^2}{4\pi^2}\cdot\bigl(-\log|s_i|^2+A\bigr)^2
 \cdot\dvol_g
+C_{i\,j}.
\end{equation}
Similarly there exist constants $C_P$ for any $P\in D_i\cap D_j$,
such that the following holds for any $R>0$:
\begin{equation}\label{eq;04.1.22.20}
 \int_{M_P^{\circ}(R)} e(\Psi_h)\cdot\dvol_g
\leq
 \int_{M_P^{\circ}(R)}
 \left(
 \frac{\rho_i^2}{4\pi^2}\cdot\bigl(-\log|s_i|^2+A\bigr)^2
+\frac{\rho_j^2}{4\pi^2}\cdot\bigl(-\log|s_j|^2+A\bigr)^2
 \right)\cdot\dvol_g
+C_P.
\end{equation}

%c17.52.tex

\subsubsection{Estimate of the energy on
$N_{i\,j}^{\circ}\setminus D$ and $N_i^{\circ}\setminus D$}
\label{subsubsection;04.1.12.6}

We use the coordinate as in the subsubsection 
\ref{subsubsection;04.1.15.1}.

\begin{lem}\label{lem;04.1.15.11}
The function $ e(\Psi_h)-\bigl|\del_{\phi}\Psi_h\bigr|^2
\cdot\bigl|\del_{\phi}\bigr|^{-2}$
is integrable.

There exists a positive constant $C_{i\,j}$
such that the following inequality holds for any $R_1<R_2$:
\begin{equation}\label{eq;04.1.15.10}
 \int_{N_{i\,j}^{\circ}(R_1,R_2)}
 \bigl|\del_{\phi}\Psi_h\bigr|^2\cdot
  \bigl|\del_{\phi}\bigr|^{-2}\cdot\dvol_g
\leq
 \int_{N_{i\,j}^{\circ}(R_1,R_2)}
 \frac{\rho_i^2}{4\pi^2}\cdot\bigl(-\log |s_i|^2+b\bigr)^2
 \cdot\dvol_g
+C_{i\,j}.
\end{equation}
\end{lem}
\pf
The first claim follows from (\ref{eq;04.1.15.4}),
(\ref{eq;04.1.15.6}) and the positivity
$e(\Psi_h)
 -\bigl|\del_{\phi}\Psi_h\bigr|^2\cdot\bigl|\del_{\phi}\bigr|^{-2}$.
The second claim follows from (\ref{eq;04.1.15.4})
and (\ref{eq;04.1.15.6}).
\hfill\qed

\begin{cor}\label{cor;04.1.30.20}
$\bigl|\del_{\phi}\Psi_h\bigr|^2$ is integrable.
\end{cor}
\pf
It is easy to derive the integrability from (\ref{eq;04.1.15.10})
and the estimate
$\bigl|\del_{\phi}\bigr|^2\sim \bigl(-\log|s_i|^2+b\bigr)^{-2}$.
\hfill\qed

\begin{lem}\label{lem;04.1.30.21}
$\bigl|\del_s\Psi_h\bigr|^2\cdot\bigl|\del_s\bigr|^{-2}$ is integrable.
\end{lem}
\pf
We put
$e_2':=\del_s
-g(\del_s,\del_{\phi})\cdot|\del_{\phi}|^{-2}\cdot\del_{\phi}$.
Then we obtain
$\bigl|\del_s\Psi_h\bigr|^2\leq
 C\cdot\bigl(
 \bigl|d\Psi_h(e_2')\bigr|^2+\bigl|\del_{\phi}\Psi_h\bigr|^2
 \bigr)$.
The function $\bigl|d\Psi_h(e_2')\bigr|^2\cdot\bigl|e_2'\bigr|^{-2}$
is integrable due to Lemma \ref{lem;04.1.15.11}.
Then it is easy to check the claim.
(Use the argument in the subsubsection \ref{subsubsection;04.1.30.10}).
\hfill\qed

\vspace{.1in}

Let $\theta=\theta_1+\theta_2$ be the orthogonal decomposition
such that $\theta_1$ is of the following form:
\[
 \theta_1=\frac{1}{4}
 H^{-1}\cdot dH\bigl(\del_{\phi}-\sqrt{-1}J\cdot\del_{\phi}\bigr)
\cdot \bigl(d\phi-\sqrt{-1}J\cdot d\phi\bigr).
\]
Here $J$ denotes the complex structure of $X$.
Then we have the following:
\[
 \bigl|\theta_1\bigr|^2
=\frac{1}{16}\bigl|
 d\Psi_h\bigl(\del_{\phi}-\sqrt{-1}J\del_{\phi}\bigr)
 \bigr|^2\cdot 2\cdot\bigl|d\phi\bigr|^2
=\frac{1}{8}\Bigl(
\bigl|
 d\Psi_h\bigl(\del_{\phi}\bigr)\bigr|^2
\cdot\bigl|\del_{\phi}\bigr|^{-2}
+\bigl|d\Psi_h
 \bigl(J\del_{\phi}\bigr)
 \bigr|^2
\cdot\bigl|J\del_{\phi}\bigr|^{-2}
\Bigr).
\]
\begin{lem}\label{lem;04.1.22.45}
We have the following estimate:
\[
 J\del_{\phi}=
 -s\del_s+O(s^2)\cdot\del_s+O(s)\cdot\del_{\phi}+O(s)\del_{x_i}.
\]
\end{lem}
\pf
Since we have $dw/w-d\zeta/\zeta=O(1)$,
we have the following with respect to the Kahler metric $g_1$ on $X$:
\[
 J\cdot d\phi=\frac{ds}{s}+O(1),
\quad
 J\cdot \frac{ds}{s}=-d\phi+O(1).
\]
Hence we have
$s^{-1}ds\bigl(J\del_{\phi}\bigr)=(J\cdot s^{-1}ds)(\del_{\phi})=-1+O(s)$
and
$d\phi\bigl(J\del_{\phi}\bigr)=\bigl(J\cdot d\phi\bigr)(\del_{\phi})
=O(s)$.
We also have $J\cdot dx_i=O(1)$ with respect to $g_1$.
Then the claim follows.
\hfill\qed

\begin{lem}\label{lem;04.1.30.50}
$\bigl|d\Psi_h\bigl(J\del_{\phi}\bigr)\bigr|^2\cdot
 \bigl|J\del_{\phi}\bigr|^{-2}$ is integrable.
\end{lem}
\pf
From Lemma \ref{lem;04.1.22.45},
we obtain
$ d\Psi_h\bigl(J\del_{\phi}\bigr)
=O(s)\cdot\del_s\Psi_h
+O(s)\cdot\del_{\phi}\Psi_h
+O(s)\cdot\del_{x_i}\Psi_h$,
and thus we have the following estimate:
\begin{multline}
 \bigl|d\Psi_h(J\del_{\phi})\bigr|^2\cdot
 \bigl|J\del_{\phi}\bigr|^{-2}
=O\Bigl(
 s^2\cdot \bigl|\del_{\phi}\Psi_h\bigr|^2\cdot
 \bigl|\del_{\phi}\bigr|^{-2}
 \Bigr)
+O\Bigl(
 \bigl|\del_{s}\Psi_h\bigr|^2\cdot
 \bigl|\del_s\bigr|^{-2}
 \Bigr)
+O\Bigl(
 s^2\cdot(-\log |s|^2)^2
\cdot\sum_i\bigl|\del_{x_i}\Psi_h\bigr|^2
 \Bigr)\\
=
 O\Bigl(
 s^2\cdot \bigl|\del_{\phi}\Psi_h\bigr|^2\cdot
 \bigl|\del_{\phi}\bigr|^{-2}
 \Bigr)
+O\Bigl(
 \bigl|\del_{s}\Psi_h\bigr|^2\cdot
 \bigl|\del_s\bigr|^{-2}
 \Bigr)
+s\cdot
 O\Bigl(e\bigl(\Psi_h\bigr)\Bigr).
\end{multline}
Then it is easy to derive the integrability of
$s\cdot e\bigl(\Psi_h\bigr)$ from 
from Lemma \ref{lem;04.1.15.11},
and thus the integrability of
$\bigl|d\Psi_h\bigl(J\del_{\phi}\bigr)\bigr|^2\cdot
 \bigl|J\del_{\phi}\bigr|^{-2}$ follows
from Corollary \ref{cor;04.1.30.20}
and Lemma \ref{lem;04.1.30.21}.
\hfill\qed

\begin{lem}\label{lem;04.1.28.60}
There are integrable functions $\breve{J}_{i\,j}$
and a positive constant $\breve{C}_{i\,j}$
such that  the following holds,
for any $R_1<R_2$:
\begin{multline}
 \int_{N_{i\,j}^{\circ}(R_1,R_2)}\left(
 \frac{\rho_i^2}{32\pi^2}\cdot\bigl(-\log |s_i|^2+A\bigr)^2
 -\breve{J}_{i\,j}
 \right)\cdot\dvol_g
\leq 
 \int_{N_{i\,j}^{\circ}(R_1,R_2)}\bigl|\theta_1\bigr|^2\cdot\dvol_g \\
\leq
 \int_{N_{i\,j}^{\circ}(R_1,R_2)}
 \frac{\rho_i^2}{32\pi^2}\cdot\bigl(-\log |s_i|^2+A\bigr)^2
\cdot\dvol_g+\breve{C}_{i\,j}.
\end{multline}
We have the following finiteness on $N_i^{\circ}$:
\begin{equation}\label{eq;04.1.12.5}
 \int_{N_{i}^{\circ}- D_i^{\circ}}
 |\theta_1|^2\cdot\frac{\dvol_g}{\bigl(-\log |s_i|^2+A\bigr)^2}
<\infty.
\end{equation}
The function $\bigl|\theta_2\bigr|^2$ is integrable.
\end{lem}
\pf
The estimate for $\bigl|\theta_1\bigr|^2$ follows from the estimate of
$\bigl|\del_{\phi}\Psi_h\bigr|\cdot\bigl|\del_{\phi}\bigr|^{-2}$
and the integrability of
$\bigl|d\Psi_h(J\del_{\phi})\bigr|\cdot \bigl|J\del_{\phi}\bigr|^{-2}$.
The integrability of $\bigl|\theta_2\bigr|^2$ follows from
the estimates of $\bigl|\theta_1\bigr|^2$ and $e(\Psi_h)$.
\hfill\qed

\begin{cor}\label{cor;04.1.30.30}
We have the finiteness:
\[
 \int_{N_{i}^{\circ}- D_i^{\circ}}
 \bigl|\theta_1\bigr|\cdot\bigl|\theta_2\bigr|
 \frac{\dvol_g}{(-\log |s_i|^2+A)}
<\infty.
\]
\end{cor}
\pf
It follows from the $L^2$-property of $\theta_2$
and $\theta_1\cdot (-\log |s_i|^2+A)^{-1}$.
\hfill\qed

%c17.6.tex

\subsubsection{Estimate of the energy on $M_P\setminus D$}
\label{subsubsection;04.1.23.1}

Let $P$ be a point of $D_i\cap D_j$.
For simplicity, we consider the case $(i,j)=(1,2)$.
Note we have $s_i=z_i$ $(i=1,2)$ on $M_P$.
We use the result in the subsubsection \ref{subsubsection;04.1.22.22}.

From the inequalities (\ref{eq;04.1.28.10}),
(\ref{eq;04.1.8.18}), (\ref{eq;04.1.22.20})
and Lemma \ref{lem;04.1.29.3},
there exist constants $C_i$ such that the following inequalities hold,
for any $R_1<R_2$:
\begin{multline}\label{eq;04.1.29.10}
 \int_{M_P^{\circ}(R_1,R_2)}
 \bigl(1+G_i\bigr)\cdot \bigl|\del_{x_i}\Psi_h\bigr|^2\cdot
 \bigl|\del_{x_i}\bigr|^{-2}_g
 \cdot\dvol_g
\leq
 \int_{M_P^{\circ}(R_1,R_2)}
 2\cdot\bigl(1+G_i\bigr)\cdot \bigl|\del_{z_i}\Psi_h\bigr|^2\cdot
 \bigl|\del_{z_i}\bigr|^{-2}_g
 \cdot\dvol_g \\
\leq
 \int_{M_P^{\circ}(R_1,R_2)}
 \frac{\rho_i^2}{4\pi^2}\cdot\bigl(-\log|z_i|^2+A\bigr)^2
 \cdot\dvol_g
+C_i.
\end{multline}
Here $G_i$ are given in Lemma \ref{lem;04.1.29.3},
and  we have the estimate
$G_i=O\bigl(|z_1|\cdot |z_2|\cdot 
  (-\log|z_1|^2+A)\cdot (-\log|z_2|^2+A)\bigr)$.

\begin{lem}\mbox{{}}\label{lem;04.1.29.21}
The functions $\bigl|\del_{y_i}\Psi_h\bigr|\cdot \bigl|\del_{y_i}\bigr|^{-2}$
$(i=1,2)$
are integrable with respect to $\dvol_g$.
\end{lem}
\pf
The integrability of
$(1+G_i)\cdot\bigl|\del_{y_i}\Psi_h\bigr|\cdot \bigl|\del_{y_i}\bigr|_g^{-2}$
follows from 
(\ref{eq;04.1.28.10}),
(\ref{eq;04.1.8.18})
(\ref{eq;04.1.29.10}) and
the relation
$2\cdot\bigl|\del_{z_i}\bigr|^{-2}\cdot\bigl|\del_{z_i}\Psi_h\bigr|^2
=\bigl|\del_{x_i}\bigr|^{-2}\cdot\bigl|\del_{x_i}\Psi_h\bigr|^2
+\bigl|\del_{y_i}\bigr|^{-2}\cdot\bigl|\del_{y_i}\Psi_h\bigr|^2$.
Then it is easy to derive the lemma.
\hfill\qed

\vspace{.1in}

Recall the argument in the subsubsection \ref{subsubsection;04.1.22.22}.
We put
$M_{P\,1}:=
   \bigl\{(y_1,z_2)\in\real_{\geq\,0}\times\Deltabarast\bigr\}$.
We have the measure
$d\mu_0:=(1+F_2')\cdot(2y_2+A)^{-2}\cdot dy_1\cdot dx_2\cdot dy_2$.
Due to Lemma \ref{lem;04.1.8.17},
there exists a function $J_{100}$ on $M_{P\,1}$
with the estimate $|J_{100}|=O(e^{-y_1/2})$,
such that the following holds:
\[
 \int_{0}^{2\pi}\bigl(1+I_1\bigr)\cdot
 \bigl|\del_{x_1}\Psi_h\bigr|^2\cdot dx_1
\geq
 \frac{\rho(f_1)^2}{2\pi}-J_{100}.
\]
We put
 $M_{P\,1}(R_1,R_2):=
 \bigl\{(y_1,z_2)\in M_{P\,1}\,\big|\,
 R_1\leq y_1+y_2\leq R_2
 \bigr\}$.
From (\ref{eq;04.1.29.10}),
we obtain the following inequality:
\[
 \int_{M_{P\,1}(R_1,R_2)}
 J_{101}\cdot d\mu
\leq \int_{M_{P\,1}(R_1,R_2)}
 J_{100}\cdot d\mu
+C_1
\leq \tilde{C}_1.
\]
Here $C_1$ and $\tilde{C}_1$ are positive constant,
which are independent of $R_1<R_2$.
Then we obtain the integrability of $J_{101}$
with respect to $d\mu$.

We put as follows:
\begin{equation}\label{eq;04.1.29.200}
 J_{102}
:=\int_{0}^{2\pi}\bigl|\del_{x_1}\Psi_h\bigr|^2\cdot dx_1
 -\frac{\rho(f_1)^2}{2\pi},
\quad
 J_{103}:=\int_{0}^{2\pi}
 \bigl|I_1\bigr|\cdot \bigl|\del_{x_1}\Psi_h\bigr|^2\cdot dx_1.
\end{equation}
From the estimate $|I_1|=O\bigl(e^{-y_1/2}\bigr)$
and the integrability of $J_{101}$,
we obtain the integrability of $J_{103}$.
Hence we also obtain the integrability of $J_{102}$.

Let us consider the following function
on $\real_{\geq\,0}^2$:
\[
 \Lambda_i:=
 \int_{0}^{2\pi}\int_{0}^{2\pi}
 \bigl|\del_{x_i}\Psi_h\bigr|^2\cdot dx_1\cdot dx_2
\cdot \bigl(2y_i+A\bigr)^2
-\rho(f_i)^2\cdot\bigl(2y_i+A\bigr)^2.
\]
\begin{lem} \label{lem;04.1.29.20}
The functions $\Lambda_i$ $(i=1,2)$
are integrable with respect to
$(2y_1+A)^{-2}\cdot (2y_2+A)^{-2}\cdot dy_1\cdot dy_2$.
\end{lem}
\pf
The case $i=1$ follows from the integrability
of $J_{102}$ above.
The case $i=2$ can be discussed similarly.
\hfill\qed

\vspace{.1in}

We decompose $\theta=\theta_1+\theta_2$,
where $\theta_i$ are of the form
$f_i\cdot dz_i/z_i=f_i\cdot (\sqrt{-1}dx_i-dy_i)$.
We have the following equality:
\[
 \bigl|f_i\bigr|^2_h
=\frac{1}{16}\cdot
\Bigl(
 \bigl|\del_{x_i}\Psi_h\bigr|^2
+\bigl|\del_{y_i}\Psi_h\bigr|^2
\Bigr).
\]
Let us consider the following functions
on $\real_{\geq\,0}^2$:
\[
 \Phi_i:=\int_0^{2\pi}\int_0^{2\pi}
 \bigl|f_i\bigr|^2_h\cdot dx_1\cdot dx_2
\cdot (2y_i+A)^{2}.
\]
\begin{lem}\label{lem;04.1.29.100}
We have the decomposition
$\Phi_i=(16)^{-1}\cdot\rho(f_i)^2\cdot(2y_i+A)^{2}+J_{4\,i}$,
where $J_{4\,i}$ are integrable with respect to the measure
$(2y_1+A)^{-2}\cdot (2y_2+A)^{-2}\cdot dy_1\cdot dy_2$.
\end{lem}
\pf
It follows from Lemma \ref{lem;04.1.29.21} and Lemma \ref{lem;04.1.29.20}.
\hfill\qed

\subsection{Preliminary integrability}
\label{subsection;04.1.30.400}

%c17.3.tex

\subsubsection{Statement and some reductions}

\begin{lem}\label{lem;04.1.11.6}
$\delbar\theta$ and $\theta^2$ are $L^2$
with respect to the measure $\dvol_g$.
\end{lem}

We will prove Lemma \ref{lem;04.1.11.6} in the next
subsubsections.
Let us take a function $\psi:\real_{\geq\,0}\lrarr \real$
satisfying the following:
\[
  0\leq \psi\leq 1,\quad\quad
  \psi(x)=1\quad (x\leq 1/2),
\quad\quad
 \psi(x)=0\quad(x\geq 2/3).
\]
For any positive number $N$,
we put as follows:
\[
 \chi_N:=\prod_i \psi\Bigl(
 -\frac{\log |s_i|^2}{N}
 \Bigr).
\]

\begin{lem}
When $N$ is sufficiently large,
we have $\chi_N=\psi\bigl(-N^{-1}\cdot\log|s_i|^2\bigr)$
on $N_i^{\circ}$,
and we have
$\chi_N
=\psi\bigl(-N^{-1}\cdot\log|s_i|^2\bigr)
 \cdot
 \psi\bigl(-N^{-1}\cdot\log|s_j|^2\bigr)$
on $M_P$ for $P\in D_i\cap D_j$.
\hfill\qed
\end{lem}

Due to Proposition \ref{prop;04.1.11.1},
we have the following equality:
\begin{multline}
 \int_{X-D}\chi_N\cdot\Bigl(
 C_1\cdot\big|[\theta,\theta]\big|^2_h
+C_2\cdot\big|\delbar\theta\big|_h^2
 \Bigr)\cdot\omega^n
=
 \int_{X-D}
 \delbar\del\chi_N\cdot
 \big\langle\theta,\,\theta\big\rangle\cdot \omega^{n-2} \\
=\sum \int_{N_i^{\circ}}
 \delbar\del\chi_N\cdot
 \big\langle\theta,\,\theta\big\rangle\cdot \omega^{n-2}
+\sum_{P}\int_{M_P}
 \delbar\del\chi_N\cdot
 \big\langle\theta,\,\theta\big\rangle\cdot \omega^{n-2}.
\end{multline}
Since the integrand of the left hand side is positive,
we have only to show that
each term in the right hand side is bounded
independently of $N$.
We will check such boundedness in the
subsubsections
\ref{subsubsection;04.1.22.25}
and \ref{subsubsection;04.1.22.26}.

\subsubsection{On $N_{i,j}^{\circ}\setminus D$
and $N_i^{\circ}\setminus D$}
\label{subsubsection;04.1.22.25}

We use the results in the subsubsection \ref{subsubsection;04.1.12.6}.
On $N_{i,j}^{\circ}$, we have the following equality:
\[
 \del\delbar\chi_N=
 \frac{1}{N^2}\cdot
 \psi''\Bigl(-\frac{\log |s_i|^2}{N}\Bigr)
 \cdot\del\log |s_i|^2\wedge\delbar\log |s_i|^2
-\frac{1}{N}
 \psi'\Bigl(\frac{-\log |s_i|^2}{N}\Bigr)
 \cdot\del\delbar\log |s_i|^2.
\]
We put $\tau=-\del\delbar\log|s_i|^2$.
It is a $C^{\infty}$-closed form on $X$, and it gives
the first Chern class of $\nbigo(D_i)$ in the cohomology level.
We put as follows:
\[
 \mu=J\cdot d\phi+\sqrt{-1}d\phi,\quad\quad
 G_1:=\del\log |s_i|^2-\mu.
\]
Then we have $\bigl|G_1\bigr|=O(1)$
with respect to the Kahler form $g_1$ of $X$.

Let $\theta=\theta_1+\theta_2$ be the orthogonal decomposition
as in the subsubsection \ref{subsubsection;04.1.12.6}.
Recall that $\theta_1$ is of the form $f\cdot \mu$.
We have the following:
\begin{multline}
 \bigl\langle \theta,\theta
 \bigr\rangle
\cdot
 \del\log |s_i|^2\cdot\delbar\log |s_i|^2 \\
=\bigl\langle\theta_1,\theta_1\bigr\rangle
 \cdot G_1\cdot \bar{G}_1
+\bigl\langle\theta_1,\theta_2\bigr\rangle
 \cdot G_1\cdot\left(\bar{\mu}+\bar{G}_1\right)
+\bigl\langle\theta_2,\theta_1\bigr\rangle
\cdot\left(\mu+G_1\right)\cdot\bar{G}_1
+\bigl\langle\theta_2,\theta_2\bigr\rangle
\cdot\left(\mu+G_1\right)
\cdot\left(\bar{\mu}+\bar{G}_1\right).
\end{multline}
It is easy to check the following estimates:
\[
 \bigl\langle
 \theta_1,\theta_1
 \bigr\rangle\cdot G_1\cdot \bar{G}_1
\cdot \bigl(-\log|s_i|^2+A\bigr)^{-2}
=O\Bigl(\bigl|\theta_1\bigr|^2 \cdot \bigl(-\log|s_i|^2+A\bigr)^{-2}
 \cdot\dvol_g
\Bigr),
\]
\[
 \bigl\langle
 \theta_1,\theta_2
 \bigr\rangle\cdot G_1\cdot
 \left(\bar{\mu}+\bar{G}_1\right)
\cdot \bigl(-\log|s_i|^2+A\bigr)^{-2}
=O\Bigl(
 \bigl|\theta_1\bigr|\cdot\bigl(-\log |s_i|^2+A\bigr)^{-1}
\cdot\bigl|\theta_2\bigr|\cdot\dvol_g
 \Bigr),
\]
\[
 \bigl\langle
 \theta_2,\theta_1
 \bigr\rangle\cdot
 \left(\mu+G_1\right)\cdot\bar{G}_1
\cdot\bigl(-\log|s_i|^2+A\bigr)^{-2}
=O\Bigl(
 \bigl|\theta_2\bigr|\cdot\bigl|\theta_1\bigr|
\cdot\bigl(-\log |s_i|^2+A\bigr)^{-1}\cdot\dvol_g
 \Bigr),
\]
\[
 \bigl\langle\theta_2,\theta_2\bigr\rangle
\cdot\left(\mu+G_1\right)
\cdot\left(\bar{\mu}+\bar{G}_1\right)
\cdot\bigl(-\log|s_i|^2+A\bigr)^{-2}
=O\Bigl(
 \bigl|\theta_2\bigr|^2
\cdot\dvol_g
 \Bigr).
\]
The right hand sides are integrable due to
Lemma \ref{lem;04.1.28.60}
and Corollary \ref{cor;04.1.30.30}.
We have the boundedness:
\[
 \frac{\bigl(-\log |s_i|^2+A\bigr)^2}{N^2}
\cdot\psi''\Bigl(\frac{-1}{N}\log |s_i|^2\Bigr)\leq C.
\]
Since the support of the function
$\psi''\bigl(-N^{-1}\cdot\log |s_i|^2\bigr)$ goes to infinity
when $N\to\infty$,
we obtain the following convergence:
\begin{equation}\label{eq;04.1.12.7}
 \lim_{N\to\infty}
 \int_{N_i^{\circ}}
 \bigl\langle
 \theta,\theta
 \bigr\rangle\cdot\psi''\left(
 \frac{-1}{N}\log |s_i|^2
 \right)\cdot\frac{1}{N^2}\cdot
 \del\log |s_i|^2\wedge \delbar\log |s_i|^2=0.
\end{equation}

We decompose as
$\bigl\langle\theta,\theta\bigr\rangle
=\bigl\langle\theta_1,\theta_1\bigr\rangle
+\bigl\langle\theta_1,\theta_2\bigr\rangle
+\bigl\langle\theta_2,\theta_1\bigr\rangle
+\bigl\langle\theta_2,\theta_2\bigr\rangle$.
It is easy to check the following estimates:
\[
 \bigl\langle\theta_1,\theta_2\bigr\rangle\cdot\tau
 \cdot\bigl(-\log|s_i|^2+A\bigr)^{-1}
=-\overline{\bigl\langle\theta_2,\theta_1\bigr\rangle}\cdot\tau
 \cdot\bigl(-\log|s_i|^2+A\bigr)^{-1}
=O\Bigl(\bigl|\theta_1\bigr|\cdot\bigl|\theta_2\bigr|
\cdot\bigl(-\log|s_i|^2+A\bigr)^{-1}
 \cdot\dvol_g\Bigr),
\]
\[
 \bigl\langle\theta_2,\theta_2\bigr\rangle\cdot\tau
\cdot\bigl(-\log|s_i|^2+A\bigr)^{-1}
=O\Bigl(
 \bigl|\theta_2\bigr|^2\cdot\dvol_g
\cdot\bigl(-\log|s_i|^2+A\bigr)^{-1}
 \Bigr).
\]
The right hand sides are integrable due to
Lemma \ref{lem;04.1.28.60}
and Corollary \ref{cor;04.1.30.30}.
We have the boundedness
of the functions
$-N^{-1}\cdot\log |s_i|^2\cdot
 \psi'\bigl(-N^{-1}\cdot\log |s_i|^2\bigr)$,
independently of $N$,
and the supports of the functions go to infinity.
Hence we obtain the following convergence:
\begin{equation}\label{eq;04.1.12.71}
\lim_{N\to\infty}
 \int_{N_i^{\circ}}
\Bigl(
 \bigl\langle\theta_2,\theta_1\bigr\rangle
+\bigl\langle\theta_1,\theta_2\bigr\rangle
+\bigl\langle\theta_2,\theta_2\bigr\rangle
\Bigr)
\cdot\frac{1}{N}\cdot
\psi'\Bigl(\frac{-1}{N}\log |s_i|^2\Bigr)
\wedge\tau=0.
\end{equation}

To estimate the remained term,
we use the coordinate as in the subsubsection \ref{subsubsection;04.1.15.1}.
We put $\tau_0:=\tau_{|D_i^{\circ}}$.
Let $\pi_i:N_i^{\circ}\lrarr D_i^{\circ}$ be the projection.
We have the decomposition $\tau=\pi_i^{\ast}\tau_0+\tau_1$.
Then we have the estimate
$\bigl|\tau_1\bigr|=O(s)$.
Recall that $\theta_1$ is of the form $f\cdot \mu$.
Then we have the following:
\begin{multline}
 \bigl\langle\theta_1,\theta_1\bigr\rangle\wedge\tau
=\bigl\langle\theta_1,\theta_1\bigr\rangle\wedge\pi_i^{\ast}\tau_0
+\bigl\langle\theta_1,\theta_1\bigr\rangle\wedge\tau_1 \\
=-2\sqrt{-1}\cdot |f|^2\cdot\frac{ds\cdot d\phi}{s}
 \wedge \pi_i^{\ast}\tau_0
+O\Bigl(
 |f|^2\cdot ds\cdot d\phi\wedge\pi_i^{\ast}\tau_0
+
 \bigl|\theta_1\bigr|^2\cdot s^2\cdot\dvol_g
 \Bigr).
\end{multline}
Since $|\theta_1\bigr|^2\cdot \bigl(-\log s-b_1\bigr)^{-2}\cdot\dvol_g$
is integrable due to Lemma \ref{lem;04.1.28.60},
the last term does not contribute to the limit when $N\to\infty$.
We have the equality:
\[
 |f|^2=
\frac{1}{16}\Bigl(
 \bigl|d\Psi_h\bigl(\del_{\phi}\bigr)\bigr|^2
+\bigl|d\Psi_h\bigl(J\del_{\phi}\bigr)\bigr|^2
 \Bigr)
\cdot (1+O(s)).
\]
We have only to consider
$\bigl|\del_{\phi}\Psi_h\bigr|^2/16$
instead of $|f|^2$, due to the integrability of the other terms.

\vspace{.1in}
Let us consider the following integral:
\[
  \int_{N_{i\,j}^{\circ}\setminus D}
 \frac{1}{N}\psi'\Bigl(\frac{-1}{N}\log |s_i|^2\Bigr)
\cdot \bigl|\del_{\phi}\Psi_h\bigr|^2\cdot
  \frac{ds\cdot d\phi}{s}\wedge\pi_i^{\ast}\tau_0
=\int_{W_{i\,j}\times\openclosed{0}{1}}
 \frac{1}{N}\psi'\Bigl(\frac{-1}{N}\log |s_i|^2\Bigr)
\cdot \Phi\frac{ds}{s}\wedge\pi^{\ast}\tau_0.
\]
Here we put as follows:
\[
 \Phi:=\int_{0}^{2\pi}\bigl|\del_{\phi}\Psi_h\bigr|^2\cdot d\phi.
\]
Due to Corollary \ref{cor;04.1.12.1},
there exists an integrable function $J_{25}$ on
$W_{i\,j}\times\openclosed{0}{1}$ with respect to
$s^{-1}ds\cdot \dvol_{W_{i\,j}}$,
such that the following holds:
\[
 \frac{\rho_i^2}{2\pi}
 -J_{25}
 \leq
 \Phi.
\]
There exists a positive constant $C$
such that the following holds, for any $0<a<b<1$:
\begin{equation}\label{eq;04.1.28.50}
 \int_{W_{i\,j}\times[a,b]}
 \Phi\cdot \frac{ds}{s}\cdot\dvol_{W_{i\,j}}
\leq
 \int_{W_{i\,j}\times[a,b]}
 \left( \frac{\rho_i^2}{2\pi}
 \right)\cdot\frac{ds}{s}\cdot\dvol_{W_{i\,j}}
+C.
\end{equation}
We put $J_{27}:=\Phi-(2\pi)^{-1}\cdot\rho_i^2+J_{25}$,
which is a positive function.
It is easy to derive the integrability of $J_{27}$
from (\ref{eq;04.1.28.50}).
Namely,
we have the decomposition
$ \Phi=\frac{\rho_i^2}{2\pi}+J_{28}$,
where $J_{28}$ is integrable
with respect to $s^{-1}\cdot ds\wedge \dvol_{W_{i\,j}}$.
We have the following:
\begin{multline}
 \int_{W_{i\,j}\times\openclosed{0}{1}}
 \frac{1}{N}\cdot  \psi'\Bigl(\frac{-1}{N}\log |s_i|^2\Bigr)
\cdot \Phi\cdot\frac{ds}{s}\wedge \pi^{\ast}\tau_0
=\int_{W_{i\,j}\times\openclosed{0}{1}}
 \frac{1}{N}
 \psi'\Bigl(\frac{-1}{N}\log |s_i|^2\Bigr)\cdot
 \left(
\frac{\rho_i^2}{2\pi}
+J_{28}
\right)
\cdot \frac{ds}{s}\wedge \pi^{\ast}\tau_0 \\
=\frac{\rho_i^2}{4\pi}\int_{W_{i\,j}}\tau_0
+\int_{W_{i\,j}\times\openclosed{0}{1}}
 \frac{1}{N}
 \psi'\Bigl(\frac{-1}{N}\log|s_i|^2\Bigr)\cdot
 J_{28}
\cdot\frac{ds}{s}\wedge\pi^{\ast}\tau_0.
\end{multline}
The second term converges to $0$ when $N\to\infty$.

Let us consider the limit of the integrals over $N_i^{\circ}$.
Then we obtain the following convergence:
\[
 \lim_{N\to\infty}
 \int_{N_i^{\circ}}
 \bigl\langle\theta_1,\theta_1\bigr\rangle\cdot
 \frac{1}{N}\cdot\psi'\Bigl(\frac{-1}{N}\log|s_i|^2\Bigr)\wedge\tau
=C\cdot\rho_i^2\cdot (D_i,D_i).
\]
Here $(D_i,D_i)$ denotes the self intersection number of $D_i$
in $X$,
and $C$ denotes the constant which is independent
of $X$, $D$, $(E,\nabla)$.

Thus we obtain the following convergence:
\begin{equation}\label{eq;04.1.30.500}
\lim_{N\to\infty}
 \int_{N_i^{\circ}}\bigl\langle\theta,\theta\bigr\rangle
 \wedge\del\delbar\chi_N
=C\cdot\rho_i^2\bigl(D_i,D_i\bigr).
\end{equation}
In particular, we obtain the uniform boundedness
on the integrals over $N_i^{\circ}$.

%c17.31.tex

\subsubsection{On $M_P\setminus D$}
\label{subsubsection;04.1.22.26}

Let $P$ be a point of $D_i\cap D_j$.
For simplicity, we consider the case $(i,j)=(1,2)$.
Note we have $s_i=z_i$ $(i=1,2)$ on $M_P$.
We have the following equality:
\begin{multline}
 \del\delbar\chi_N=
\psi''\bigl(-N^{-1}\log|z_1|^2\bigr)
\cdot\psi\bigl(-N^{-2}\log|z_2|^2\bigr)
\cdot N^{-2}\cdot \frac{dz_1\cdot d\bar{z}_1}{|z_1|^2} \\
+ 
\psi'\bigl(-N^{-1}\log|z_1|^2\bigr)
\cdot\psi'\bigl(-N^{-2}\log|z_2|^2\bigr)
\cdot N^{-2}\cdot \frac{dz_1\cdot d\bar{z}_2}{z_1\cdot \bar{z}_2} \\
+
\psi'\bigl(-N^{-1}\log|z_1|^2\bigr)
\cdot\psi'\bigl(-N^{-2}\log|z_2|^2\bigr)
\cdot N^{-2}\cdot \frac{dz_2\cdot d\bar{z}_1}{z_2\cdot \bar{z}_1}\\
+ 
\psi\bigl(-N^{-1}\log|z_1|^2\bigr)
\cdot\psi''\bigl(-N^{-2}\log|z_2|^2\bigr)
\cdot N^{-2}\cdot \frac{dz_2\cdot d\bar{z}_2}{|z_2|^2}.
\end{multline}

We decompose $\theta=\theta_1+\theta_2$,
where $\theta_i$ are of the form $f_i\cdot dz_i/z_i$.
Then we have the following equality:
\begin{multline}
 \bigl\langle\theta,\theta\bigr\rangle
\wedge\del\delbar\chi_N 
=\bigl\langle\theta_1,\theta_1\bigr\rangle
\cdot
 \psi\bigl(-N^{-1}\log|z_1|^2 \bigr)
\cdot
 \psi''\bigl(-N^{-1}\log|z_2|^2\bigr)
\cdot
 \frac{1}{N^2}\frac{dz_2\cdot d\bar{z}_2}{|z_2|^2} \\
+\bigl\langle\theta_1,\theta_2\bigr\rangle
\cdot
 \psi'\bigl(-N^{-1}\log|z_1|^2 \bigr)
\cdot
 \psi'\bigl(-N^{-1}\log|z_2|^2\bigr)
\cdot
 \frac{1}{N^2}\frac{dz_2\cdot d\bar{z}_1}{z_2\cdot\bar{z}_1}\\
+\bigl\langle\theta_2,\theta_1\bigr\rangle
\cdot
 \psi'\bigl(-N^{-1}\log|z_1|^2 \bigr)
\cdot
 \psi'\bigl(-N^{-1}\log|z_2|^2\bigr)
\cdot
 \frac{1}{N^2}\frac{dz_1\cdot d\bar{z}_2}{z_1\cdot\bar{z}_2} \\
+\bigl\langle\theta_2,\theta_2\bigr\rangle
\cdot
 \psi\bigl(-N^{-1}\log|z_1|^2 \bigr)
\cdot
 \psi''\bigl(-N^{-1}\log|z_2|^2\bigr)
\cdot
 \frac{1}{N^2}\frac{dz_1\cdot d\bar{z}_1}{|z_1|^2}.
\end{multline}

We use the real coordinate $z_i=\exp\bigl(\sqrt{-1}x_i-y_i\bigr)$
as usual,
and we use the results in the subsubsection \ref{subsubsection;04.1.23.1}.
Let us estimate the following integral:
\begin{multline}
 \int_{M_P}\bigl\langle\theta_1,\theta_1\bigr\rangle
\cdot\psi\bigl(2N^{-1}y_1\bigr)\cdot\psi''\bigl(2N^{-1}y_2\bigr)\cdot
\frac{1}{N^2}\frac{dz_2\cdot d\bar{z}_2}{|z_2|^2} \\
=\int_{M_P}\bigl|f_1\bigr|^2_h\cdot
 \psi\bigl(2N^{-1}y_1\bigr)\cdot\psi''\bigl(2N^{-1}y_2\bigr)\cdot
\frac{1}{N^2}dx_1\cdot dy_1\cdot dx_2\cdot dy_2 \\
=\int_{\real_{\geq\,0}^2}
 \Phi_1\cdot \frac{1}{N^2}
\cdot
 \psi\bigl(2N^{-1}y_1\bigr)\cdot\psi''\bigl(2N^{-1}y_2\bigr)\cdot
 \frac{dy_1}{(2y_1+A)^2}
\cdot dy_2.
\end{multline}
Due to Lemma \ref{lem;04.1.29.100},
the right hand side can be rewritten as follows:
\begin{multline} \label{eq;04.1.29.101}
\int_{\real_{\geq\,0}^2}
 \rho_1^2\cdot
 \frac{1}{N^2}
\cdot
 \psi\bigl(2N^{-1}y_1\bigr)\cdot\psi''\bigl(2N^{-1}y_2\bigr)\cdot
 dy_1\cdot dy_2 \\
+\int_{\real_{\geq\,0}^2}
 J_{4\,1}\cdot
 \frac{(2y_2+A)^2}{N^2}
\cdot
 \psi\bigl(2N^{-1}y_1\bigr)\cdot\psi''\bigl(2N^{-1}y_2\bigr)\cdot
\frac{dy_1}{(2y_1+A)^2}\frac{dy_2}{(2y_2+A)^2}.
\end{multline}
Since $\psi'(0)=\psi'(\infty)=0$,
the first term in (\ref{eq;04.1.29.101}) vanishes.
We have the boundedness of the functions
$N^{-2}\cdot (2y_2+A)^{2}\cdot\psi''\bigl(2N^{-1}y_2\bigr)$
independently of $N$,
and the supports of the functions go to infinity
when $N\to\infty$.
Thus the second term converges to $0$
when $N\to\infty$.
Namely we obtain the following convergence:
\[
 \lim_{N\to\infty}
 \bigl\langle
 \theta_1,\theta_1
 \bigr\rangle\cdot
 \psi\bigl(-N^{-1}\log|z_1|^2\bigr)
\cdot\psi''\bigl(-N^{-1}\log|z_2|^2\bigr)
\cdot\frac{1}{N^2}\cdot\frac{dz_2\cdot d\bar{z}_2}{|z_2|^2}
=0.
\]
Similarly, we obtain the following:
\[
 \lim_{N\to\infty}
 \bigl\langle
 \theta_2,\theta_2
 \bigr\rangle\cdot
 \psi''\bigl(-N^{-1}\log|z_1|^2\bigr)
\cdot\psi\bigl(-N^{-1}\log|z_2|^2\bigr)
\cdot\frac{1}{N^2}\cdot\frac{dz_1\cdot d\bar{z}_1}{|z_1|^2}
=0.
\]

It is easy to check the following estimate:
\[
 \bigl\langle
 \theta_1,\theta_2
 \bigr\rangle\cdot\frac{dz_2}{z_2}\cdot\frac{d\bar{z}_1}{\bar{z}_1}
=O\Bigl(
 \bigl|f_1\bigr|\cdot\bigl|f_2\bigr|
\cdot
dx_1\cdot dy_1\cdot dx_2\cdot dy_2
 \Bigr).
\]
We put 
$ M_N:=\bigl\{(z_1,z_2)\,\big|\,
 N/2\leq -\log|z_i|^2\leq 2\cdot N/3
 \bigr\}$.
Then the support of
the function $\prod_{i=1,2}\psi'\bigl(-N^{-1}\log|z_i|^2\bigr)$
is contained in $M_N$.
Hence we obtain the following estimate:
\begin{multline}
 \int\bigl\langle\theta_1,\theta_2\bigr\rangle
 \cdot\psi'\bigl(-N^{-1}\log|z_1|^2\bigr)
 \cdot\psi'\bigl(-N^{-1}\log|z_2|^2\bigr)
 \cdot \frac{1}{N^2}\cdot
 \frac{dz_2\cdot d\bar{z}_1}{z_2\cdot \bar{z}_1}
=O\left(
 \int_{M_N}\bigl|f_1\bigr|\cdot \bigl|f_2\bigr|\cdot
 \frac{d\mu_1}{N^2}
 \right) \\
=O\left(
 \left(\int_{M_N}\bigl|f_1\bigr|^2\frac{d\mu_1}{N^2}
 \right)^{1/2}
\cdot
 \left(\int_{M_N}\bigl|f_2\bigr|^2\frac{d\mu_1}{N^2}\right)^{1/2}
 \right).
\end{multline}
Here we put $d\mu_1=dx_1\cdot dy_1\cdot dx_2\cdot dy_2$.
Due to Lemma \ref{lem;04.1.29.100},
we have the following inequalities:
\[
 \int_{M_N}\bigl|f_1\bigr|^2\cdot \frac{d\mu}{N^2}
\leq
\int_{N/2}^{2N/3}dy_1
\int_{N/2}^{2N/3}dy_2\cdot
 \frac{\rho_1^2}{16\cdot N^2}
+
\int_{N/2}^{2N/3}\frac{dy_1}{(2y_1+A)^2}
\int_{N/2}^{2N/3}\frac{dy_2}{N^2}\cdot
J_{4\,i}.
\]
The second term in the right hand side converges to $0$
when $N\to\infty$,
due to the integrability of $J_{4\,i}$
with respect to the measure
$(2y_1+A)^{-2}\cdot(2y_2+A)^{-2}\cdot dy_1\cdot dy_2$.
The first term is as follows:
\[
\frac{\rho_i^2}{16N^2}\cdot 
 \frac{N^2}{6^2}
=\frac{\rho_i^2}{16\times 6^2}.
\]
Hence we obtain the boundedness of the following integrals,
independently of $N$:
\[
 \int \bigl\langle\theta_1,\theta_2\bigr\rangle\cdot
 \psi'\bigl(-N^{-1}\log|z_1|^2\bigr)
\cdot
 \psi'\bigl(-N^{-1}\log|z_2|^2\bigr)\cdot
\frac{1}{N^2}
 \frac{dz_2\cdot d\bar{z}_1}{z_2\cdot \bar{z}_1}.
\]
Similarly we obtain the boundedness of the following,
independent of $N$:
\[
 \int \bigl\langle\theta_2,\theta_1\bigr\rangle\cdot
  \psi'\bigl(-N^{-1}\log|z_1|^2\bigr)
\cdot
 \psi'\bigl(-N^{-1}\log|z_2|^2\bigr)\cdot
\frac{1}{N^2}
 \frac{dz_1\cdot d\bar{z}_2}{z_1\cdot \bar{z}_2}.
\]
Hence we obtain the boundedness of 
$\int_{M_P} \bigl\langle\theta,\theta\bigr\rangle\cdot\del\delbar\chi_N$,
independently of $N$.

Thus the proof of Lemma \ref{lem;04.1.11.6} is accomplished.
\hfill\qed

\subsection{Pluri-harmonicity}
\label{subsection;04.1.30.401}

%c17.4.tex

\subsubsection{Statement and some reduction}

\begin{prop}\label{prop;04.1.12.30}
The harmonic metric $h$ is pluri-harmonic.
\end{prop}

We use the Bochner type formula in Proposition \ref{prop;04.1.22.40}.
To show Proposition \ref{prop;04.1.12.30},
we have only to show the following vanishing of the limit:
\[
 \int_{X-D} d \bigl\langle\delbar\theta,
 \theta-\theta^{\dagger}\bigr\rangle
=\lim_{N\to\infty}
 \int_{X-D} \chi_N\cdot
 d \bigl\langle\delbar\theta,
 \theta-\theta^{\dagger}\bigr\rangle
=\lim_{N\to\infty}
 \int_{X-D} \!\!-d\chi_N\wedge
 \bigl\langle\delbar\theta,
 \theta-\theta^{\dagger}\bigr\rangle=0.
\]
We have only to see the vanishing of the limit of the integrals
over $N_i^{\circ}-D_i^{\circ}$ and $M_P\setminus D$.

\subsubsection{On $N^{\circ}_{i\,j}\setminus D$ and
$N^{\circ}_i\setminus D$}

We use the orthogonal decomposition $\theta=\theta_1+\theta_2$
as in the subsubsection \ref{subsubsection;04.1.22.25}.
Namely we have the following:
\[
 \theta_1=\frac{1}{4}
 H^{-1}\cdot dH\bigl(\del_{\phi}-\sqrt{-1}J\cdot\del_{\phi}\bigr)
 \cdot \bigl(d\phi-\sqrt{-1}J\cdot d\phi\bigr).
\]
Hence we have the following:
\[
 \theta_1^{\dagger}=\frac{1}{4}
 H^{-1}\cdot dH\bigl(\del_{\phi}+\sqrt{-1}J\cdot\del_{\phi}\bigr)
\cdot\bigl(d\phi+\sqrt{-1}J\cdot d\phi\bigr).
\]
Therefore we have the following:
\[
 \theta_1-\theta_1^{\dagger}
=-\frac{\sqrt{-1}}{2}
 \Bigl(
 H^{-1}dH\bigl(J\del_{\phi}\bigr)\cdot d\phi
+H^{-1}dH\bigl(\del_{\phi}\bigr)\cdot\bigl(s^{-1}ds +G_2\bigr)
 \Bigr).
\]
Here we put $G_2:=Jd\phi-d \log s$.
We also put $G_3:= d\log |s_i|-d\log s$.
We have the estimates $\bigl|G_i\bigr|=O(1)$ $(i=2,3)$
with respect to the Kahler metric $g_1$ of $X$.
We also have the following equality:
\[
 d\chi_N
=-\psi'\bigl(-N^{-1}\log |s_i|^2\bigr)
\cdot \frac{2}{N}\cdot d\log |s_i|.
\]
Therefore we have the following equality:
\begin{multline}\label{eq;04.1.12.15}
d\chi_N\wedge
 \bigl\langle
 \delbar\theta,\theta-\theta^{\dagger}
 \bigr\rangle
=
-\psi'\bigl(-N^{-1}\log |s_i|^2\bigr)\cdot
 \frac{\sqrt{-1}}{N}
 \left(
 \frac{ds}{s}+G_3
\right)
  \wedge
 \bigl\langle
 \delbar\theta, H^{-1}dH\bigl(J\del_{\phi}\bigr)\cdot d\phi
 \bigr\rangle \\
-\psi'\bigl(-N^{-1}\log s^2\bigr)\cdot
 \frac{\sqrt{-1}}{N}\left(
 \frac{ds}{s}+G_3\right)
 \wedge
 \bigl\langle
 \delbar\theta,H^{-1}dH\bigl(\del_{\phi}\bigr)\cdot G_2
 \bigr\rangle\\
-\psi'\bigl(-N^{-1}\log s^2\bigr)\cdot
 \frac{\sqrt{-1}}{N}
 \cdot G_3
 \wedge
 \bigl\langle
 \delbar\theta,H^{-1}dH\bigl(\del_{\phi}\bigr)\cdot s^{-1}ds
 \bigr\rangle.
\end{multline}
Due to Lemma \ref{lem;04.1.22.45},
the first term is dominated by the following,
which is integrable (see the proof of Lemma \ref{lem;04.1.30.50}):
\[
 \bigl|\delbar\theta\bigr|
\cdot
\Bigl(
 \bigl|\del_s\Psi_h\bigr|\cdot\bigl|\del_s\bigr|^{-2}
+s\cdot\bigl|\del_{\phi}\Psi_h\bigr|\cdot \bigl|\del_{\phi}\bigr|^{-2}
+s^{1/2}\sum \bigl|\del_{x_i}\Psi_h\bigr|^2\cdot\bigl|\del_{x_i}\bigr|^2
\Bigr)
 \cdot\dvol_g.
\]
The second term of (\ref{eq;04.1.12.15} )is dominated by
$\bigl|\delbar\theta\bigr|\cdot
 \bigl|\del_{\phi}\Psi_h\bigr|\cdot\dvol_g$,
which is also integrable
(Corollary \ref{cor;04.1.30.20} and 
 Lemma \ref{lem;04.1.11.6}).
The third term of (\ref{eq;04.1.12.15}) can be 
dominated similarly.
Hence the right hand side of (\ref{eq;04.1.12.15})
converges $0$ when $N\to\infty$,
because the support of $\psi'\bigl(-N^{-1}\cdot\log s^2\bigr)$ goes 
to infinity.

%c17.41.tex

\subsubsection{On $M_P\setminus D$}

We have
$\theta-\theta^{\dagger}
=\theta_1-\theta_1^{\dagger}
+\theta_2-\theta_2^{\dagger}$.
We use the real coordinate $z_i=\exp\bigl(\sqrt{-1}x_i-y_i\bigr)$,
and we use the result in the subsubsection \ref{subsubsection;04.1.23.1}.
We have the following formula:
\begin{equation}\label{eq;04.1.29.150}
 \theta_i-\theta_i^{\dagger}
=\frac{\sqrt{-1}}{2}
 \left(
 h^{-1}\frac{\del h}{\del y_i}\cdot dx_i
-h^{-1}\frac{\del h}{\del x_i}\cdot dy_i
 \right).
\end{equation}
We also have the following on $M_P$:
\[
 d\chi_N=
-\psi'\bigl(2N^{-1}y_1\bigr)\cdot \psi\bigl(2N^{-1}y_2\bigr)\cdot
 \frac{2\cdot dy_1}{N}
-\psi\bigl(2N^{-1}y_1\bigr)\cdot \psi'\bigl(2N^{-1}y_2\bigr)\cdot
 \frac{2\cdot dy_2}{N}.
\]
Hence we have the following formula:
\begin{multline}
 d\chi_N\wedge\bigl\langle\delbar\theta,\theta-\theta^{\dagger}\bigr\rangle
=-\psi'\bigl(2N^{-1}y_1\bigr)\cdot\psi'\bigl(2N^{-1}y_2\bigr)
\cdot
 \frac{2\cdot dy_1}{N}
\wedge
 \Bigl\langle \delbar\theta,\,\,
 \frac{\sqrt{-1}}{2}h^{-1}\del_{y_1}h\cdot dx_1
 +\theta_2-\theta_2^{\dagger}
 \Bigr\rangle \\
-\psi'\bigl(2N^{-1}y_1\bigr)\cdot\psi\bigl(2N^{-1}y_2\bigr)
 \frac{2\cdot dy_2}{N}
\wedge
 \Bigl\langle \delbar\theta,\,\,
 \theta_1-\theta_1^{\dagger}
+\frac{\sqrt{-1}}{2}h^{-1}\del_{y_2}h\cdot dx_2
 \Bigr\rangle.
\end{multline}
We have the following estimate independently of $N$:
\[
 \psi'\bigl(2N^{-1}y_1\bigr)\cdot\psi\bigl(2N^{-1}y_2\bigr)
\cdot\frac{2\cdot dy_1}{N}
\wedge
 \Bigl\langle\delbar\theta,\,\,
 \frac{\sqrt{-1}}{2}h^{-1}\del_{y_1}h\cdot dx_1
 \Bigr\rangle
=O\Bigl(
 \bigl|\delbar\theta\bigr|\cdot
 \bigl|\del_{y_1}\Psi_h\bigr|\cdot (2y_1+A)\cdot \dvol_g
 \Bigr).
\]
Note $\delbar\theta$ and $\bigl|\del_{y_1}\Psi_h\bigr|\cdot (2y_1+A)$
are $L^2$ with respect to the measure $\dvol_g$
(Lemma \ref{lem;04.1.29.21} and Lemma \ref{lem;04.1.11.6}).

We have the following estimate, independently of $N$:
\begin{equation}\label{eq;04.1.12.16}
-\psi'\bigl(2N^{-1}y_1\bigr)\cdot
 \psi\bigl(2N^{-1}y_2\bigr)\cdot \frac{2}{N}dy_1
\wedge \bigl\langle \delbar\theta,\theta_2-\theta_2^{\dagger}
 \bigr\rangle
=O\Bigl(
 \psi'\bigl(2N^{-1}y_1\bigr)
\cdot
 \psi\bigl(2N^{-1}y_2\bigr)
\cdot\bigl|\delbar\theta\bigr|\cdot
 \bigl|\theta_2-\theta_2^{\dagger}\bigr|\cdot\dvol_g
 \Bigr).
\end{equation}
The support of $\psi'\bigl(2N^{-1}y_1\bigr)$
is contained in
$\bigl\{2^{-1}N\leq y_1\leq 3^{-1}N\bigr\}$.
Hence (\ref{eq;04.1.12.16}) is dominated by the following:
\begin{multline} \label{eq;04.1.12.20}
 \int dx_1\cdot dx_2
 \int_{N/4}^{N/3}
 \frac{dy_1}{(2y_1+A)^2}
 \int_{0}^{N/3} \frac{dy_2}{(2y_2+A)^2}
\cdot\bigl|\delbar\theta\bigr|\cdot
 \bigl|\theta_2-\theta_2^{\dagger}\bigr|\\
\leq
 \left(
 \int dx_1\cdot dx_2
 \int_{N/4}^{N/3}\frac{dy_1}{(2y_1+A)^2}
 \int_{0}^{N/3}\frac{dy_2}{(2y_2+A)^2}
 \cdot\bigl|\delbar\theta\bigr|^2
 \right)^{1/2}\\
\times
 \left(
 \int dx_1\cdot dx_2
 \int_{N/4}^{N/3} \frac{dy_1}{(2y_1+A)^2}
 \int_{0}^{N/3}\frac{dy_2}{(2y_2+A)^2}\cdot
 \bigl|\theta_2-\theta_2^{\dagger}\bigr|^2
 \right)^{1/2}.
\end{multline}
Due to the integrability of $\bigl|\delbar\theta\bigr|^2$,
the first term in the right hand side goes to $0$
when $N\to\infty$.
The square of the second term is dominated by the following,
due to (\ref{eq;04.1.29.150}),
Lemma \ref{lem;04.1.29.21}
and
Lemma \ref{lem;04.1.29.100}:
\[
 \int_{N/4}^{N/3} \frac{dy_1}{(2y_1+A)^2}
 \int_0^{N/3}\frac{dy_2}{(2y_2+A)^2}
 \left(
 \rho_2^2\cdot (2y_2+A)^2+J
 \right).
\]
Here $J$ denotes an integrable function
with respect to the measure
$(2y_1+A)^{-2}\cdot(2y_2+A)^{-2}\cdot dy_1\cdot dy_2$.
The contributions of $J$ go to $0$ when $N\to\infty$,
due to the integrability of $J$.
On the other hand, 
$\int_{N/4}^{N/3} (2y_1+A)^{-2}\cdot dy_1\times
 \int_0^{N/3} \rho_2^2\cdot dy_2$
are bounded independently of $N$.
Thus the second term in the right hand side of (\ref{eq;04.1.12.20})
is bounded.
Hence the right hand side of (\ref{eq;04.1.12.20}) converges to $0$
when $N\to\infty$.
Thus we obtain the following convergence:
\[
 \lim_{N\to\infty}
 \int_{M_P}
 -\psi'\bigl(2N^{-1}y_1\bigr)\cdot\psi'\bigl(2N^{-1}y_2\bigr)
 \frac{2\cdot dy_1}{N}
\wedge
 \Bigl\langle \delbar\theta,\,\,
 \frac{\sqrt{-1}}{2}h^{-1}\del_{y_1}h\cdot dx_1+\theta_2-\theta_2^{\dagger}
 \Bigr\rangle
=0.
\]
Similarly, we obtain the following convergence:
\[
 \lim_{N\to\infty}
\int_{M_P}
 -\psi'\bigl(2N^{-1}y_1\bigr)\cdot\psi\bigl(2N^{-1}y_2\bigr)
 \frac{2\cdot dy_2}{N}
\wedge
 \Bigl\langle \delbar\theta,\,\,
 \theta_1-\theta_1^{\dagger}
+\frac{\sqrt{-1}}{2}h^{-1}\del_{y_2}h\cdot dx_2
 \Bigr\rangle=0.
\]
Thus we obtain the following convergence:
\[
\lim_{N\to\infty}
\int_{M_P}
 d\chi_N\wedge\bigl\langle\delbar\theta,
 \theta-\theta^{\dagger}\bigr\rangle
=0.
\]

Hence the proof of Proposition \ref{prop;04.1.12.30}
is accomplished.
\hfill\qed

\subsection{Tameness and pure imaginary property}

\label{subsection;04.1.30.411}
%c17.21

Let $P$ be a point of $D_i\cap D_j$.
For simplicity, we consider the case $i=1$, $j=2$.
Recall the integrability of $J_{102}$ on $M_{P\,1}$
with respect to the measure
$(y_2+A)^{-2}dy_1\cdot dx_2\cdot dy_2$
in the subsubsection \ref{subsubsection;04.1.23.1}.
(see the page \pageref{eq;04.1.29.200}).
For any point $Q\in \Deltabarast$,
we put
$M_{P\,1\,Q}:=\bigl\{(y_1,Q)\,\big|\,y_1\in\real_{\geq 0}\bigr\}
\subset M_{P\,1}$.
Then the restrictions of $J_{102}$ to $M_{P\,1\,Q}$
are integrable with respect to the measure $dy_1$
for almost all $Q\in\Deltabarast$
by the theorem of Fubini.

On the other hand, 
we obtain the integrability of the restriction
of $\bigl|\del_{y_1}\Psi_h\bigr|^2\cdot \bigl|\del_{y_1}\bigr|^{-2}$
to $\Deltabarast\times Q$
with respect to the measure
$(2y_1+A)^{-2}\cdot dx_1\cdot dy_1$
for almost every $Q\in \Deltabarast$,
from Lemma \ref{lem;04.1.29.21}.
Thus we obtain the following lemma.
\begin{lem}\label{lem;04.1.29.250}
For almost every $Q\in\Deltabarast$,
there exists an integrable function
$J_{Q}$ on $\Deltabarast$ with respect to
$\dvol_Q=(2y_1+A)^{-2}\cdot dx_1\cdot dy_1$,
such that the following holds:
\[
 \int_{T(R)\times Q}
 \bigl|\theta_1\bigr|^2\cdot\dvol_Q
\leq
 \int_{T(R)\times Q}
 \left(
 \frac{\rho_1^2}{4\pi^2}\cdot(2y_1+A)^2
+J_Q
 \right)\cdot\dvol_Q.
\]
Here we put
$T(R):=\bigl\{z\in\cnum\,\big|\,0\leq -\log|z|\leq R\bigr\}$.
\hfill\qed
\end{lem}

\begin{cor}\label{cor;04.1.29.301}
For almost every $Q\in\Deltabarast$,
the restriction $h_{|\Deltabarast\times Q}$
is tame and pure imaginary.
\end{cor}
\pf
It follows from Lemma \ref{lem;04.1.29.250}
and Proposition \ref{prop;04.1.25.10}.
\hfill\qed

\vspace{.1in}

Let $W$ be a compact subregion
contained in $U_P\cap D_1$,
and we put $Y:=\Deltabar\times W$.
\begin{lem}
The restriction $h_{|Y\setminus D}$
is tame and pure imaginary.
\end{lem}
\pf
On $Y\setminus D$, we describe $\theta$ as follows:
\[
 \theta=f\cdot \frac{dz_2}{z_2}+g\cdot dz_1.
\]
Then we can easily show $\det(t-f)$ is holomorphic
on $Y$, due to Lemma \ref{lem;04.1.21.1} and
Corollary \ref{cor;04.1.29.301}.
Since the roots of the polynomial $\det(t-f)_{|Q}$ are pure imaginary
for almost every $Q\in W$ due to Corollary \ref{cor;04.1.29.301},
the roots of the polynomial $\det(t-f)_{|Q}$ are pure imaginary
for every $Q\in W$.

Let us consider $\det(t-g)$.
We have $|g|^2=O\bigl(\bigl|\del_{z_1}\Psi_h\bigr|^2\bigr)$.
By using the maximum principle for the family
of tame pure imaginary harmonic bundles
(Lemma \ref{lem;a12.29.5}),
we obtain the boundedness of $|g|$.
Then we obtain the boundedness of $\det(t-g)$,
which implies that $\det(t-g)$ is holomorphic.
Thus we are done.
\hfill\qed

\begin{lem}
The restriction $h_{|U_P}$ is tame
and pure imaginary
\end{lem}
\pf
Let us describe $\theta$ as follows, on $U_P$:
\[
 \theta=f_1\cdot\frac{dz_1}{z_1}+f_2\cdot\frac{dz_2}{z_2}.
\]
Let us consider $\det(t-f_i)$.
By using the previous consideration,
we have already known that
$\det(t-f_i)$ are holomorphic on
$U_P-\{P\}$.
Hence we obtain that they are holomorphic on $U_P$.
\hfill\qed

\begin{thm}\label{thm;04.1.29.305}
The pluri-harmonic metric $h$ of $(E,\nabla)$
is tame and pure imaginary.
\end{thm}
\pf
Let $Q$ be any  point of $D_i-\bigcup_{j\neq i}D_i\cap D_j$.
Let $W_0$ be a coordinate neighbourhood of $Q$,
and $(z_1,z_2)$ be a coordinate of $W_0$
such that $z_2^{-1}(0)=W_0\cap D_i$.
Then we can take a sequence of coordinate neighbourhoods
$W_0,W_1\ldots,W_l$ such that
$W_a\cap W_{a+1}\cap D_i\neq \emptyset$
and $W_l\subset U_P$ for some $P$.
On each $W_a$,
we develop $\theta=f^{(a)}\cdot dz_2^{(a)}/z_2^{(a)}+g^{(a)}\cdot dz_1^{(a)}$.
By an inductive argument using Lemma \ref{lem;04.1.21.2},
we can show that
$\det(t-f^{(a)})$ and $\det(t-g^{(a)})$ are holomorphic.
Hence we obtain that 
$(E,\nabla,h)$ is tame and pure imaginary.
\hfill\qed

\begin{rem}
For the given proof of Theorem {\rm\ref{thm;04.1.29.305}},
the sets $D_i\cap \bigcup_{j\neq i}D_j$ have to contain some point.
It is easy to modify the argument
in the case $D_i\cap \bigcup_{j\neq i}D_j=\emptyset$
for some component $D_i$.
For example, we have only to add some extra smooth divisor 
which intersects $D_i$.
Or, we have only to take a point $P$ 
in $D_i$, and we take a good coordinate around $P$
(see the subsubsection {\rm\ref{subsubsection;04.1.29.350}}).
\hfill\qed
\end{rem}

\subsection{The existence of pluri-harmonic metric
 for the higher dimensional projective case}

\label{subsection;04.1.30.412}

%c19.tex

Let $X$ be a smooth projective variety over $\cnum$,
and $D$ be a normal crossing divisor.
Let $(E,\nabla)$ be a flat simple bundle on $X-D$.
\begin{thm}\label{thm;04.1.23.20}
There exists a tame pure imaginary pluri-harmonic metric $h$
of $(E,\nabla)$,
which is unique up to positive constant multiplication.
\end{thm}
\pf
We use an induction on $\dim(X)$.

Let us take a sufficiently ample bundle $L$ of $X$.
We have the vector space $H^0(X,L)$.
We have the subspace $V_P:=\{f\in H^0(X,L)\,|\,f(P)=0\}$.
We put as follows:
\[
 \proj:=\proj(H^0(X,L)^{\lor}),
\quad
 \proj(P):=\proj(V_P^{\lor}).
\]
For any element $s\in \proj$,
we put $Y_s:=s^{-1}(0)$.
Let $U_0$ denote the subset of $\proj$,
which consists of the elements $s$
satisfying that  $Y_s$ are smooth
and that $Y_s\cap D$ are normal crossing.
It is Zariski dense subset of $\proj$.
We put $U(P)=U_0\cap \proj(P)$,
which is Zariski dense subset of $\proj(P)$.

For any element $s\in U_0$,
let $U(s)$ denote the subset of $U_0$,
which consists of the elements $s'$
satisfying that $Y_{s'}$ are transversal with $Y_s$
and that $D\cap Y_s\cap Y_{s'}$ are normal crossing.
We put $U_0^{(1)}:=\bigl\{(s,s')\,|\,s\in U_0,\,s'\in U(s)\bigr\}$.

For any point $P\in X$ and any element $s\in U_0$,
we put $U(s,P):=\bigl\{s'\in U_s\,|\,P\in Y_{s'}\bigr\}$.

Let $s$ be any element of $U_0$.
Due to the hypothesis of the induction,
we can take a tame pure imaginary pluri-harmonic metric $h_s$
of $(E,\nabla)_{|Y_s\setminus D}$.

\begin{lem}
Let $(s,s')$ be an element of $U^{(1)}_0$.
We take tame pure imaginary pluri-harmonic metrics $h_{s}$
and $h_{s'}$ of $(E,\nabla)_{|Y_s\setminus D}$ and
$(E,\nabla)_{|Y_{s'}\setminus D}$ respectively.
Then there exists a positive constant $a$ such that
$h_{s\,|\,Y_{s}\cap Y_{s'}\setminus D}
=a\cdot h_{s'\,|\,Y_s\cap Y_{s'}\setminus D}$.
\end{lem}
\pf
It follows from the uniqueness up to positive constant multiplication
(Proposition \ref{prop;04.1.3.4}).
Remark $\dim(Y_s\cap Y_{s'})\geq 1$.
\hfill\qed

\vspace{.1in}
Let us fix an element $s\in U_0$
and a tame pure imaginary pluri-harmonic metric $h_{s}$.

\begin{lem}\label{lem;04.2.2.10}
Let $(s_1,s_2)$ be an element of $U_0^{(1)}$
such that $s_i\in U(s)$.
Let us take $h_{s_i}$ such as
$h_{s_i\,|\,Y_{s}\cap Y_{s_i}\setminus D}
=h_{s\,|\,Y_{s}\cap Y_{s_i}\setminus D}$.
Assume $Y_s\cap Y_{s_1}\cap Y_{s_2}\setminus D\neq\emptyset$.
Then we have 
$h_{s_1\,|\,Y_{s_1}\cap Y_{s_2}\setminus D}
=h_{s_2\,|\,Y_{s_1}\cap Y_{s_2}\setminus D}$.
\end{lem}
\pf
There exists a positive constant $a$ such that
$h_{s_1\,|\,Y_{s_1}\cap Y_{s_2}\setminus D}
=a\cdot h_{s_2\,|\,Y_{s_1}\cap Y_{s_2}\setminus D}$.
On $Y_{s}\cap Y_{s_1}\cap Y_{s_2}\setminus D$,
we have
$h_{s_i\,|\,Y_{s_1}\cap Y_{s_2}\cap Y_s\setminus D}
=h_{s\,|\,Y_{s_1}\cap Y_{s_2}\cap Y_s\setminus D}$.
Hence we obtain $a=1$.
Thus we are done.
\hfill\qed

\vspace{.1in}

\begin{lem} \label{lem;04.1.3.10}
Let $Q$ be any point $X-D$.
For any elements $s_i\in U(s,Q)$ $(i=1,2)$,
there exists an elements $s_3\in U(s,Q)$ such that
$(s_i,s_3)\in U_0^{(1)}$ and 
$Y_s\cap Y_{s_i}\cap Y_3\setminus D\neq \emptyset$.
Moreover the set of such $s_3$ is Zariski dense in $U(s,Q)$.
\end{lem}
\pf
It follows from an easy argument using Zariski density.
\hfill\qed

\vspace{.1in}

Let $s_1$ be a section of $U(s,Q)$.
Let us take $h_{s_1}$ such as
$h_{s_1\,|\,Y_s\cap Y_{s_1}\setminus D}
=h_{s\,|\,Y_s\cap Y_{s_1}\setminus D}$.
We put $h_{Q}:=h_{s_1\,|\,Q}$.

\begin{lem}
We have the $C^{\infty}$-hermitian metric $h$
of $(E,\nabla)$ such that the following holds:
\[
h_{|Q}=h_Q,
\quad
 h_{|Y_s}=h_{s},
\quad
 h_{|Y_{s_1}}=h_{s_1},\,\,(s_1\in U(s,Q)).
\]
\end{lem}
\pf
It follows from Lemma \ref{lem;04.2.2.10} and
Lemma \ref{lem;04.1.3.10}.
\hfill\qed

\begin{lem}
The metric $h$ is pluri-harmonic metric.
\end{lem}
\pf
We denote the corresponding $(1,0)$-form by $\theta$
(the subsubsection \ref{subsubsection;04.1.5.10}).
We have only to show that
$\delbar\theta=\theta^2=0$. 
Let $Q$ be any point of $X-D$ and $H$
be a $\cnum$-subspace of $T_QX$ of codimension one.
We can take $s_1\in U(s,Q)$ such that $T_QY_{s_1}=H$.
Since the restriction of $h$ to $Y_{s_1}$ is pluri-harmonic,
we have $\delbar\theta_{|H}=\theta^2_{|H}=0$.
Then we obtain $\delbar\theta=\theta^2=0$.
\hfill\qed

\begin{lem}\label{lem;04.1.23.10}
$h$ is tame and pure imaginary.
\end{lem}
\pf
Once we know the tameness,
we obtain the pure imaginary property by considering
the restriction of $h$ to any $Y_s$.
Let us show the tameness.

Let $Q$ be a smooth point of $D$.
Let us take a neighbourhood $U$ of $Q$
with a coordinate $(z_1,\ldots,z_n)$
such that $U\cap D=z_1^{-1}(0)$.
We have the description:
\[
 \theta=f_1\cdot\frac{dz_1}{z_1}+\sum_{j=2}^n g_j\cdot dz_j.
\]
Let us see that the coefficients of
the characteristic polynomials $\det(t-f_1)$ and $\det(t-g_j)$
are holomorphic on $U$.

We put $S_i=\{z_i=0\}$.
We have the naturally defined projection
$\pi_2:U\lrarr S_1\cap S_2$.
For any point $P\in S_1\cap S_2$,
let us consider the restriction to $\pi_2^{-1}(P)$.
Here we may assume that $\pi_2^{-1}(P)$
is an intersection $Y_s\cap U$ for some $s$.
Then we have already known that
$\det(t-f_1)_{|\pi_2^{-1}(P)}$
and $\det(t-g_2)_{|\pi_2^{-1}(P)}$ are holomorphic
on $\pi_2^{-1}(P)$.
Then it is easy to derive that
$\det(t-f_1)$ and $\det(t-g_2)$ are holomorphic on $U$.
Similarly, we can derive that
$\det(t-g_i)$  are holomorphic on $U$.

Let $Q$ be any point of $D$.
Let $U$ be a neighbourhood of $Q$
with coordinate $(z_1,\ldots,z_n)$
such that $D\cap U=\bigcup_{i=1}^l\{z_i=0\}$.
We have the description:
\[
 \theta=\sum _{j=1}^lf_j\cdot\frac{dz_j}{z_j}+\sum_{j=l+1}^ng_j\cdot dz_j.
\]
By applying the consideration above,
we have already known that
$\det(t-f_j)$ and $\det(t-g_j)$ are holomorphic
outside the subset of codimension two.
Then we obtain that they are holomorphic on $U$
due to Hartogs' theorem.
Thus we obtain Lemma \ref{lem;04.1.23.10},
and hence the proof of Theorem \ref{thm;04.1.23.20}
is accomplished.
\hfill\qed

\section{An application}
\label{section;04.1.30.230}

\subsection{Preliminary (pull back of tame harmonic bundle)}
%c20.tex

Let $X$ and $Y$ be smooth projective varieties over $\cnum$.
Let $D_X$ and $D_Y$ be normal crossing divisors of $X$ and $Y$
respectively.
Let $F:X\lrarr Y$ be a morphism such that
$F^{-1}(D_Y)\subset D_X$.
Recall that we have the natural morphism
$F^{\ast}\Omega_Y^{1,0}(\log D_Y)\lrarr 
 \Omega_X^{1,0}(\log D_X)$.

\begin{lem}
Let $\harmonicbundle$ be a tame harmonic bundle on $Y-D_Y$.
Then the pull back $F^{\ast}\harmonicbundle$ is
a tame harmonic bundle on $X-D$.
If $\harmonicbundle$ is pure imaginary,
then $F^{\ast}\harmonicbundle$ is also pure imaginary.
\end{lem}
\pf
We take a prolongment $(\tilde{E},\tilde{\theta})$
of $(E,\theta)$ (Lemma \ref{lem;04.1.23.30}),
and then the eigenvalues of the residues of $\theta$
are pure imaginary, by definition.
Then we obtain the section
$F^{\ast}(\tilde{\theta})\in
 End\bigl(F^{\ast}\tilde{E}\bigr)\otimes F^{\ast}\Omega^{1,0}_Y(\log D_Y)$.
It naturally induces the regular Higgs field
$\theta_1\in \End\bigl(F^{\ast}\tilde{E}\bigr)
 \otimes\Omega^{1,0}_X(\log D_X)$.
The restriction of
$(F^{\ast}\tilde{E},\theta_1)$
to $X-D_X$ obviously coincides with $F^{\ast}_{|X-D_X}(E,\theta)$.
It implies that
$F^{\ast}_{|X-D_X}(E,\theta)$ is tame.

Assume that $\harmonicbundle$ is pure imaginary.
We have the irreducible decompositions
$D_X=\bigcup_i D_{X,i}$ and $D_Y=\bigcup_i D_{Y,i}$.
Let $P$ be a point of $D_{X,k}-\bigcup_{j\neq k}D_{X,k}\cap D_{X,j}$.
Let $D_{Y,i_1},\ldots,D_{Y,i_l}$ be the irreducible components of
$D_{Y}$, which contain $F(P)$.
Then the residue $\Res_{D_{X,k}}(\theta_1)_{|P}$ 
can be described as the linear combination of
$\Res_{D_{Y,i_j}}(\tilde{\theta})_{|F(P)}$ $(j=1,\ldots,l)$
with positive integer coefficients.
We also note that
$\Res_{D_{Y,i_j}}(\tilde{\theta})_{|F(P)}$ $(j=1,\ldots,l)$
are commutative, and that their eigenvalues are pure imaginary.
Thus we obtain that
the eigenvalues of $\Res_{D_{X,k}}(\theta_1)_{|P}$ 
are pure imaginary.
Thus we can conclude that
$F_{|X-D_X}^{\ast}\harmonicbundle$ is also pure imaginary.
\hfill\qed
\subsection{Pull back of semisimple local system}
%c20.1.tex

The following theorem is the answer to a question posed by Kashiwara.
\begin{thm}\label{thm;04.1.30.100}
Let $X$ and $Y$ be irreducible quasi projective varieties over $\cnum$.
Let $L$ be a semisimple local system on $Y$.
Let $F:X\lrarr Y$ be a morphism.
Then $F^{-1}(L)$ is also semisimple.
\end{thm}
\pf
We may assume that $X$ and $Y$ are smooth.
By a standard argument,
we can take smooth projective varieties $\overline{X}$
and $\overline{Y}$ such that the following holds:
\begin{itemize}
\item
 We have the inclusions $X\subset\overline{X}$
 and $Y\subset\overline{Y}$.
 The complements $D_{\overline{X}}:=\overline{X}-X$
 and $D_{\overline{Y}}:=\overline{Y}-Y$ are normal crossing divisors.
\item
 We have the morphism $\overline{F}:\overline{X}\lrarr\overline{Y}$
 such that $\overline{F}(X)\subset Y$ and
 $\overline{F}_X=F$.
 Note that we have
 $\overline{F}^{-1}(D_{\overline{Y}})\subset D_{\overline{X}}$.
\end{itemize}
Let $(E,\nabla)$ be a flat bundle on $Y$ corresponding to $L$.
Then we can take a tame pure imaginary pluri-harmonic metric $h$
of $(E,\nabla)$.
Then $F^{-1}(E,\nabla,h)$ is a tame pure imaginary harmonic bundle
on $X=\overline{X}-D_{\overline{X}}$.
Hence $F^{-1}(E,\nabla)$ is semisimple.
Thus we are done.
\hfill\qed

\noindent
{\it Address\\
Department of Mathematics,
Osaka City University,
Sugimoto, Sumiyoshi-ku,
Osaka 558-8585, Japan.\\
takuro@sci.osaka-cu.ac.jp\\
}

\end{document}